\pdfoutput=1

\RequirePackage{fix-cm}

\documentclass[leqno, fleqn, centertags, 12pt]{article}

\usepackage{latexsym}
\usepackage{amsmath}
\usepackage{amssymb}
\usepackage{verbatim}

\usepackage[T1]{fontenc}

\usepackage[upright, widespace]{fourier}

\usepackage{bm}

\usepackage{fullpage}

\usepackage{tikz-cd}

\tikzcdset{row sep/boom/.initial=5em}
\tikzcdset{column sep/boom/.initial=5.1em}

\tikzcdset{row sep/boomm/.initial=5.6em}
\tikzcdset{column sep/boomm/.initial=4.25em}

\tikzcdset{row sep/sboom/.initial=5.6em}
\tikzcdset{column sep/sboom/.initial=5.712em}

\tikzcdset{column sep/sma/.initial=3em}
\tikzcdset{row sep/sma/.initial=4.25em}

\tikzcdset{row sep/tri/.initial=6em}
\tikzcdset{column sep/tri/.initial=2.4em}

\newskip\nineskipamount \nineskipamount=9pt plus 0pt minus 0pt
\newskip\zeroskipamount \zeroskipamount=0pt plus 0pt minus 0pt
\usepackage[noorphans,vskip=\nineskipamount]{quoting}

\newcommand{\dis}{\displaystyle}

\makeatletter
\renewcommand{\@makefntext}[1]{\vspace*{0.5ex}\parindent=0em
\hspace*{-0.4em}
\hbox to 0.4em{\hss\@makefnmark}\hspace*{0.4em}{#1}
}
\makeatother

\newcounter{mysectionnumber}
\newcommand{\mysection}[2]{\setcounter{footnote}{0}
\setcounter{myparnum}{0}
\refstepcounter{mysectionnumber}
\vspace{21pt}{\Large {\themysectionnumber.} {#1}}\label{#2}\vspace*{15pt}}

\newcommand{\myuppar}[1]{\vspace{\medskipamount}\textbf{#1}\hspace*{0.5em}}

\newcommand{\myit}[1]{\textbf{\textit{#1}}\hspace{0.0em}}

\newcounter{myparnum}[mysectionnumber]
\renewcommand{\themyparnum}{\arabic{mysectionnumber}.\arabic{myparnum}}
\newcommand{\mypar}[2]{\refstepcounter{myparnum}{\vspace{\medskipamount}\textbf{{\themyparnum. #1}\label{#2}}\hspace{0.5em}}}

\newcounter{mylemmanum}[myparnum]
\newcommand{\proof}{\vspace{\medskipamount}{\textbf{{\emph{Proof}.}}\hspace*{1em}}}

\newcommand{\eproof}{ $\blacksquare$}

\def\sss{\hspace{0.05em}\ }
\def\dss{\hspace{0.1em}\ }
\def\trs{\hspace{0.15em}\ }
\def\qss{\hspace{0.2em}\ }
\def\pss{\hspace{0.3em}\ }
\def\oss{\hspace{0.4em}\ }

\def\halfff{\hspace*{0.025em}}
\def\fff{\hspace*{0.05em}}
\def\dff{\hspace*{0.1em}}
\def\trf{\hspace*{0.15em}}
\def\qff{\hspace*{0.2em}}
\def\pff{\hspace*{0.3em}}
\def\off{\hspace*{0.4em}}

\newcommand{\nsp}{\hspace*{-0.1em}}
\newcommand{\nnsp}{\hspace*{-0.15em}}

\newcommand{\dnsp}{\hspace*{-0.2em}}

\renewcommand{\leq}{\leqslant}
\renewcommand{\geq}{\geqslant}

\newcommand{\zzz}{\mathbf{Z}}

\newcommand{\ccc}{\mathbf{C}}
\newcommand{\rrr}{\mathbf{R}}

\newcommand{\image}{\operatorname{Im}\trf}

\newcommand{\kernel}{\operatorname{Ker}\trf}

\newcommand{\id}{\operatorname{id}}

\newcommand{\fred}{^{\dff \operatorname{Fred}}}

\newcommand{\ai}{\operatorname{a-ind}}

\def\omult{{%
    \setbox0\hbox{${=\mathrel{\mkern-5mu}=}$}%
    \rlap{\hbox to \wd0{\hss$|\mathrel{\mkern-7.5mu}|$\hss}}\box0
    }}

\newcommand{\num}[1]{|\qff #1 \qff|}
\newcommand{\bnum}[1]{\bigl|\qff #1 \qff\bigr|}

\newcommand{\norm}[1]{\|\qff #1 \qff\|}

\newcommand{\sco}[1]{\langle\trf #1 \trf\rangle}
\newcommand{\bsco}[1]{\left\langle\trf #1 \trf\right\rangle}

\newcommand{\comp}{^{\dff \mathrm{comp}}}

\newcommand{\ttoo}{\hspace*{0.2em}\longrightarrow\hspace*{0.2em}}

\begin{document}

\setlength{\baselineskip}{12pt plus 0pt minus 0pt}
\setlength{\parskip}{12pt plus 0pt minus 0pt}
\setlength{\abovedisplayskip}{12pt plus 0pt minus 0pt}
\setlength{\belowdisplayskip}{12pt plus 0pt minus 0pt}

\newskip\smallskipamount \smallskipamount=3pt plus 0pt minus 0pt
\newskip\medskipamount   \medskipamount  =6pt plus 0pt minus 0pt
\newskip\bigskipamount   \bigskipamount =12pt plus 0pt minus 0pt

\author{Nikolai\qss V.\qss Ivanov}
\title{Boundary\qss triplets\qss and\qss the\qss index\qss of\qss families\qss of\pss
self-adjoint\qss elliptic\qss boundary\qss problems} 
\date{}

\footnotetext{\hspace*{-0.65em}\copyright\qss 
Nikolai\qss V.\qss Ivanov,\qss 2023.\oss 
}

\footnotetext{\hspace*{-0.65em}The author\dss is\dss grateful\dss to\dss M.\dss Prokhorova\sss 
for careful\sss reading of\dss this paper and\sss pointing out\sss some infelicities.}

\maketitle

\renewcommand{\baselinestretch}{1}
\selectfont

\vspace*{0ex}

\vspace*{5.4ex}

\myit{\hspace*{0em}\large Contents}  \vspace*{1ex} \vspace*{\bigskipamount}\\ 
\hbox to 0.8\textwidth{\myit{1.}\hspace*{0.5em} Introduction \hfil 1}\hspace*{0.5em} \vspace*{0.25ex}\\
\hbox to 0.8\textwidth{\myit{2.}\hspace*{0.5em} Boundary\dss triplets and self-adjoint\sss extensions\hfil 5}  \hspace*{0.5em} \vspace*{0.25ex}\\
\hbox to 0.8\textwidth{\myit{3.}\hspace*{0.5em} Gelfand\dss triples \hfil 9} \hspace*{0.5em} \vspace*{0.25ex}\\
\hbox to 0.8\textwidth{\myit{4.}\hspace*{0.5em} Abstract\dss boundary\dss problems \hfil 12} \hspace*{0.5em} \vspace*{0.25ex}\\
\hbox to 0.8\textwidth{\myit{5.}\hspace*{0.5em} Families of\dss abstract\dss boundary\dss problems \hfil 22} \hspace*{0.5em} \vspace*{0.25ex}\\  
\hbox to 0.8\textwidth{\myit{6.}\hspace*{0.5em} Differential\dss boundary\sss problems of\dss order one \hfil 28} \hspace*{0.5em} \vspace*{0.25ex}\\
\hbox to 0.8\textwidth{\myit{7.}\hspace*{0.5em} Dirac-like\sss boundary\sss problems \hfil 33}  \hspace*{0.5em} \vspace*{0.25ex}\\ 
\hbox to 0.8\textwidth{\myit{8.}\hspace*{0.5em} Comparing\sss two boundary conditions \hfil 38}  \hspace*{0.5em}
\vspace*{0.25ex}\\ 
\hbox to 0.8\textwidth{\myit{9.}\hspace*{0.5em} Rellich\dss example \hfil 43}  \hspace*{0.5em}  
\vspace*{1.5ex}\\
\hbox to 0.8\textwidth{\myit{References} \hfil 44}\hspace*{0.5em}  \hspace*{0.5em}  \vspace*{0.25ex}

\renewcommand{\baselinestretch}{1}
\selectfont

\vspace*{5.4ex}

%\newpage
\mysection{Introduction}{introduction}

\myuppar{The\dss Lagrange\trs identity\qss (Green\dss formula).}
Let\sss $X$\sss be a compact\sss manifold\sss with\sss the boundary\sss $Y\off =\off \partial\dff X$\nnsp.\oss
Let\sss $D$ be a differential\sss operator acting on section of\dss a\dss Hermitian\dss bundle $E$ over $X$\nnsp,\oss
and\sss $D'$\sss be\sss the operator\sss formally adjoint\sss to $D$\nnsp.\oss
The\sss theory of\dss boundary\sss problems for $D$ crucially depends on an\sss identity 
of\dss the form\vspace{1.5pt}
\begin{equation}
\label{lagrange-intro}
\quad
\sco{\dff D\fff u\fff,\qff v \dff}
\pff -\pff
\sco{\dff u\dff,\qff D'\fff v \dff}
\off =\off
\bsco{\dff \gamma_1\dff u\dff,\qff 
\gamma_0\dff v \dff}_{\dff \partial}
\pff -\pff
\bsco{\dff \gamma_0\dff u\dff,\qff 
\gamma_1\dff v \dff}_{\dff \partial}
\pff,
\end{equation}

\vspace{-12pt}\vspace{1.5pt}
where $u\fff,\qff v$ are sections of\dss $E$\nnsp.\oss
It\dss is\dss known as\sss the\qss \emph{Lagrange\trs identity}\pss or\dss the\qss \emph{Green\dss formula}.\oss
The scalar\sss products\sss 
$\sco{\dff \bullet\fff,\qff \bullet \dff}$\sss
and\sss $\sco{\dff \bullet\fff,\qff \bullet \dff}_{\dff \partial}$\sss
are obtained\dss by\sss taking\sss the scalar product\sss in $E$ and\sss
integrating over $X$ and\sss $Y$ respectively.\oss
The operators $\gamma_0\dff,\pff \gamma_1$ are\sss the\qss \emph{boundary operators}\pss
taking\sss sections of\dss $E$ over $X$\sss to sections of\dss the restriction\sss
$E\trf|\trf Y$ of\dss $E$\sss to $Y$\dnsp.\oss
The identity\qss (\ref{lagrange-intro})\qss is\dss established\sss first\sss
for smooth sections $u\fff,\qff v$\sss and\sss then extended\sss to sections\sss 
in\dss Sobolev\dss spaces.\oss

When $D$\sss is\dss a differential\sss operator of\dss order $2$\nnsp,\oss
it\dss is\dss only\sss natural\dss to take as\sss $\gamma_0$\sss the restriction of\dss sections\sss to $Y$
and as\sss $\gamma_1$\sss the normal\sss derivative along\sss $Y$\dnsp.\oss
M.I.\dss Vishik\qss \cite{v}\qss discovered\dss that\sss this\qss ``na\"{i}ve''\qss  choice of\dss
$\gamma_0\dff,\pff \gamma_1$\sss is\dss not\sss the most\sss efficient\sss one.\oss
Namely,\oss it\dss is\dss advantageous\sss to adjust\sss the operator $\gamma_1$ 
by\sss replacing\sss it\dss by $\gamma_1\qff -\qff P$\dnsp,\oss
where\sss $P$\dnsp,\oss nowadays\sss known as\sss the\qss
\emph{Dirichlet--to--Neumann}\pss operator{},\oss takes a section\sss $f$ over\sss $Y$\sss
to\sss $\gamma_1\dff u$\nnsp,\oss where\sss $u$\sss is\dss the solution of\dss boundary\sss problem\sss
$D\dff u\off =\off 0$\nnsp,\dss $\gamma_0\dff u\off =\off f$\nnsp.\oss
Then\sss the identity\qss (\ref{lagrange-intro})\qss still\dss holds,\oss
and,\oss moreover{},\oss is\dss valid\sss for $u\fff,\qff v$\sss belonging\sss
to\sss the maximal\dss domains of\dss definition of\dss $D\fff,\qff D'$ respectively.

G.\dss Grubb\qss \cite{g1}\qss extended\trs Vishik's\dss theory\sss to operators
of\dss even\sss order $2 m$\qss (and strengthened\dss it\sss in some\sss respects).\oss
Both\dss Vishik\dss and\dss Grubb\dss 
considered only operators acting\sss on\sss functions,\oss 
but\sss the generalization\sss to sections of\dss bundles\dss is\dss routine.\oss
At\sss the same\sss time\sss the assumption\sss of\sss the even order was used\dss throughout\sss
and was even\sss build\sss into\sss the notations.\oss
There are $2 m$\sss boundary operators
and\sss they are split\sss into\sss to\sss two groups
of\dss $m$\sss operators each,\oss 
one leading\sss to $\gamma_0$ and\sss the other\sss to $\gamma_1$\nsp.\oss
Much\sss later\dss
B.M.\dss Brown,\oss G.\dss Grubb\dss and\trs I.G.\dss Wood\qss \cite{bgw}\qss
indicated\dss that\trs Grubb's\trs theory can\sss be adapted\dss to 
some matrix operators of\dss order $1$\nnsp.\oss

\myuppar{Self-adjoint\sss operators and\dss boundary\sss triplets.}
In\sss the present\sss paper we are interested only\sss in\sss formally self-adjoint\sss operators,\oss
and\sss their realizations by self-adjoint\sss operators\sss in\dss Hilbert\dss spaces.\oss
The\sss theory of\dss boundary\sss triplets provides an abstract\sss
axiomatic framework for using\sss the\dss Lagrange\dss identity\dss
to construct\sss self-adjoint\sss operators in\dss Hilbert\dss spaces.\oss
Let\sss $H_{\dff 0}$\sss and\sss $K^{\dff \partial}$\sss be Hilbert spaces and $T$ be
a densely defined\sss symmetric operator\sss in $H_{\dff 0}$\nsp,\oss
$T^{\dff *}$ be its adjoint\sss operator,\qss
and\sss $\mathcal{D}\off =\off \mathcal{D}\dff(\trf T^{\dff *}\trf)$
be\sss its domain.\oss
Let\sss  
$\gamma_0\dff,\pff \gamma_1\dff \colon\dff \mathcal{D}\qff \ttoo\qff K^{\dff \partial}$\sss
be\sss two\sss linear maps.\oss
The\sss triple $(\trf \mathcal{D}\fff,\pff \gamma_0\dff,\pff \gamma_1\trf)$\sss is\dss
a\qss \emph{boundary\dss triplet}\pss for $T^{\dff *}$\sss if\dss
$\gamma_0\dff \oplus\dff \gamma_1\dff \colon\dff \mathcal{D}\qff \ttoo\qff 
K^{\dff \partial}\dff \oplus\dff K^{\dff \partial}$\sss
is\dss surjective and\vspace{3pt}
\[
\quad
\sco{\dff T^{\dff *} u\fff,\qff v \dff}
\pff -\pff
\sco{\dff u\dff,\qff T^{\dff *} v \dff}
\off =\off
\bsco{\dff \gamma_1\dff u\dff,\qff 
\gamma_0\dff v \dff}_{\dff \partial}
\pff -\pff
\bsco{\dff \gamma_0\dff u\dff,\qff 
\gamma_1\dff v \dff}_{\dff \partial}
\]

\vspace{-12pt}\vspace{3pt} 
for every\sss $u\fff,\qff v\qff \in\qff \mathcal{D}$\dnsp.\oss
Here\sss $\sco{\dff \bullet\fff,\qff \bullet \dff}$\sss
is\dss the scalar product\sss in $H_{\dff 0}$ and\sss
$\sco{\dff \bullet\fff,\qff \bullet \dff}_{\dff \partial}$\sss
is\dss the scalar product\sss in $K^{\dff \partial}$\dnsp.\oss
A boundary\sss triplet\sss for $T^{\dff *}$ exists if\dss and\dss only\trs if\dss
$T$ admits a self-adjoint\sss extension,\oss 
and\sss then\sss the self-adjoint\sss extensions of\dss $T$ are in a natural\sss one-to-one
correspondence with\qss \emph{self-adjoint\sss relations}\dss
$\mathcal{B}\qff \subset\qff K^{\dff \partial}\dff \oplus\dff K^{\dff \partial}$\dnsp,\oss
which should\dss be understood as boundary conditions for $T$ or $T^{\dff *}$\dnsp.\oss
Such self-adjoint\sss relations are a minor\sss but\sss essential\dss generalization\sss
of\dss self-adjoint\sss operators\sss $K^{\dff \partial}\qff \ttoo\qff K^{\dff \partial}$\dnsp.\oss 
Also,\oss if\dss a boundary\sss triplet\sss exists,\oss then\sss 
$\mathcal{D}\dff(\qff \overline{T}\qff)\off =\off \kernel\fff \gamma_0\dff \oplus\dff \gamma_1$\nsp.\oss
We refer\sss to\sss the book of\trs K.\dss Schm\"{u}dgen\qss \cite{s},\oss Chapter\qss 14,\oss
for\sss the details.\oss
An outline of\dss this\sss theory,\oss sufficient\sss for our\sss purposes,\oss
is\dss contained\sss in\qss \cite{i2},\oss Sections\qss 11\qss and\qss 12.\oss

Let\sss $D$\sss be a formally self-adjoint\sss differential\sss operator{}.\oss
In\sss this case\sss the\dss Lagrange\dss identity\qss (\ref{lagrange-intro})\qss leads\sss to a description
of\dss some self-adjoint\dss boundary conditions,\oss 
i.e.\oss of\dss relations between\sss $\gamma_0\dff u$\sss
and\sss $\gamma_1\dff u$\sss leading\sss to self-adjoint\sss realisations of\dss $D$\sss
in\sss a\dss Hilbert\dss space,\oss
usually\sss in\sss the\dss Sobolev\dss space\sss 
$H_{\dff 0}\dff(\trf X^{\dff \circ}\fff,\qff E\trf)$\nnsp,\oss
where\sss $X^{\dff \circ}\off =\off X\qff \smallsetminus\qff Y$\dnsp.\oss
In\sss fact,\oss the most\sss interesting\sss self-adjoint\dss 
boundary conditions are defined\sss in\sss terms of\pss
``na\"{i}ve''\qss boundary operators,\oss 
not\dss the ones adjusted according\sss to\dss Vishik and\dss Grubb.\oss
But\sss the\qss ``na\"{i}ve''\qss boundary operators do not\sss form
a boundary\sss triplet\sss and by\sss this reason are not\sss good enough.\oss
The results of\trs Vishik\qss \cite{v}\qss and\trs Grubb\qss \cite{g1}\qss concerned\sss with\sss adjusting\sss
the\qss ``na\"{i}ve''\qss boundary operators for $D$ can\sss be understood\dss
as a construction of\dss a boundary\sss triplet\sss for\sss $D^{\fff *}$
starting\sss with\sss the\qss ``na\"{i}ve''\qss Lagrange\dss identity\qss (\ref{lagrange-intro}).\oss
But\sss they\sss were proved\dss before\sss 
the notion of\dss boundary\sss triplets appeared.\oss

\myuppar{The index of\dss families of\dss self-adjoint\sss operators.}
The\sss theory of\dss boundary\sss triplets\sss turns out\sss be a very efficient\sss tool\dss
to study\dss the index of\dss families of\dss self-adjoint\sss operators.\oss
It\sss was already used\dss by\sss the author\sss in\qss \cite{i2}\qss to explain
why\sss the index\sss theorem of\pss \cite{i2}\qss applies only\sss to
bundle-like boundary conditions.\oss
See\qss \cite{i2},\oss Section\qss 13.\oss
The key result\sss for\sss the applications\sss to\sss the index\sss theory\dss
is\dss a corollary of\dss one of\dss the main\sss results 
of\dss the\sss theory of\dss boundary\sss triplets,\oss
namely,\oss of\dss the\qss \emph{Krein--Naimark\trs resolvent\dss formula}.\oss
See\sss the identities\qss (\ref{cayley-t-i})\qss and\qss (\ref{cayley-t-minus-i})\qss below.\oss
Originally\dss the\dss Krein--Naimark\trs resolvent\dss formula and\dss this corollary
are proved\dss for\sss the boundary\sss triplets constructed\sss abstractly\sss
in\sss terms of\dss the operator\sss theory.\oss
But\sss the boundary\sss triplets are unique in a very strong sense,\oss
and\dss this allows\sss to apply\sss this corollary\sss to families of\dss differential\sss operators\qss
\emph{after\sss the\trs Lagrange\dss identity\trs is\dss adjusted}.\oss

The\dss Krein--Naimark\trs resolvent\dss formula\dss is\dss concerned\sss with extensions
of\dss symmetric operators,\oss but\sss in all\sss sources known\sss to\sss the author\dss
it\dss is\dss proved without\sss referring\sss to\dss von\dss Neumann\dss theory of\dss extensions.\oss
For\sss the sake of\dss the readers who may be,\oss like\sss the author{},\oss
not\sss quite comfortable about\sss such state of\dss affairs,\oss
we provided\sss in\dss Section\qss \ref{abstract}\qss a direct\sss proof\dss
of\dss the identities\qss (\ref{cayley-t-i})\qss and\qss (\ref{cayley-t-minus-i})\qss
based on\dss von\dss Neumann\trs theory and\dss bypassing\sss the
Krein--Naimark\trs resolvent\dss formula.\oss
For a proof\trs based on\sss the\dss Krein--Naimark\trs resolvent\dss formula
see\qss \cite{i2},\oss Section\qss 12.\oss

\myuppar{An abstract\dss construction of\dss boundary\dss triplets.} 
In\dss Section\qss \ref{abstract-index}\qss we develop an abstract\sss axiomatic version
of\dss the adjustment\sss procedure of\trs Vishik\dss and\trs Grubb.\oss
The starting\sss point\dss is\dss an axiomatic version of\dss the\qss ``na\"{i}ve''\qss
Lagrange\dss identity,\oss closely\sss related\dss to\sss the
abstract\dss boundary\sss problems\sss framework of\pss \cite{i2},\oss Section\qss 5.\oss
Another\sss key\sss ingredient\dss is\dss a\qss \emph{reference operator},\oss
a self-adjoint\sss and\sss invertible extension $A$ of\dss $D$\dnsp,\oss
which\dss is\dss assumed\dss to be na\"{i}vely defined\dss by\sss the boundary\sss
condition\sss $\gamma_0\off =\off 0$\nnsp.\oss
Under appropriate assumptions we construct\sss a boundary\sss triplet\sss
generalizing\sss the boundary\sss triplet\sss implicitly\sss present\sss in\sss the results
of\trs Vishik\dss and\trs Grubb.\oss
An\sss important\sss role in\sss the construction of\trs
Grubb\qss \cite{g1}\qss is\dss played\dss by\sss two results of\qss
J.L.\dss Lions\dss and\trs E.\dss Magenes\qss \cite{lm2},\qss \cite{lm3}.\oss
We prove an abstract\sss version of\dss them.\oss
See\qss Theorems\qss \ref{extending-gamma}\qss and\qss \ref{extended-isomorphism}.\oss
These\sss two\sss theorems are proved\dss by an adaptation of\dss arguments of\qss
Lions\dss and\trs Magenes.\oss
Most\sss of\dss the other\sss proofs in\dss Section\qss \ref{abstract-index}\qss
are adapted\dss from\dss arguments of\trs Grubb\qss \cite{g1}.

Our axiomatic approach\sss involves\sss the notion of\pss 
\emph{Gelfand\trs triples}\pss in an essential\sss manner{}.\oss
We need\sss only\sss the relevant\sss definitions and 
a couple of\dss the most\sss basic properties,\oss
and we included\sss in\dss Section\qss \ref{gelfand-triples}\qss 
a self-contained exposition of\dss what\dss is\dss needed.\oss

\myuppar{Families of\dss abstract\dss boundary\sss problems.}
In\dss Section\qss \ref{families}\qss we add\sss parameters\sss to\sss
the\sss theory developed\dss in\dss Sections\qss \ref{abstract}\sss --\sss \ref{abstract-index}\qss
and apply\sss the parameterized\dss theory\sss to\sss the index.\oss
Let\sss $W$\sss be a reasonable\sss topological\sss space.\oss
Suppose\sss that\sss everything\sss in\sss sight\sss depends on\sss the parameter\sss $w\qff \in\qff W$\dnsp.\oss
We indicate\sss this dependence by a subscript.\oss
The index of\dss the family\sss $A_{\dff w}\dff,\pff w\qff \in\qff W$\sss of\dss
reference operators\dss is\dss equal\dss to $0$\sss because\sss these operators are invertible.\oss
By\sss our assumptions\sss these operators are defined\sss by\sss the boundary conditions
$\gamma_{w\trf 0}\off =\off 0$\nnsp.\oss
Let\sss $\mathcal{B}_{\fff w}\dff,\pff w\qff \in\qff W$\sss be a family of\dss
na\"{i}ve\dss boundary\sss conditions.\oss
Replacing\sss the boundary conditions\sss $\gamma_{w\trf 0}\off =\off 0$\sss
by\sss the boundary conditions\sss $\mathcal{B}_{\fff w}$ results in a family of\dss self-adjoint\sss
operators with\sss the index equal\dss to\sss the index of\dss
$\mathcal{B}_{\fff w}\qff -\qff M_{\dff w}\dff,\pff w\qff \in\qff W$\dnsp,\oss
where\sss $M_{\dff w}$ are\sss the adjusting operators,\oss
analogues of\dss the\dss Dirichlet--to--Neumann\dss ones.\oss
See\sss the\sss last\sss subsection of\trs Section\qss \ref{families}\qss for\sss
the precise statement.\oss
The appearance of\dss the summand\sss $-\qff M_{\dff w}$\sss
is\dss fairly surprising,\oss but\sss the operators\sss $\mathcal{B}_{\fff w}$
are usually\sss invertible and\dss then\sss the index of\dss the uncorrected\sss
family\sss $\mathcal{B}_{\fff w}\dff,\pff w\qff \in\qff W$\sss is\sss $0$\dnsp.\vspace{-0.125pt}

\myuppar{Differential\dss boundary\sss problems of\dss order one.}
Section\qss \ref{differential-boundary-problems}\qss is\dss devoted\dss
to\sss the application of\dss the abstract\sss theory\sss to
differential\sss operators of\dss order one.\oss
The key\sss role\dss is\dss played\dss by\sss the results related\dss to\sss the\qss
\emph{Calder\'{o}n\dss projector{}.}\oss
They\sss are used\sss to prove\sss that\sss 
the adjusting operators $M_{\dff w}$ are pseudo-differential\sss
operators of\dss order zero continuously depending on $w\qff \in\qff W$\dnsp.\oss
We identify\dss the symbols of\dss the adjusting operators $M_{\dff w}$
and\sss prove\sss that\sss the graph of\dss $M_{\dff w}$\sss is\dss
the space of\dss the\dss Cauchy\dss data of\dss solutions of\dss
the equation\sss $D^{\fff *}\dff u\off =\off 0$\nnsp.\oss

In\dss Section\qss \ref{dirac-like}\qss we use\sss the results of\trs 
Section\qss \ref{differential-boundary-problems}\qss to give an analytic
and\sss fairly\sss elementary\qss proof\dss of\dss 
the index\sss theorem\sss for\dss Dirac-like boundary problems
from\qss \cite{i2},\oss Section\qss 15.\oss
See\dss Theorem\qss \ref{index-theorem}.\oss
In\qss \cite{i2}\qss this\sss theorem\sss was deduced\dss from\sss
its analogue for\sss the\sss topological\dss index and\dss
the general\dss index\dss theorem\sss for operators of\dss order one.\oss
In a sense\sss the proof\dss from\qss \cite{i2}\qss provides\sss
topological\dss reasons for\sss the need\dss to adjust\sss
the boundary operators and\sss predicts\sss the symbols of\dss adjusting\sss operators.\oss
The desire\sss to find\sss a direct\sss analytic proof\dss of\dss this\sss
theorem was\sss the starting\sss point\sss of\dss the present\sss paper{}.\oss 

In\dss Section\qss \ref{comparing}\qss we modify\dss the arguments of\trs
Section\qss \ref{dirac-like}\qss in order\sss to prove an\dss
Agranovich--Dynin\dss type\sss theorem computing\sss the difference of\dss indices
of\dss two families of\dss self-adjoint\sss problems differing only\sss by\sss 
the boundary conditions.\oss
The formula for\sss the difference\dss is\dss more complicated\dss than\sss
the classical\dss one and\dss involves an appropriate form of\dss the adjusting\sss
operators $M_{\dff w}$\nsp.\oss

We\sss limited ourselves\sss by\sss the differential\sss operators of\dss order one
in\sss order\sss to avoid\dss technical\sss complications and\sss present\sss
the main\sss ideas in\sss the most\sss transparent\sss form.\oss
Most\sss of\dss the results can\sss be extended\dss to pseudo-differential\sss
operators satisfying\sss the\sss transmission condition.\oss\vspace{-0.125pt}

\myuppar{Rellich\dss example.}
In\dss Section\qss \ref{rellich}\qss we illustrate\sss our abstract\sss
theory by an example,\oss due\sss to\dss F.\dss Rellich,\oss 
of\dss a family of\dss self-adjoint\dss boundary problems
parameterized\dss by\sss the circle for\sss a fixed operator\sss of\dss order $2$\nnsp.\oss
The operator\dss is\sss $-\qff d^{\dff 2}/\dff d x^{\dff 2}$ 
on\sss the interval\sss $[\trf 0\fff,\qff 1\dff]$\nnsp.\oss
The index of\dss this family can be computed\dss by\sss brute force
and\dss is\dss equal\dss to $1$\nnsp.\oss
We prove\sss this result\sss as an application of\dss our abstract\sss theory.\oss
In\sss dimension one\sss there\dss is\dss no need\dss to adjust\sss
the\dss Lagrange\dss identity,\oss but\sss one can still\dss do\sss this.\oss
Each way\sss leads\sss to a proof\dss that\sss the index\dss is\dss equal\dss to $1$\nnsp.\oss

\newpage
\mysection{Boundary\qss triplets\qss and\qss self-adjoint\qss extensions}{abstract}

\myuppar{The main abstract\sss example.}
Let\sss us consider\sss the following abstract\sss situation\sss
taken\sss from\qss \cite{s},\oss Example\qss 14.5.\oss
Let $H$ be a separable\dss Hilbert\dss space,\pss
$T$ be a densely defined closed symmetric operator in $H$\nnsp,\oss
and $A$ be a fixed self-adjoint\sss extension of\dss $T$\dnsp.\oss
By\sss $\dotplus$ we will\sss denote\sss the direct,\oss 
but\sss not\sss necessarily orthogonal,\oss sum of\dss subspaces of\dss $H$\nnsp.\oss
Let\sss us fix a number\sss $\mu\qff \in\qff \ccc\qff \smallsetminus\qff \rrr$\nnsp.\oss 
Then\sss the numbers $\mu\dff,\off \overline{\mu}$\dss  belong\sss to\sss the resolvent\sss set\sss of\dss $A$\nnsp.\oss
In particular,\oss the operators\sss $(\trf A\qff -\qff \mu\trf)^{\dff -\dff 1}$\sss
and\sss
$(\trf A\qff -\qff \overline{\mu}\qff)^{\dff -\dff 1}$\sss
are well\sss defined.\oss
Let\vspace{1.5pt}
\[
\quad
\mathcal{K}_{\dff +}
\off =\off
\kernel\fff (\trf T^{\dff *}\qff -\qff\halfff \mu\qff)
\off =\off 
\image\fff (\trf T\qff -\qff \overline{\mu}\qff)^{\dff \perp}
\quad
\mbox{and}\quad
\]

\vspace{-37.5pt}
\[
\quad
\mathcal{K}_{\dff -}
\off =\off
\kernel\fff (\trf T^{\dff *}\qff -\qff \overline{\mu}\qff)
\off =\off 
\image\fff (\trf T\qff -\qff \mu\qff)^{\dff \perp}
\pff.
\]

\vspace{-12pt}\vspace{1.5pt}
Then\sss the domain\sss $\mathcal{D}\dff(\trf T^{\dff *}\trf)$\sss
of\dss $T^{\dff *}$\sss is\dss equal\sss to\vspace{3pt}
\begin{equation}
\label{st-example}
\quad
\mathcal{D}\dff(\trf T^{\dff *}\trf)
\off =\off
\mathcal{D}\dff(\trf T\trf)
\qff \dotplus\qff
A\trf (\trf A\qff -\qff \mu\trf)^{\dff -\dff 1}\dff \mathcal{K}_{\dff -}
\qff \dotplus\qff
(\trf A\qff -\qff \mu\trf)^{\dff -\dff 1}\dff \mathcal{K}_{\dff -}
\off.
\end{equation}

\vspace{-12pt}\vspace{3pt}
See\qss \cite{s},\oss Proposition\qss 14.11.\oss
Hence every\sss $z\qff \in\qff \mathcal{D}\dff(\trf T^{\dff *}\trf)$\sss
can be uniquely\sss written as\vspace{3pt}
\begin{equation}
\label{z-sum}
\quad
z
\off =\off
z_{\trf T}
\qff +\qff
A\trf (\trf A\qff -\qff \mu\trf)^{\dff -\dff 1}\dff z_{\dff 0}
\qff +\qff
(\trf A\qff -\qff \mu\trf)^{\dff -\dff 1}\dff z_{\dff 1}
\off
\end{equation}

\vspace{-12pt}\vspace{3pt}
with $z_{\trf T}\qff \in\qff \mathcal{D}\dff(\trf T\trf)$
and\sss $z_{\dff 0}\dff,\qff z_{\dff 1}\qff \in\qff \mathcal{K}_{\dff -}$\nsp.\oss
Then $K\off =\off \mathcal{K}_{\dff -}$ and\dss the maps\sss
$\Gamma_{i}\dff \colon\dff z\off \longmapsto\off z_{\dff i}$\nsp,\dss
$i\off =\off 0\fff,\qff 1$\sss form a boundary\dss triplet\sss for $T^{\dff *}$\dnsp.\oss
This\dss is\dss the boundary\sss triplet\sss from\qss \cite{s},\oss Example\qss 14.5.\oss

\myuppar{The main example from\sss the point\sss of\dss view of\trs
von\dss Neumann\dss theory.}
By\dss von\dss Neuman\dss theory\vspace{3pt}\vspace{-0.75pt}
\[
\quad
\mathcal{D}\dff(\trf T^{\dff *}\trf)
\off =\off
\mathcal{D}\dff(\trf T\trf)
\qff \dotplus \qff
\mathcal{K}_{\dff +}
\qff \dotplus \qff 
\mathcal{K}_{\dff -}
\pff
\]

\vspace{-12pt}\vspace{3pt}\vspace{-0.75pt}
and\sss $A$ defines an\sss isometry
$V\dff \colon \mathcal{K}_{\dff +}\qff \ttoo\qff \mathcal{K}_{\dff -}$
such\sss that\sss $A$\sss is\dss  
the restriction of\dss 
$T^{\dff *}$\sss to\vspace{3pt}
\[
\quad
\bigl\{\pff x\qff +\qff y\qff -\qff V y\pff \bigl|\pff 
x\qff \in\qff \mathcal{D}\dff(\trf T\trf)\fff,\off
y\qff \in\qff  \mathcal{K}_{\dff +}\pff \bigr\}
\off.
\]

\vspace{-12pt}\vspace{3pt}
\mypar{Lemma.}{unitary-op}
$\dis
V y
\off =\off
\frac{A\qff -\qff \mu}{\dff A\qff -\qff \overline{\mu}\dff}\qff y$\oss
\emph{for every\sss
$y\qff \in\qff \mathcal{K}_{\dff +}$\nsp.\oss}

\proof
If\dss $y\qff \in\qff \mathcal{K}_{\dff +}$\nsp,\oss
then\vspace{3pt}
\[
\quad
(\trf A\qff -\qff \overline{\mu}\qff)\qff
(\trf y\qff -\qff V y\qff)
\off =\off
\mu\qff y
\qff -\qff
\overline{\mu}\pff V y
\qff -\qff
\overline{\mu}\qff y
\qff +\qff
\overline{\mu}\pff V y
\off =\off
(\trf \mu\qff -\qff \overline{\mu}\qff)\qff y
\pff,
\]

\vspace{-33pt}
\[
\quad
y\qff -\qff V y
\off =\off
(\trf A\qff -\qff \overline{\mu}\qff)^{\dff -\dff 1}\dff
(\trf \mu\qff -\qff \overline{\mu}\qff)\qff y
\pff,
\]

\vspace{-12pt}\vspace{3pt}
and\dss hence\qss
$\dis
V y
\off =\off
y\qff -\qff
(\trf A\qff -\qff \overline{\mu}\qff)^{\dff -\dff 1}\dff
(\trf \mu\qff -\qff \overline{\mu}\qff)\qff y
\off =\off
\frac{A\qff -\qff \mu}{\dff A\qff -\qff \overline{\mu}\dff}\qff y$\nsp.\oss  \eproof

\mypar{Lemma.}{von-N}
\emph{Suppose\sss that\sss
$z
\off =\off 
z_{\trf T}\qff +\qff z_{\dff +}\qff +\qff z_{\dff -}$\sss
with\sss
$z_{\trf T}\qff \in\qff \mathcal{D}\dff(\trf T\trf)$\nnsp,\qss
$z_{\dff +}\qff \in\qff \mathcal{K}_{\dff +}$\nsp,\oss and\sss
$z_{\dff -}\qff \in\qff \mathcal{K}_{\dff -}$\nsp.\oss
Let\sss
$z_{\trf 0}
\off =\off 
z_{\dff -}\qff +\pff V\fff z_{\dff +}$
and\sss
$z_{\dff 1}
\off =\off 
-\qff \mu\trf z_{\dff -}\qff -\pff \overline{\mu}\pff V\fff z_{\dff +}$\nsp.\oss
Then\sss the equality\pss \textup{(\ref{z-sum})}\qss holds.\oss}

\proof
Lemma\qss \ref{unitary-op}\qss implies\sss that\vspace{1.5pt}
\[
\quad
z_{\trf 0}
\off =\off 
z_{\dff -}\qff +\pff V\fff z_{\dff +}
\off =\off
z_{\dff -}
\qff +\qff
\frac{A\qff -\qff \mu}{\dff A\qff -\qff \overline{\mu}\dff}\off z_{\dff +}
\quad
\mbox{and}
\]

\vspace{-30pt}
\[
\quad
z_{\dff 1}
\off =\off 
-\qff \mu\trf z_{\dff -}\qff -\pff \overline{\mu}\pff V\fff z_{\dff +}
\off =\off
-\qff \mu\trf z_{\dff -}
\qff -\pff
\overline{\mu}\off
\frac{A\qff -\qff \mu}{\dff A\qff -\qff \overline{\mu}\dff}\off z_{\dff +}
\pff.
\]

\vspace{-12pt}\vspace{3pt}
It\sss follows\sss that\vspace{3pt}
\[
\quad
A\trf (\trf A\qff -\qff \mu\trf)^{\dff -\dff 1}\dff z_{\dff 0}
\off =\off
A\trf (\trf A\qff -\qff \mu\trf)^{\dff -\dff 1}\dff z_{\dff -}
\qff +\qff
A\trf (\trf A\qff -\qff \overline{\mu}\qff)^{\dff -\dff 1}\dff z_{\dff +}
\pff,
\]

\vspace{-30pt}
\[
\quad
(\trf A\qff -\qff \mu\trf)^{\dff -\dff 1}\dff z_{\dff 1}
\off =\off
-\qff (\trf A\qff -\qff \mu\trf)^{\dff -\dff 1}\dff \mu\trf z_{\dff -}
\qff -\qff
\overline{\mu}\qff (\trf A\qff -\qff \overline{\mu}\qff)^{\dff -\dff 1}\dff z_{\dff +}
\pff,
\]

\vspace{-12pt}\vspace{1.5pt}
and\dss hence\vspace{1.5pt}
\[
\quad
A\trf (\trf A\qff -\qff \mu\trf)^{\dff -\dff 1}\dff z_{\dff 0}
\qff +\qff
(\trf A\qff -\qff \mu\trf)^{\dff -\dff 1}\dff z_{\dff 1}
\]

\vspace{-30pt}
\[
\quad
=\off
A\trf (\trf A\qff -\qff \mu\trf)^{\dff -\dff 1}\dff z_{\dff -}
\qff -\qff
(\trf A\qff -\qff \mu\trf)^{\dff -\dff 1}\dff \mu\trf z_{\dff -}
\off +\off
A\trf (\trf A\qff -\qff \overline{\mu}\qff)^{\dff -\dff 1}\dff z_{\dff +}
\qff -\qff
\overline{\mu}\qff (\trf A\qff -\qff \overline{\mu}\qff)^{\dff -\dff 1}\dff z_{\dff +}
\]

\vspace{-31.5pt}
\[
\quad
=\off
z_{\dff -}\qff +\qff z_{\dff +}
\pff.
\]

\vspace{-12pt}\vspace{3pt}
The equality\pss \textup{(\ref{z-sum})}\qss follows.\oss  \eproof

\myuppar{The\dss Lagrange\dss identity.}
In\sss the notations of\qss Lemma\qss \ref{von-N},\oss
let\sss\vspace{3pt}
\[
\quad
\Gamma_0\dff z
\off =\off 
z_{\dff -}\qff +\pff V\fff z_{\dff +}
\quad
\mbox{and}\quad
\Gamma_1\dff z
\off =\off 
-\qff \mu\trf z_{\dff -}\qff -\pff \overline{\mu}\pff V\fff z_{\dff +}
\pff.
\]

\vspace{-12pt}\vspace{3pt}
We claim\sss that\sss then\sss for every\sss
$x\fff,\qff y\qff \in\qff \mathcal{D}\dff(\trf T^{\dff *}\trf)$\vspace{3pt}
\begin{equation}
\label{l-von-N}
\quad
\sco{\dff T^{\dff *} x\fff,\qff y \dff}
\pff -\pff
\sco{\dff x\dff,\qff T^{\dff *} y \dff}
\off =\off
\bsco{\dff \Gamma_1\dff x\dff,\qff 
\Gamma_0\dff y \dff}
\pff -\pff
\bsco{\dff \Gamma_0\dff x\dff,\qff 
\Gamma_1\dff y \dff}
\pff,
\end{equation}

\vspace{-12pt}\vspace{3pt}
where all\sss scalar products are\sss taken\sss in\sss $H$\nnsp.\oss
In order\sss to prove\sss this,\oss let\sss us write $x$\sss and\sss $y$\sss
in\sss the form of\qss Lemma\qss \ref{von-N},\pss
$x
\off =\off 
x_{\trf T}\qff +\qff x_{\dff +}\qff +\qff x_{\dff -}$\sss
and\sss
$y
\off =\off 
y_{\trf T}\qff +\qff y_{\dff +}\qff +\qff y_{\dff -}$\nsp.\oss
Since\sss $T$\sss is\dss symmetric,\pss $x_{\trf T}$ and\sss $y_{\trf T}$
do not\sss affect\sss the validity of\dss the\dss Lagrange\dss identity and\dss hence
we can assume\sss that\sss $x_{\trf T}\off =\off y_{\trf T}\off =\off 0$\nnsp.\oss
Then\sss
$T^{\dff *}\dff x\off =\off \mu\dff x_{\dff +}\qff +\pff \overline{\mu}\qff x_{\dff -}$
and\dss
$T^{\dff *}\dff y\off =\off \mu\dff y_{\dff +}\qff +\pff \overline{\mu}\qff y_{\dff -}$\nsp.\oss
Therefore\sss the\sss left\dss hand side
of\pss (\ref{l-von-N})\qss is\dss equal\dss to\vspace{3pt}
\[
\quad
\mu\dff \sco{\dff x_{\dff +}\fff,\qff y_{\dff +} \dff}
\qff +\qff
\mu\dff \sco{\dff x_{\dff +}\fff,\qff y_{\dff -} \dff}
\qff +\qff
\overline{\mu}\qff \sco{\dff x_{\dff -}\fff,\qff y_{\dff +} \dff}
\qff +\qff
\overline{\mu}\qff \sco{\dff x_{\dff -}\fff,\qff y_{\dff -} \dff}
\]

\vspace{-33pt}
\[
\quad
-\off
\overline{\mu}\qff \sco{\dff x_{\dff +}\fff,\qff y_{\dff +} \dff}
\qff -\qff
\mu\dff \sco{\dff x_{\dff +}\fff,\qff y_{\dff -} \dff}
\qff -\qff
\overline{\mu}\qff  \sco{\dff x_{\dff -}\fff,\qff y_{\dff +} \dff}
\qff -\qff
\mu\dff \sco{\dff x_{\dff -}\fff,\qff y_{\dff -} \dff}
\]

\vspace{-33pt}
\[
\quad
=\off
(\trf \mu\qff -\pff \overline{\mu}\qff)\dff \sco{\dff x_{\dff +}\fff,\qff y_{\dff +} \dff}
\qff -\qff
(\trf \mu\qff -\pff \overline{\mu}\qff)\dff \sco{\dff x_{\dff -}\fff,\qff y_{\dff -} \dff}
\pff.
\]

\vspace{-12pt}\vspace{3pt}
The right\sss hand side of\pss (\ref{l-von-N})\qss is\dss equal\dss to\vspace{3pt}
\[
\quad
\bsco{\trf -\qff \mu\trf x_{\dff -}\qff -\pff \overline{\mu}\pff V\fff x_{\dff +}
\dff,\off
y_{\dff -}\qff +\pff V\fff y_{\dff +} \trf}
\pff -\pff
\bsco{\trf
x_{\dff -}\qff +\pff V\fff x_{\dff +}
\dff,\off
-\qff \mu\trf y_{\dff -}\qff -\pff \overline{\mu}\pff V\fff y_{\dff +} \trf}
\]

\vspace{-33pt}
\[
\quad
=\off
-\qff 
\mu\dff \sco{\dff x_{\dff -}\dff,\qff y_{\dff -}\dff }
\qff -\qff 
\mu\dff \sco{\dff x_{\dff -}\dff,\qff V y_{\dff +}\dff }
\qff -\qff 
\overline{\mu}\qff \sco{\dff V x_{\dff -}\dff,\qff y_{\dff -}\dff }
\qff -\qff 
\overline{\mu}\qff \sco{\dff V x_{\dff +}\dff,\qff V y_{\dff +}\dff }
\]

\vspace{-33pt}
\[
\quad
\phantom{=\off}
+\qff
\overline{\mu}\qff \sco{\dff x_{\dff -}\dff,\qff y_{\dff -}\dff}
\qff +\qff
\mu\dff \sco{\dff x_{\dff -}\dff,\qff V y_{\dff +}\dff}
\qff +\qff
\overline{\mu}\qff \sco{\dff V x_{\dff +}\dff,\qff y_{\dff -}\dff}
\qff +\qff
\mu\dff \sco{\dff V x_{\dff +}\dff,\qff V y_{\dff +}\dff}
\]

\vspace{-33pt}
\[
\quad
=\off
-\qff
(\trf \mu\qff -\pff \overline{\mu}\qff)\dff \sco{\dff x_{\dff -}\fff,\qff y_{\dff -} \dff}
\qff +\qff
(\trf \mu\qff -\pff \overline{\mu}\qff)\dff \sco{\dff V x_{\dff +}\fff,\qff V y_{\dff +} \dff}
\]

\vspace{-33pt}
\[
\quad
=\off
-\qff
(\trf \mu\qff -\pff \overline{\mu}\qff)\dff \sco{\dff x_{\dff -}\fff,\qff y_{\dff -} \dff}
\qff +\qff
(\trf \mu\qff -\pff \overline{\mu}\qff)\dff \sco{\dff x_{\dff +}\fff,\qff y_{\dff +} \dff}
\]

\vspace{-12pt}\vspace{3pt}
because\sss $V$\sss is\dss an\sss isometry.\oss
The\sss identity\qss (\ref{l-von-N})\qss follows.\oss

\myuppar{The boundary\sss triplet.}
Clearly,\oss the map\sss $z\off \longmapsto\off (\trf z_{\dff -}\dff,\qff z_{\dff +}\trf)$\sss
from $\mathcal{D}\dff(\trf T^{\dff *}\trf)$\sss to\sss
$\mathcal{K}_{\dff -}\dff \oplus\dff \mathcal{K}_{\dff +}$\sss
is\dss surjective.\oss
Since\sss $\mu\qff \not\in\qff \rrr$\nnsp,\oss
this implies\sss that\sss the map\sss\vspace{3pt}
\[
\quad
\Gamma_0\dff \oplus\dff \Gamma_1\dff \colon\dff
\mathcal{D}\dff(\trf T^{\dff *}\trf)
\qff \ttoo\qff
\mathcal{K}_{\dff -}\dff \oplus\dff \mathcal{K}_{\dff -}
\]

\vspace{-12pt}\vspace{3pt}
is\dss also surjective.\oss
In view of\pss (\ref{l-von-N})\qss
this implies\sss that\sss
$(\qff \mathcal{K}_{\dff -}\dff,\off \Gamma_0\dff,\off \Gamma_1 \trf)$\sss
is\dss a boundary\sss triplet\sss for\sss $T^{\dff *}$\dnsp.\oss
Lemma\qss \ref{von-N}\qss implies\sss that\sss this\dss is\dss the same boundary\sss
triplet\sss as\sss in\qss \cite{s},\oss Example\qss 14.5.\oss

\myuppar{Isometries\sss between subspaces.}
Let\sss $W\dff \colon\dff \mathcal{K}\qff \ttoo\qff \mathcal{K}\fff'$\sss
be an\sss isometry\sss between\sss two closed subspaces 
$\mathcal{K}\dnsp,\off \mathcal{K}\fff'\qff \subset\qff H$\nnsp.\oss
We will\sss denote by\sss $W_{\dff 0}$\sss the corresponding\qss
\emph{partial\dss isometry}\pss of\dss $H$\nnsp,\oss
i.e.\qss the operator equal\dss to $W$ on\sss $\mathcal{K}$ and\dss to $0$
on\sss $\mathcal{K}^{\dff \perp}$\dnsp.\oss
When\sss $\mathcal{K}\pff =\off \mathcal{K}\fff'$\dnsp,\oss
we will\sss denote by\sss $W_{\dff H}$\sss the operator\sss
$H\qff \ttoo\qff H$ equal\dss to $W$ on\sss $\mathcal{K}$ and\dss to\sss the identity
on\sss $\mathcal{K}^{\dff \perp}$\dnsp.\oss

\myuppar{The\dss Cayley\dss transforms.}
As above,\oss let\sss $\mu\qff \in\qff \ccc\qff \smallsetminus\qff \rrr$
and\sss $T$\sss be a symmetric operator in\sss $H$\nnsp.\oss
The\dss \emph{$\mu$\dnsp-Cayley\trs transform}\trs $U_{\dff \mu}\dff(\trf T\trf)$
of\dss $T$\sss is\dss the partial\dss isometry
equal\dss to\vspace{3pt}
\[
\quad
\frac{T\qff -\qff \mu}{\dff T\qff -\qff \overline{\mu}\dff}
\]

\vspace{-12pt}\vspace{3pt}
on\sss $\image\fff (\trf T\qff -\qff \overline{\mu}\qff)$
and as $0$ on\sss the orthogonal\sss complement\sss
$\image\fff (\trf T\qff -\qff \overline{\mu}\qff)^{\dff \perp}
\off =\off
\mathcal{K}_{\dff +}$\nnsp.\oss
It\dss induces an\sss isometry\sss
$\image\fff (\trf T\qff -\qff \overline{\mu}\qff)
\qff \ttoo\qff
\image\fff (\trf T\qff -\qff \mu\trf)$\nnsp.\oss
If\dss $A$\sss is\dss a self-adjoint\sss extension of\dss $T$ as above,\oss
then\sss $U_{\dff \mu}\dff(\trf A\trf)$ agrees with\sss $U_{\dff \mu}\dff(\trf T\trf)$
on\sss $\image\fff (\trf T\qff -\qff \overline{\mu}\qff)$\nnsp.\oss
In\sss this case\sss
$\image\fff (\trf A\qff -\qff \overline{\mu}\qff)\off =\off H$
and\sss $U_{\dff \mu}\dff(\trf A\trf)$\sss is\dss an\sss isometry\sss $H\qff \ttoo\qff H$\nnsp.\oss
Let\sss $V$ be related\sss to $A$ as above.\oss

\mypar{Lemma.}{cayley-t}
$U_{\dff \mu}\dff(\trf A\trf)
\off =\off
U_{\dff \mu}\dff(\trf T\trf)
\qff +\qff V_{\fff 0}$\nsp.\oss

\proof
This immediately\sss follows\sss from\trs Lemma\qss \ref{unitary-op}.\oss  \eproof

\myuppar{Changing\sss the extension $A$\nnsp.}
Let\sss
$B\dff \colon\dff \mathcal{K}_{\dff -}\qff \ttoo\qff \mathcal{K}_{\dff -}$\sss
be a self-adjoint\sss operator.\oss
Together with our boundary\sss triplet\sss $B$ defines another
self-adjoint\sss extension $A'$ of\dss $T$\dnsp,\oss
namely\sss the restriction of\dss $T^{\dff *}$\sss to\sss
$\kernel\fff (\trf \Gamma_{1}\qff -\qff B\trf \Gamma_{0}\trf)$\nnsp.\oss
The domain\sss $\mathcal{D}\dff(\trf A'\trf)$\sss is\dss described\dss by\sss the equation\vspace{3pt}
\[
\quad
-\qff \mu\trf z_{\dff -}\qff -\pff \overline{\mu}\pff V\fff z_{\dff +}
\off =\off
B\trf(\trf z_{\dff -}\qff +\pff V\fff z_{\dff +}\trf)
\pff,
\]

\vspace{-12pt}\vspace{3pt}
or,\oss equivalently,\oss by either of\dss the equations\vspace{1.5pt}\vspace{-0.25pt}
\[
\quad
(\trf B\qff +\qff \mu\trf)\trf z_{\dff -}
\qff +\qff
(\trf B\qff +\pff \overline{\mu}\qff)\trf V\fff z_{\dff +}
\off =\off
0
\pff,
\quad
z_{\dff -}
\pff +\pff
\frac{B\qff +\pff \overline{\mu}\fff}{B\qff +\qff \mu}\pff V\fff z_{\dff +}
\off =\off
0
\pff.
\]

\vspace{-12pt}\vspace{1.5pt}\vspace{-0.25pt}
Let\sss $V\fff'$ be\sss related\dss to $A'$\sss in\sss
the same way as $V$\sss to $A$\nnsp.\oss
The\sss last\sss formula shows\sss that\vspace{1.5pt}\vspace{-0.25pt}
\[
\quad
V\fff'
\off =\off
\frac{B\qff +\pff \overline{\mu}\fff}{B\qff +\qff \mu}\pff V
\]

\vspace{-12pt}\vspace{1.5pt}\vspace{-0.25pt}
Together\sss with\trs Lemma\qss \ref{cayley-t}\qss 
this\sss implies\sss that\vspace{1.5pt}
\[
\quad
U_{\dff \mu}\dff(\trf A'\trf)
\off\fff =\off\dff
\left(\qff
\frac{B\qff +\pff \overline{\mu}\fff}{B\qff +\qff \mu}
\qff\right)_{\dff H}
U_{\dff \mu}\dff(\trf A\trf)
\pff.
\]

\vspace{-12pt}\vspace{1.5pt}
Now,\oss let\sss us\sss take\sss $\mu\off =\off i$\sss
and observe\sss that\sss $U_{\fff i}\trf(\trf \bullet\trf)$\sss
is\dss the usual\trs Cayley\trs transform,\oss
which we will\sss denote by\sss $U\trf(\trf \bullet\trf)$\nnsp.\oss
Hence\sss the\sss last\sss formula\sss implies\sss that\vspace{3pt}
\begin{equation}
\label{cayley-t-i}
\quad
U\trf(\trf A'\trf)
\off =\off
U\trf(\trf B\trf)_{\dff H}\trf
U\trf(\trf A\trf)
\pff.
\end{equation}

\vspace{-12pt}\vspace{3pt}
Since\sss 
$U_{\fff -\dff i}\trf(\trf \bullet\trf)
\off =\off
U\trf(\trf \bullet\trf)^{\dff -\dff 1}$\dnsp,\oss
for\sss $\mu\off =\off -\qff i$\sss
we get\sss
$U\trf(\trf A'\trf)^{\dff -\dff 1}
\off =\off
U\trf(\trf B\trf)_{\dff H}^{\dff -\dff 1}\trf
U\trf(\trf A\trf)^{\dff -\dff 1}$\sss
and\dss hence\vspace{3pt}
\begin{equation}
\label{cayley-t-minus-i}
\quad
U\trf(\trf A'\trf)
\off =\off
U\trf(\trf A\trf)\trf
U\trf(\trf B\trf)_{\dff H}
\pff.
\end{equation}

\vspace{-12pt}\vspace{3pt}
Here\sss $U\trf(\trf A\trf)\fff,\pff
U\trf(\trf A'\trf)$ have\sss the same meaning as before,\oss
but\sss $B$\sss is\dss a operator\sss in\sss
$\kernel\fff (\trf T^{\dff *}\qff -\qff i\qff)$
and\sss not\sss in\sss
$\kernel\fff (\trf T^{\dff *}\qff +\qff i\qff)$
as before.\oss
The identity\qss (\ref{cayley-t-minus-i})\qss is\dss
a minor\sss generalization of\dss the\sss last\sss
identity\sss in\qss \cite{s},\oss Theorem\qss 14.20.\oss
The same arguments work\sss for self-adjoint\sss relations\sss
$\mathcal{B}
\qff \subset\qff
\mathcal{K}_{\dff -}\qff \oplus\qff \mathcal{K}_{\dff -}$\nsp.\oss
In\sss this case we get\sss
$U\trf(\trf A'\trf)
\off =\off
U\trf(\trf \mathcal{B}\trf)_{\dff H}\trf
U\trf(\trf A\trf)$\sss
for\sss $\mu\off =\off i$\nnsp.\oss

\newpage
\mysection{Gelfand\qss triples}{gelfand-triples}

\myuppar{Gelfand\dss triples.}
Let\sss $H$\sss be a separable\dss Hilbert\dss space and\sss
$K\qff \subset\qff H$\sss be a dense vector subspace which\dss
is\dss also a\dss Hilbert\dss space in\sss its own\sss right.\oss
Let\sss $\iota\dff \colon\dff K\qff \ttoo\qff H$\sss be\sss the inclusion.\oss
It\dss is\dss assumed\sss that\sss $\iota$\sss is\dss bounded.\oss
Let\sss $K\fff'$\sss be\sss the space of\dss anti-linear\sss maps\dss $K\qff \ttoo\qff \ccc$\nnsp.\oss
Then\sss the dual\sss of\dss $\iota$\sss is\dss the map\sss
$\iota\fff'\dff \colon\dff H\fff'\qff \ttoo\qff K\fff'$\nnsp.\oss
Since $\iota$\sss is\dss injective,\pss $\iota\fff'\dff(\trf H\trf)$\sss is\dss
dense in\sss $K\fff'$\dnsp.\oss
Since $\iota\trf(\trf K\trf)$\sss is\dss dense in\sss $H$\nnsp,\oss
the map\sss $\iota\fff'$\sss is\dss injective.\oss
We will\dss identify\sss $H$\sss with\sss $H\fff'$ in\sss the usual\sss manner,\oss
but\sss not\sss $K$\sss with\sss $K\fff'$\dnsp,\oss
despite\sss the fact\sss that\sss $K$\sss is\dss a\dss Hilbert\dss space.\oss
In\sss fact,\oss it\dss is\dss impossible\sss to satisfactory\sss identify\sss
both\sss $H$\sss with\sss $H\fff'$ and\sss $K$\sss with\sss $K\fff'$\sss simultaneously.\oss 
Still,\oss since $K$\sss is\dss a\dss Hilbert\dss space,\pss $K\fff'$\sss is\dss
also a\dss Hilbert\dss space.\oss
Usually\sss we\sss will\dss treat\sss the maps\sss $\iota$\sss
and\sss $\iota\fff'$ as inclusions\qss (of\dss sets).\oss
Then we\sss get\sss a\sss triple\sss
$K\qff \subset\qff H\qff \subset\qff K\fff'$\sss of\trs Hilbert\sss spaces,\oss
called\dss the\qss \emph{Gelfand\trs triple}\pss associated\sss with\sss the pair\sss $K\fff,\qff H$\nnsp.\oss
We will\sss denote\sss the\dss Hilbert\dss 
scalar\sss product\sss in\sss $K$\sss
by\sss $\sco{\dff \bullet\fff,\qff \bullet\dff}_{\dff K}$
and use similar notations for other\dss Hilbert\dss spaces.\oss
Let\sss\vspace{3pt}
\[
\quad
\sco{\dff \bullet\fff,\qff \bullet\dff}_{\dff K\fff'\fff,\qff K}
\qff \colon\qff
K\fff'\dff \times\dff K
\qff \ttoo\qff 
\ccc
\]

\vspace{-12pt}\vspace{3pt}
be\sss the canonical\dss pairing 
$\sco{\dff y\fff,\qff x\dff}_{\dff K\fff'\fff,\qff K}
\off =\off
y\trf(\trf x\trf)$\nnsp,\oss
If\dss
$(\trf y\fff,\qff x\trf)\qff \in\qff H\dff \times\dff K$\nnsp,\oss
then\vspace{3pt}
\[
\quad
\sco{\dff y\fff,\qff x\dff}_{\dff K\fff'\fff,\qff K}
\off =\off
\iota\fff'\fff y\trf(\trf x\trf)
\off =\off
\bsco{\dff y\fff,\pff \iota\dff x\dff}_{\dff H}
\off =\off 
\sco{\dff y\fff,\qff x\dff}_{\dff H}
\pff.
\]

\vspace{-12pt}\vspace{3pt}
In other words,\pss
$\sco{\dff \bullet\fff,\qff \bullet\dff}_{\dff K\fff'\fff,\qff K}$\sss
is\dss equal\dss to\sss
$\sco{\dff \bullet\fff,\qff \bullet\dff}_{\dff H}$\sss
on\sss $H\dff \times\dff K$\nnsp.\oss
Let\vspace{3pt}
\[
\quad
\sco{\dff \bullet\fff,\qff \bullet\dff}_{\dff K\fff,\qff K\fff'}
\qff \colon\qff
K\dff \times\dff K\fff'
\qff \ttoo\qff 
\ccc
\]

\vspace{-12pt}\vspace{3pt}
be\sss the pairing\sss
$\sco{\dff x\fff,\qff y\dff}_{\dff K\fff,\qff K\fff'}
\off =\off
\overline{y\trf(\trf x\trf)}$\nsp.\oss
If\dss
$(\trf x\fff,\qff y\trf)\qff \in\qff K\dff \times\dff H$\nnsp,\oss
then\vspace{3pt}
\[
\quad
\sco{\dff x\fff,\qff y\dff}_{\dff K\fff,\qff K\fff'}
\off =\off
\overline{\iota\fff'\fff y\trf(\trf x\trf)}
\off =\off
\overline{\sco{\dff y\fff,\qff x\dff}}_{\dff H}
\off =\off
\sco{\dff x\fff,\qff y\dff}_{\dff H}
\pff.
\]

\vspace{-12pt}\vspace{3pt}
In other words,\pss
$\sco{\dff \bullet\fff,\qff \bullet\dff}_{\dff K\fff,\qff K\fff'}$\sss
is\dss equal\dss to\sss
$\sco{\dff \bullet\fff,\qff \bullet\dff}_{\dff H}$\sss
on\sss $K\dff \times\dff H$\nnsp.\oss
The\dss Hilbert\dss space $K\fff'$ can be also constructed as\sss the completion
of\dss $H$\sss with\sss respect\sss to\sss a new scalar\sss product.\oss
Let\sss $\iota^{\dff *}\dff \colon\dff H\qff \ttoo\qff K$\sss be\sss
the operator adjoint\sss to $\iota$\nnsp,\oss
i.e.\qss such\sss that\vspace{3pt}
\[
\quad
\bsco{\dff \iota\trf x\fff,\pff y\dff}_{\dff H}
\off =\off
\bsco{\dff x\fff,\pff \iota^{\dff *} y\dff}_{\dff K}
\]

\vspace{-12pt}\vspace{3pt}
for every\sss $x\qff \in\qff K$\nnsp,\qss $y\qff \in\qff H$\nnsp.\oss
For $y\qff \in\qff H$ the\sss linear\sss functional\sss
$\iota\fff'\fff y\dff \colon\dff K\qff \ttoo\qff \ccc$\dss is\dss
equal\dss to\vspace{3pt}
\[
\quad
\iota\fff'\fff y\dff \colon\dff
x
\off \longmapsto\off
\bsco{\dff y\fff,\pff \iota\dff x\dff}_{\dff H}
\off =\off
\bsco{\dff \iota^{\dff *} y\fff,\pff x\dff}_{\dff K}
\pff.
\]

\vspace{-12pt}\vspace{3pt}
It\sss follows\sss that\sss for every\sss
$u\fff,\qff v\qff \in\qff H$\sss
the scalar\sss product\sss in\sss $K\fff'$\sss is\vspace{3pt}
\[
\quad
\bsco{\dff \iota\fff'\fff u\fff,\pff \iota\fff'\fff v\dff}_{\dff K\fff'}
\off =\off
\bsco{\dff \iota^{\dff *} u\fff,\pff \iota^{\dff *} v\dff}_{\dff K}
\pff.
\]

\vspace{-12pt}\vspace{3pt}
Therefore\sss $K\fff'$\sss
can\sss be identified with\sss the completion of\dss $H$\sss
with respect\sss to\sss the scalar\sss product\vspace{3pt}
\begin{equation}
\label{isometry}
\quad
\bsco{\dff u\fff,\qff v\dff}'
\off =\off
\bsco{\dff \iota^{\dff *} u\fff,\pff \iota^{\dff *} v\dff}_{\dff K}
\pff.
\end{equation}

\vspace{-12pt}\vspace{3pt}
In\sss particular,\pss
$\bsco{\dff \bullet\fff,\qff \bullet\dff}_{\dff K\fff'}
\off =\off
\bsco{\dff \bullet\fff,\qff \bullet\dff}'$\nnsp.\oss
The equality\qss (\ref{isometry})\qss implies\sss that\sss the operator\sss $\iota^{\dff *}$\sss defines
an\sss isometry\sss from\sss the subspace\sss $H$\sss of\dss $K\fff'$\sss into\sss $K$\nnsp.\oss
Extending\sss $\iota^{\dff *}$\sss by continuity,\oss we get\sss an\sss isometric
operator\sss $\mathbf{I}\dff \colon\dff K\fff'\qff \ttoo\qff K$\nnsp.\oss
In\sss fact,\pss $\mathbf{I}$\sss is\dss surjective.\oss
Indeed,\oss the image\sss $\image\fff \mathbf{I}$\sss is\dss closed,\oss
and\dss if\dss $x\qff \in\qff K$\sss is\dss orthogonal\dss to\sss this image,\oss
then $x$\sss is\dss orthogonal\dss to\sss $\image\fff \iota^{\dff *}$
and\dss hence $x\off =\off 0$\nnsp.\oss

It\sss turns out\sss that\sss $\mathbf{I}$\sss admits a canonical\dss
presentation as\sss the composition of\dss two isometries\sss
$K\fff'\qff \ttoo\qff H\qff \ttoo\qff K$\nnsp.\oss
In order\sss to see\sss this,\oss let\sss us consider\sss
the operator\sss
$j\off =\off \iota\dff \circ\dff \iota^{\dff *}
\dff \colon\dff
H\qff \ttoo\qff H$\nnsp.\oss
Clearly,\pss $j$\sss is\dss self-adjoint\sss and\sss
$\sco{\dff j\dff u\fff,\qff u\dff}_{\dff H}
\off =\off
\sco{\dff \iota^{\dff *} u\fff,\qff \iota^{\dff *}u\dff}_{\dff K}
\qff \geq\qff
0$\sss
for every\sss $u\qff \in\qff H$\nnsp,\oss
i.e\trs $j$\sss is\dss a non-negative operator.\oss
Hence\sss the square root\sss $\Lambda\off =\off \sqrt{j}$\sss is\dss well\sss defined.\oss
Clearly,\vspace{3pt}
\[
\quad
\bsco{\dff \Lambda\dff u\fff,\qff \Lambda\dff v \dff}_{\dff H}
\off =\off
\bsco{\dff \Lambda^{2}\dff u\fff,\qff v \dff}_{\dff H}
\off =\off
\bsco{\dff j\fff u\fff,\qff v \dff}_{\dff H}
\off =\off
\bsco{\dff \iota^{\dff *} u\fff,\qff \iota^{\dff *} v\dff}_{\dff K}
\off =\off
\bsco{\dff u\fff,\qff v\dff}_{\dff K\fff'}
\]

\vspace{-12pt}\vspace{3pt}
for every\sss $u\fff,\qff v\qff \in\qff H$\nnsp.\oss
Hence\sss $\Lambda$ defines an\sss isometry\sss from\sss the subspace\sss
$H\qff \subset\qff K\fff'$\sss into\sss $H$\nnsp.\oss
Extending $\Lambda$\sss by continuity\sss to\sss $K\fff'$\sss
we get\sss an\sss isometric operator\sss 
$\Lambda'\dff \colon\dff K\fff'\qff \ttoo\qff H$\nnsp.\oss
We claim\sss that\sss it\dss is\dss surjective.\oss
Indeed,\oss the image\sss $\image\fff \Lambda'$\sss is\dss closed,\oss
and\dss if\dss $u\qff \in\qff H$\sss is\dss orthogonal\dss to\sss this image,\oss
then $u$\sss is\dss orthogonal\sss to\sss 
$\image\fff \Lambda\off =\off \Lambda\dff(\trf H\trf)$\nnsp.\oss
Since $\Lambda$\sss is\dss self-adjoint,\oss in\sss this case\sss 
$\Lambda\dff u\off =\off 0$ and\dss hence\sss $\iota^{\dff *} u\off =\off 0$\nnsp.\oss
In\sss turn,\oss this implies\sss that\sss $u\off =\off 0$\nnsp.\oss
It\sss follows\sss that\sss $\Lambda'$\sss is\dss surjective,\oss
and\dss hence\dss is\dss an\sss isomorphism\sss $K\fff'\qff \ttoo\qff H$\nnsp.\oss
Next,\oss we claim\sss that\vspace{3pt}
\[
\quad
\Lambda\dff \circ\dff \Lambda'
\off =\off
\iota\dff \circ\dff \mathbf{I}
\pff.
\]

\vspace{-12pt}\vspace{3pt}
Indeed,\oss on\sss the subspace $H\qff \subset\qff K\fff'$\sss both compositions
are equal\dss to\sss $j\off =\off \iota\dff \circ\dff \iota^{\dff *}$\dnsp.\oss
By continuity\sss this equality extends\sss to\sss the whole space $K\fff'$\dnsp.\oss
This equality\sss implies\sss that\sss the image\sss $\image\fff \Lambda$\sss
is\dss contained\sss in\sss $\image\fff \iota\off =\off K$\nnsp.\oss
Hence we may consider\sss $\Lambda$ as an operator\sss $H\qff \ttoo\qff K$\nnsp.\oss
Then\sss the above equality\sss turns into\sss
$\Lambda\dff \circ\dff \Lambda'
\off =\off
\mathbf{I}$\nnsp,\oss
where $\Lambda$\sss is\dss a map $H\qff \ttoo\qff K$ and\sss
$\Lambda'$\sss is\dss a map\sss $K\fff'\qff \ttoo\qff H$\nnsp.\oss
Since\sss $\Lambda'$ and\sss $\mathbf{I}$ are isometric\sss isomorphisms,\oss
this implies\sss that\sss $\Lambda$\sss is\dss also an\sss isometric\sss isomorphism.\oss
The equality\sss $\mathbf{I}\off =\off \Lambda\dff \circ\dff \Lambda'$\sss
is\dss the promised\sss presentation of\dss $\mathbf{I}$ as\sss 
the composition of\dss two isometries\sss
$K\fff'\qff \ttoo\qff H\qff \ttoo\qff K$\nnsp.\oss

\mypar{Lemma.}{shift}
\emph{For every\sss $x\qff \in\qff K$\nnsp,\qss $y\qff \in\qff K\fff'$\dnsp.\oss}\vspace{3pt}
\[
\quad
\bsco{\dff y\dff,\qff x\dff}_{\dff K\fff'\fff,\qff K}
\off =\off
\bsco{\dff \Lambda'\dff y\dff,\qff \Lambda^{\fff -\dff 1}\dff x\dff}_{\dff H}
\pff.
\]

\vspace{-12pt}\vspace{3pt}
\proof
Let\sss $u\off =\off \Lambda^{\fff -\dff 1}\dff x$\nnsp.\oss
If\dss $y\qff \in\qff H$\nnsp,\oss then\vspace{3pt}
\[
\quad
\bsco{\dff y\dff,\qff x\dff}_{\dff K\fff'\fff,\qff K}
\off =\off
\bsco{\dff y\dff,\qff x\dff}_{\dff H}
\off =\off
\bsco{\dff y\dff,\qff \Lambda\dff u\dff}_{\dff H}
\off =\off
\bsco{\dff \Lambda\dff y\dff,\qff u\dff}_{\dff H}
\off =\off
\bsco{\dff \Lambda'\dff y\dff,\qff u\dff}_{\dff H}
\pff.
\]

\vspace{-12pt}\vspace{3pt}
This proves\sss the\sss lemma for $y\qff \in\qff H$\nnsp.\oss
The\sss general\sss case follows by continuity.\oss  \eproof

\mypar{Lemma.}{lambda-lambda}
\emph{The operators\sss
$\Lambda\dff \colon\dff H\qff \ttoo\qff K$\sss
and\sss
$\Lambda'\dff \colon\dff K\fff'\qff \ttoo\qff H$\sss
are adjoint\sss to each other.\oss}

\proof
We need\sss to check\sss that\sss
$\Lambda'\dff y\trf(\trf u\trf)
\off =\off
y\trf(\trf \Lambda\dff u\trf)$\sss
for every\sss $y\qff \in\qff K\fff'$\dnsp,\qss $u\qff \in\qff H$\nnsp.\oss
Recall\dss that\sss we identify $H$ with $H\fff'$ in\sss the standard way.\oss
In view of\dss this identification\vspace{1.5pt}
\[
\quad
\Lambda'\dff y\trf(\trf u\trf)
\off =\off
\bsco{\dff \Lambda'\dff y\fff,\qff u\dff}_{\dff H}
\pff.
\]

\vspace{-12pt}\vspace{1.5pt}
Also,\pss
$y\trf(\trf \Lambda\dff u\trf)
\off =\off
\bsco{\dff y\fff,\qff \Lambda\dff u\dff}_{\dff K\fff'\fff,\qff K}$\sss
by\sss the definition.\oss 
If\dss $y\qff \in\qff H$\nnsp,\oss then\sss $\Lambda'\dff y\off =\off \Lambda\dff y$
and\sss\vspace{1.5pt}
\[
\quad
\bsco{\dff y\fff,\qff \Lambda\dff u\dff}_{\dff K\fff'\fff,\qff K}
\off =\off
\bsco{\dff y\fff,\qff \Lambda\dff u\dff}_{\dff H}
\pff.
\]

\vspace{-12pt}\vspace{1.5pt}
The operator $\Lambda$\nnsp,\oss considered as an operator\sss $H\qff \ttoo\qff H$\nnsp,\oss
is\dss the square root\sss of\dss a self-adjoint\sss positive operator and\dss hence\dss
is\dss self-adjoint.\oss
Therefore,\oss if\dss $y\qff \in\qff H$\nnsp,\oss
then\vspace{2.5pt}
\[
\quad
\bsco{\dff \Lambda'\dff y\fff,\qff u\dff}_{\dff H}
\off =\off
\bsco{\dff \Lambda\dff y\fff,\qff u\dff}_{\dff H}
\off =\off
\bsco{\dff y\fff,\qff \Lambda\dff u\dff}_{\dff H}
\off =\off
\bsco{\dff y\fff,\qff \Lambda\dff u\dff}_{\dff K\fff'\fff,\qff K}
\pff.
\]

\vspace{-12pt}\vspace{2.5pt}
It\sss follows\sss that\sss
$\Lambda'\dff y\trf(\trf u\trf)
\off =\off
y\trf(\trf \Lambda\dff u\trf)$\sss
for every\sss $y\fff,\qff u\qff \in\qff H$\nnsp.\oss
By continuity\sss this equality\sss holds for every\sss 
$y\qff \in\qff K\fff'$\dnsp,\qss $u\qff \in\qff H$\nnsp.\oss
The\sss lemma\sss follows.\oss  \eproof

\myuppar{Adjoints\sss in\sss the context\sss of\trs Gelfand\trs triples.}
Let\sss us define\sss the\qss \emph{adjoint}\pss of\dss
a closed densely defined operator\sss
$B\dff \colon\dff K\fff'\qff \ttoo\qff K$ as\sss the operator 
$B^{\fff *}\dff \colon\dff K\fff'\qff \ttoo\qff K\fff''\off =\off K$\sss 
having as its domain of\dss definition\sss
$\mathcal{D}\trf(\trf B^{\fff *}\trf)$\sss
the subspace of\sss $x\qff \in\qff K\fff'$ such\sss that\sss 
the\sss linear\sss functional\sss
$a\qff \longmapsto\qff 
\sco{\dff B\dff a\dff,\qff x \dff}_{\dff K\fff,\qff K\fff'}$\sss
on\sss $\mathcal{D}\trf(\trf B\trf)$\sss
extends\sss to a continuous functional\sss on $K\fff'$ and 
such\sss that\vspace{1.5pt}
\[
\quad
\sco{\dff B\dff a\dff,\qff 
x \dff}_{\dff K\fff,\qff K\fff'}
\off =\off
\sco{\dff a\dff,\qff 
B^{\fff *} x \dff}_{\dff K\fff'\fff,\qff K}
\pff
\]

\vspace{-12pt}\vspace{1.5pt}
for every\sss $a\qff \in\qff \mathcal{D}\trf(\trf B\trf)\fff,\qff 
x\qff \in\qff \mathcal{D}\trf(\trf B^{\fff *}\trf)$\nnsp.\oss
Since $\mathcal{D}\trf(\trf B\trf)$ is\dss dense,\oss 
this uniquely determines $B^{\fff *} x$\nnsp.\oss
This notion\dss is\dss different\sss from\sss the usual\sss notion
of\dss the adjoint\sss of\dss $B$ as an operator between\sss two\dss Hilbert\dss spaces,\oss
the\sss latter\dss being an operator\sss $K\qff \ttoo\qff K\fff'$\dnsp.\oss
An operator $B\dff \colon\dff K\fff'\qff \ttoo\qff K$\sss
is\dss said\sss to be\qss \emph{self-adjoint}\pss 
if\dss it\dss is\dss equal\dss to\sss its adjoint.\oss
Lemma\qss \ref{lambda-lambda}\qss implies\sss that\sss $B$\sss is\dss self-adjoint\dss
if\dss and\dss only\trs if\dss the operator\sss
$\Lambda^{\fff -\dff 1}\dff \circ\dff
B\dff \circ\dff 
(\trf \Lambda'\qff)^{\dff -\dff 1}
\dff \colon\dff
H\qff \ttoo\qff H$\sss
is\dss self-adjoint\sss in\sss the usual\sss sense.\oss

These definitions naturally extend\sss to relations\sss
$\mathcal{B}\qff \subset\qff K\fff'\dff \oplus\dff K$\nnsp.\oss
Namely,\oss 
the\qss \emph{adjoint\sss relation}\dss
$\mathcal{B}^{\fff *}$\sss is\dss defined as\sss the space of\dss
pairs\sss $(\trf x\fff,\qff y\trf)\qff \in\qff K\fff'\dff \oplus\dff K$
such\sss that\sss\vspace{1.5pt}
\[
\quad
\sco{\dff b\dff,\qff 
x \dff}_{\dff K\fff,\qff K\fff'}
\off =\off
\sco{\dff a\dff,\qff 
y \dff}_{\dff K\fff'\fff,\qff K}
\]

\vspace{-12pt}\vspace{1.5pt}
for every\sss $(\trf a\fff,\qff b\trf)\qff \in\qff \mathcal{B}$\dnsp.\oss
Lemma\qss \ref{shift}\qss implies\sss that\sss 
$(\trf x\fff,\qff y\trf)\qff \in\qff \mathcal{B}^{\fff *}$\sss
if\dss and only\dss if\vspace{3pt}
\[
\quad
\bsco{\dff \Lambda^{\fff -\dff 1}\dff b\dff,\qff 
\Lambda'\dff x \dff}_{\dff H}
\off =\off
\bsco{\dff \Lambda'\dff a\dff,\qff 
\Lambda^{\fff -\dff 1}\dff y \dff}_{\dff H}
\]

\vspace{-12pt}\vspace{3pt}
for every\sss \sss $(\trf a\fff,\qff b\trf)\qff \in\qff \mathcal{B}$\dnsp.\oss
Then $\mathcal{B}$\sss is\dss a\qss \emph{self-adjoint\sss relation},\oss 
i.e.\qss is\dss equal\dss to its adjoint,\oss 
if\dss and\dss only\trs if\trs its\sss image\sss
$\Lambda\fff'\dff \oplus\dff \Lambda^{\fff -\dff 1}\trf
(\trf \mathcal{B}\trf)$\sss in\sss $H\dff \oplus\dff H$\sss
is\dss a self-adjoint\sss relation.\oss

\newpage
\mysection{Abstract\qss boundary\qss problems}{abstract-index}

\myuppar{The operator $A$\nnsp.}
Let\sss $H_{\dff 0}$\sss be a separable\dss Hilbert\dss space\qss
and\dss let\sss $T$\sss be a closed densely defined symmetric operator in\sss $H_{\dff 0}$\nsp.\oss
We need\dss to\sss fix a closed self-adjoint\sss extension $A$ of\dss $T$\dnsp,\oss
which we will\sss call\dss the\qss \emph{reference operator}{}.\oss
The operator $A$\sss is contained\sss in\dss $T^{\dff *}$\dnsp.\oss
We will\sss say\sss that\sss $A$\sss is\qss \emph{invertible}\pss 
if\dss $A$ has a bounded everywhere defined\sss inverse\sss
$A^{\fff -\dff 1}\dff \colon\dff H_{\dff 0}\qff \ttoo\qff H_{\dff 0}$\nsp.\oss

\mypar{Lemma.}{basic-decomp}
\emph{If\trs the operator $A$\sss is\dss invertible,\oss 
then\dss there\dss is\dss a\sss topological\sss direct\sss sum decomposition}\sss
$\mathcal{D}\dff(\trf T^{\dff *}\dff)
\off =\off
\mathcal{D}\dff(\trf A\dff)
\qff \dotplus\qff
\kernel\fff T^{\dff *}$\nsp\dnsp,\oss
\emph{where\sss the domains\sss $\mathcal{D}\dff(\trf T^{\dff *}\dff)$ and\sss $\mathcal{D}\dff(\trf A\dff)$ 
are equipped\sss with\sss graph\sss topologies.\oss
The associated\sss projection\sss
$\mathcal{D}\dff(\trf T^{\dff *}\dff)
\qff \ttoo\qff
\mathcal{D}\dff(\trf A\dff)$\sss
is\dss equal\dss to\sss $A^{\fff -\dff 1}\trf T^{\dff *}$\dnsp.\oss}

\proof
See\dss Grubb\qss \cite{g1},\oss Lemma\qss II.1.1.\oss
Let\sss us reproduce\sss the key part\sss of\dss the proof.\oss
Clearly,\oss the right\sss hand side\dss is\dss contained\sss 
in\sss the\sss left\sss hand side.\oss
If\dss $u\qff \in\qff \mathcal{D}\dff(\trf T^{\dff *}\dff)$
and $u\off =\off a\qff +\qff z$\sss
with $a\qff \in\qff \mathcal{D}\dff(\trf A\dff)$
and\sss $z\qff \in\qff \kernel\fff T^{\dff *}$\dnsp,\oss
then\sss $T^{\dff *} u\off =\off A\dff a$
and\dss hence $a\off =\off A^{\fff -\dff 1}\trf T^{\dff *} u$\nnsp.\oss
It\sss follows\sss that\sss the presentation\sss
$u\off =\off a\qff +\qff z$\nnsp,\oss if\dss exists,\oss is\dss unique.\oss
If\dss $a\off =\off A^{\fff -\dff 1}\trf T^{\dff *} u$
and\sss $z\off =\off u\qff -\qff a$\nnsp,\oss
then\sss \vspace{3pt}
\[
\quad
T^{\dff *} z
\off =\off 
T^{\dff *} u\qff -\qff A\dff a
\off =\off 
T^{\dff *} u\qff -\qff T^{\dff *} a
\off =\off
0
\]

\vspace{-12pt}\vspace{3pt}
and\dss hence\sss $z\qff \in\qff \kernel\fff T^{\dff *}$\dnsp.\oss
This proves\sss the above decomposition on\sss the algebraic\sss level.\oss
Passing\sss to\sss the\sss topological\sss decomposition\dss 
is\dss fairly\sss routine.\oss  \eproof

\myuppar{Notations.}
When\sss the operator\sss $A$\sss is\dss invertible,\oss
we will\sss denote by\sss
$p\dff \colon\dff
\mathcal{D}\dff(\trf T^{\dff *}\dff)
\qff \ttoo\qff
\mathcal{D}\dff(\trf A\dff)$\sss
and\sss
$k\dff \colon\dff
\mathcal{D}\dff(\trf T^{\dff *}\dff)
\qff \ttoo\qff
\kernel\fff T^{\dff *}$\sss
the projections associated with\sss the decomposition of\trs Lemma\qss \ref{basic-decomp}.\oss
We will\sss denote\sss 
$\sco{\dff \bullet\fff,\qff \bullet \dff}_{\dff H_{\dff 0}}$\sss
simply\sss by\sss $\sco{\dff \bullet\fff,\qff \bullet \dff}$\nnsp.\oss

\mypar{Lemma.}{modified-l}
\emph{Suppose\sss that\sss $A$\sss is\dss a reference operator{}.\oss
If\dss $u\fff,\qff v\qff \in\qff \mathcal{D}\dff(\trf T^{\dff *}\dff)$\nnsp,\oss
then}\vspace{1.5pt}
\[
\quad
\sco{\dff T^{\dff *} u\fff,\qff v \dff}
\pff -\pff
\sco{\dff u\dff,\qff T^{\dff *} v \dff}
\off =\off
\bsco{\dff T^{\dff *} u\dff,\qff 
k\trf(\trf v\trf) \dff}
\pff -\pff
\bsco{\dff k\trf(\trf u\trf)\dff,\qff 
T^{\dff *}  v \dff}
\pff.
\]

\vspace{-12pt}\vspace{1.5pt}
\proof
Since\sss 
$u\off =\off p\trf(\trf u\trf)\qff +\qff k\trf(\trf u\trf)$
and\sss
$v\off =\off p\trf(\trf v\trf)\qff +\qff k\trf(\trf v\trf)$\nnsp,\oss
the\sss left\sss hand side\dss is\dss equal\dss to\vspace{3pt}
\[
\quad
\bsco{\dff T^{\dff *} u\dff,\qff 
p\trf(\trf v\trf)\qff +\qff k\trf(\trf v\trf) \dff}
\pff -\pff
\bsco{\dff p\trf(\trf u\trf)\qff +\qff k\trf(\trf u\trf)\dff,\qff 
T^{\dff *}  v \dff}
\]

\vspace{-34.25pt}
\[
\quad
=\off
\bsco{\dff T^{\dff *} u\dff,\qff 
p\trf(\trf v\trf)\dff}
\qff +\qff 
\bsco{\dff T^{\dff *} u\dff,\qff
k\trf(\trf v\trf) \dff}
\pff -\pff
\bsco{\dff p\trf(\trf u\trf)\dff,\qff
T^{\dff *}  v \dff}
\qff -\qff 
\bsco{\dff k\trf(\trf u\trf)\dff,\qff 
T^{\dff *}  v \dff}
\pff.
\]

\vspace{-12pt}\vspace{3pt}
Since\sss 
$T^{\dff *} k\trf(\trf u\trf)
\off =\off
T^{\dff *} k\trf(\trf v\trf)
\off =\off 
0$\nnsp,\oss
the\sss last\sss expression\dss is\dss equal\dss to\vspace{3pt}
\[
\quad
\bsco{\dff T^{\dff *} p\trf(\trf u\trf)\dff,\qff 
p\trf(\trf v\trf)\dff}
\qff +\qff 
\bsco{\dff T^{\dff *} u\dff,\qff
k\trf(\trf v\trf) \dff}
\pff -\pff
\bsco{\dff p\trf(\trf u\trf)\dff,\qff
T^{\dff *}  p\trf(\trf v\trf) \dff}
\qff -\qff 
\bsco{\dff k\trf(\trf u\trf)\dff,\qff 
T^{\dff *}  v \dff}
\pff
\]

\vspace{-34.25pt}
\[
\quad
=\off
\bsco{\dff A\dff p\trf(\trf u\trf)\dff,\qff 
p\trf(\trf v\trf)\dff}
\qff -\qff
\bsco{\dff p\trf(\trf u\trf)\dff,\qff
A\dff p\trf(\trf v\trf) \dff}
\qff +\qff 
\bsco{\dff T^{\dff *} u\dff,\qff
k\trf(\trf v\trf) \dff}
\qff -\qff 
\bsco{\dff k\trf(\trf u\trf)\dff,\qff 
T^{\dff *} v \dff}
\pff.
\]

\vspace{-12pt}\vspace{3pt}
Since $A$\sss is\dss a self-adjoint\sss operator 
and\sss 
$p\trf(\trf u\trf)\fff,\pff p\trf(\trf v\trf)
\qff \in\qff
\mathcal{D}\dff(\trf A\dff)$\nnsp,\oss\vspace{0pt}
\[
\quad
\bsco{\dff A\dff p\trf(\trf u\trf)\dff,\qff 
p\trf(\trf v\trf)\dff}
\qff -\qff
\bsco{\dff p\trf(\trf u\trf)\dff,\qff
A\dff p\trf(\trf v\trf) \dff}
\off =\off
0
\pff.
\]

\vspace{-12pt}\vspace{0pt}
The\sss lemma follows.\oss  \eproof

\myuppar{Boundary\sss operators.}
Let\sss $H_{\fff 1}$\sss
be a dense subspace of\dss $H_{\dff 0}$\nsp,\oss
which\dss is\dss a\sss Hilbert\sss space in\sss its own\sss right.\oss
Let\sss $K^{\dff \partial}$\sss be another separable\dss Hilbert\dss space
and\sss 
and\sss $K$\sss be a dense subspace of\dss $K^{\dff \partial}$\dnsp,\oss
which\dss is\dss a\sss Hilbert\sss space in\sss its own\sss right.\oss
Suppose\sss that\qss \emph{the inclusion maps\sss 
$H_{\fff 1}\qff \ttoo\qff H_{\dff 0}$\sss
and\sss
$K\qff \ttoo\qff K^{\dff \partial}$\sss
are bounded operators with\sss respect\sss
to\sss these\sss Hilbert\sss space structures.}\oss
Let\sss
$K\qff \subset\qff K^{\dff \partial}\qff \subset\qff K\fff'$\sss
be\sss the\dss Gelfand\trs triple associated\sss with\sss the pair\sss
$K\fff,\pff K^{\dff \partial}$\dnsp.\oss

We will\sss denote by\sss
$\sco{\dff \bullet\dff,\qff \bullet \dff}_{\dff \partial}$\sss
the scalar product\sss in\sss  $K^{\dff \partial}$\dnsp.\oss
\emph{Suppose\sss that\dss
$H_{\dff 1}\qff \subset\qff \mathcal{D}\dff(\trf T^{\dff *}\trf)$\nnsp,\qss
$H_{\dff 1}$\sss is\dss dense in\sss $\mathcal{D}\dff(\trf T^{\dff *}\trf)$\nnsp,\oss
the operator\sss
$H_{\dff 1}\qff \ttoo\qff H_{\dff 0}$\sss
induced\dss by\dss $T^{\dff *}$ is\dss bounded,\oss
and}\vspace{0pt}
\[
\quad
\gamma_0\dff,\pff \gamma_1\dff \colon\dff
H_{\dff 1}\qff \ttoo\qff K
\off \subset\off
K^{\dff \partial}
\]

\vspace{-12pt}\vspace{0pt}
\emph{are bounded operators such\sss that\dss 
$\gamma
\off =\off
\gamma_0\dff \oplus\dff \gamma_1\dff \colon\dff
H_{\dff 1}\qff \ttoo\qff K\dff \oplus\dff K$\sss
is\dss surjective and}\vspace{3pt}
\begin{equation}
\label{lagrange-identity}
\quad
\sco{\dff T^{\dff *} u\fff,\qff v \dff}
\pff -\pff
\sco{\dff u\dff,\qff T^{\dff *} v \dff}
\off =\off
\bsco{\dff \gamma_1\dff u\dff,\qff 
\gamma_0\dff v \dff}_{\dff \partial}
\pff -\pff
\bsco{\dff \gamma_0\dff u\dff,\qff 
\gamma_1\dff v \dff}_{\dff \partial}
\end{equation}

\vspace{-12pt}\vspace{3pt}
\emph{for every $u\fff,\qff v\qff \in\qff H_{\dff 1}$\nsp.}\oss
Then $\gamma\dff \colon\dff H_{\fff 1}\qff \ttoo\qff K\dff \oplus\dff K$\sss
admits a continuous section,\oss i.e.\qss there exists a bounded operator\sss
$\kappa\dff \colon\dff K\dff \oplus\dff K\qff \ttoo\qff H_{\fff 1}$\sss
such\sss that\sss $\gamma\dff \circ\dff \kappa$\sss is\dss equal\dss to\sss
the identity\sss map.\oss
\emph{Suppose further\sss that\sss $\kernel\fff \gamma$\sss is\dss dense in\sss $H_{\dff 0}$\nsp.}\oss
\emph{We will\sss also assume\sss that\qss
$\mathcal{D}\dff(\trf T\trf)
\off =\off 
\kernel\fff \gamma
\off =\off
\kernel\fff \gamma_0\dff \oplus\dff \gamma_1$\nsp.}\oss
The subspace\sss $H_{\dff 1}$\sss is\dss usually
strictly smaller\sss than $\mathcal{D}\dff(\trf T^{\dff *}\trf)$\nnsp,\oss
and\sss the map\sss
$\gamma
\off =\off
\gamma_0\dff \oplus\dff \gamma_1$ 
considered as a map\sss
$H_{\dff 1}\qff \ttoo\qff K^{\dff \partial}\dff \oplus\dff K^{\dff \partial}$\sss
is\dss not\sss surjective unless $K\off =\off K^{\dff \partial}$\dnsp.\oss
By\sss these reasons\sss the\qss \emph{boundary\sss operators}\qss 
$\gamma_0\dff,\pff \gamma_1$\sss do not\sss
define a boundary\sss triplet\sss for $T^{\dff *}$\dnsp.\oss

\mypar{Theorem.}{extending-gamma}
\emph{The operators\sss $\gamma_0\dff,\pff \gamma_1$\sss extend\dss by continuity\sss
to bounded operators}\vspace{1.5pt}
\[
\quad
\Gamma_0\dff,\pff \Gamma_1\dff \colon\dff
\mathcal{D}\dff(\trf T^{\dff *}\trf)
\qff \ttoo\qff 
K\fff'
\]

\vspace{-12pt}\vspace{1.5pt}
\emph{(where\sss $\mathcal{D}\dff(\trf T^{\dff *}\trf)$\sss
is\dss equipped\sss with\sss the graph\dss topology)\qss 
such\sss that\sss the\trs Lagrange\dss identity}\vspace{3pt}
\begin{equation}
\label{lagrange-asymm}
\quad
\sco{\dff T^{\dff *} u\fff,\qff v \dff}
\pff -\pff
\sco{\dff u\dff,\qff T^{\dff *} v \dff}
\off =\off
\bsco{\dff \Gamma_1\dff u\dff,\qff 
\Gamma_0\dff v \dff}_{\dff K\fff'\fff,\qff K}
\pff -\pff
\bsco{\dff \Gamma_0\dff u\dff,\qff 
\Gamma_1\dff v \dff}_{\dff K\fff'\fff,\qff K}
\end{equation}

\vspace{-12pt}\vspace{3pt}
\emph{holds\sss for every\sss $u\qff \in\qff \mathcal{D}\dff(\trf T^{\dff *}\trf)\dff,\qff
v\qff \in\qff H_{\dff 1}$\nsp.}\oss

\proof
The proof\dss is\dss based on\sss ideas of\trs Lions\dss and\dss Magenes\qss \cite{lm3}.\oss
See\qss \cite{lm3},\oss the proof\dss of\trs Theorem\qss 3.1.\oss
Let\sss us\sss temporarily\sss fix some\sss 
$u\qff \in\qff \mathcal{D}\dff(\trf T^{\dff *}\trf)$\nnsp.\oss
Given\sss $\varphi\qff \in\qff K$\nnsp,\oss let\sss us\sss choose some\sss $w\qff \in\qff H_{\dff 1}$
such\dss that\sss $\gamma_0\dff w\off =\off 0$ and\sss $\gamma_1\dff w\off =\off \varphi$\sss
and\sss set\vspace{3pt}
\[
\quad
Y^{\dff w}\dff(\trf \varphi\trf)
\off =\off
\sco{\dff u\dff,\qff T^{\dff *} w \dff}
\pff -\pff
\sco{\dff T^{\dff *} u\fff,\qff w \dff}
\pff.
\]

\vspace{-12pt}\vspace{3pt}
We claim\sss that\sss $Y^{\dff w}\dff(\trf \varphi\trf)$ does not\sss depend on\sss
the choice of\dss $w$\nnsp.\oss
Indeed,\oss if\dss $w_{\dff 1}$\sss is\dss some other choice and\sss
$d\off =\off w\qff -\qff w_{\dff 1}$\nsp,\oss
then\sss $\gamma_0\dff d\off =\off \gamma_1\dff d\off =\off 0$\sss
and\dss hence\sss $d\qff \in\qff \mathcal{D}\dff(\trf T\trf)$\nnsp.\oss
It\sss follows\sss that\sss
$\sco{\dff u\dff,\qff T\dff d \dff}
\pff -\pff
\sco{\dff T^{\dff *} u\fff,\qff d \dff}
\off =\off
0$\nnsp.\oss
At\sss the same\sss time\vspace{1.5pt}
\[
\quad
\sco{\dff T^{\dff *} u\fff,\qff d \dff}
\pff -\pff
\sco{\dff u\dff,\qff T\dff d \dff}
\off =\off
Y^{\dff w}\dff(\trf \varphi\trf)
\qff -\qff
Y^{\dff w_{\dff 1}}\dff(\trf \varphi\trf)
\pff,
\]

\vspace{-12pt}\vspace{1.5pt}
and\dss therefore\sss
$Y^{\dff w}\dff(\trf \varphi\trf)
\off =\off
Y^{\dff w_{\dff 1}}\dff(\trf \varphi\trf)$\nnsp.\oss
The claim\sss follows.\oss
Now we can\sss set\sss 
$Y\trf(\trf \varphi\trf)
\off =\off 
Y^{\dff w}\dff(\trf \varphi\trf)$\sss
for an arbitrary choice of\dss $w$\nnsp.\oss
Clearly,\oss the map\sss $\varphi\off \longmapsto\off Y\trf(\trf \varphi\trf)$\sss
is\dss anti-linear.\oss
Moreover,\oss it\dss is\dss continuous because\sss the section\sss $\kappa$\sss
allows\sss to choose $w$ continuously depending on $\varphi$\nnsp.\oss
Therefore\sss the map\sss $\varphi\off \longmapsto\off Y\trf(\trf \varphi\trf)$\sss
belongs\sss to\sss $K\fff'$\dnsp.\oss
In other\sss terms,\vspace{1.5pt}
\[
\quad
Y\trf(\trf \varphi\trf)
\off =\off
\sco{\dff \tau\dff u\fff,\qff \varphi\dff}_{\dff K\fff'\fff,\qff K}
\]

\vspace{-12pt}\vspace{1.5pt}
for some\sss $\tau\dff u\qff \in\qff K\fff'$\dnsp.\oss
Let\sss us denote by\sss $\norm{\bullet}_{\trf T^{\dff *}}$\sss
the graph\sss norm\sss in\sss $\mathcal{D}\dff(\trf T^{\dff *}\trf)$\nnsp,\oss
by\sss $\norm{\bullet}$\sss the norm\sss in\sss $H_{\dff 0}$\nsp,\oss
and\dss by\sss $\norm{\bullet}_{\dff K}$\sss the norm\sss in\sss $K$\nnsp.\oss 
If\dss $w\off =\off \kappa\trf(\trf \varphi\trf)$\nnsp,\oss
then\vspace{3pt}
\[
\quad
\num{\sco{\dff \tau\dff u\fff,\qff \varphi\dff}_{\dff K\fff'\fff,\qff K}}
\off \leq\off
\num{\sco{\dff u\dff,\qff T^{\dff *} w \dff}}
\off +\off
\num{\sco{\dff T^{\dff *} u\fff,\qff w \dff}}
\]

\vspace{-34.5pt}
\[
\quad
\leq\off
\norm{\dff u\dff}\dff \cdot\dff \norm{\dff T^{\dff *} w \dff}
\off +\off
\norm{\dff T^{\dff *} u\dff}\dff \cdot\dff \norm{w}
\off \leq\off
\norm{\dff u\dff}_{\trf T^{\dff *}}\qff 
\left(\qff 
\norm{\dff T^{\dff *} w \dff}
\off +\off
\norm{w}
\qff\right)
\]

\vspace{-34.5pt}
\[
\quad
=\off
\norm{\dff u\dff}_{\trf T^{\dff *}}\qff 
\left(\qff 
\norm{\dff T^{\dff *} \kappa\trf(\trf 0\fff,\qff \varphi\trf) \dff}
\off +\off
\norm{\kappa\trf(\trf 0\fff,\qff \varphi\trf)}
\qff\right)
\off \leq\off
\norm{\dff u\dff}_{\trf T^{\dff *}}\qff 
\left(\qff 
C\qff \norm{\varphi}_{\dff K}
\off +\off
C\fff'\qff \norm{\varphi}_{\dff K}
\qff\right)
\qff,
\]

\vspace{-12pt}\vspace{3pt}
where\sss $C$\sss and\sss $C\fff'$ are\sss the norms of\dss the composition of\dss 
the section\sss
$\kappa\dff \colon\dff K\dff \oplus\dff K\qff \ttoo\qff H_{\dff 1}$\sss
with\sss the map\sss $H_{\dff 1}\qff \ttoo\qff H_{\dff 0}$
induced\dss by\sss $T^{\dff *}$\sss and\sss with\sss the inclusion\sss
$H_{\dff 1}\qff \ttoo\qff H_{\dff 0}$\sss respectively.\oss
It\sss follows\sss that\sss $u\off \longmapsto\off \tau\dff u$\sss
is\dss a continuous map\sss
$\mathcal{D}\dff(\trf T^{\dff *}\trf)\qff \ttoo\qff K\fff'$\dnsp.\oss

Let\sss us\sss prove now\sss that\sss $\tau$\sss extends $\gamma_0$\nsp.\oss
If\dss $u\qff \in\qff H_{\dff 1}$\nsp,\qss $\varphi\qff \in\qff K$\nnsp,\oss
and $w$\sss is\dss as above,\oss
then\vspace{1.5pt}
\[
\quad
\sco{\dff \tau\dff u\fff,\qff \varphi\dff}_{\dff K\fff'\fff,\qff K}
\off =\off
\sco{\dff u\dff,\qff T^{\dff *} w \dff}
\pff -\pff
\sco{\dff T^{\dff *} u\fff,\qff w \dff}
\]

\vspace{-34.5pt}
\[
\quad
\off =\off
\bsco{\dff \gamma_0\dff u\dff,\qff 
\gamma_1\dff w \dff}_{\dff \partial}
\pff -\pff
\bsco{\dff \gamma_1\dff u\dff,\qff 
\gamma_0\dff w \dff}_{\dff \partial}
\off =\off
\bsco{\dff \gamma_0\dff u\dff,\qff 
\varphi \dff}_{\dff \partial}
\off =\off
\bsco{\dff \gamma_0\dff u\dff,\qff 
\varphi \dff}_{\dff K\fff'\fff,\qff K}
\off,
\]

\vspace{-12pt}\vspace{1.5pt}
where\sss the\sss last\sss equality\sss holds
because\sss $\gamma_0\dff u\qff \in\qff K^{\dff \partial}$
and\sss $\varphi\qff \in\qff K$\nnsp.\oss
Since\sss $\varphi\qff \in\qff K$\nnsp,\oss
it\sss follows\sss that\sss $\tau\dff u\off =\off \gamma_0\dff u$\nnsp.\oss
Therefore\sss $\tau$\sss is\dss a continuous extension of\dss $\gamma_0$\nsp.\oss
Let\sss us set\sss $\Gamma_0\off =\off \tau$\sss
and\dss let\sss us construct\sss the extension\sss $\Gamma_1$\sss
in\sss the same way.\oss
Since $H_{\dff 1}$\sss is\dss dense in\sss $\mathcal{D}\dff(\trf T^{\dff *}\trf)$\sss
by our assumptions,\oss these extensions are unique and are\sss the extensions by\sss the continuity.\oss
Since\sss $\sco{\dff \bullet\fff,\qff \bullet\dff}_{\dff K\fff,\qff K\fff'}$\sss
is\dss equal\dss to\sss
$\sco{\dff \bullet\fff,\qff \bullet\dff}_{\dff \partial}$\sss
on\sss $K\dff \times\dff K^{\dff \partial}$\dnsp,\oss
the identity\qss (\ref{lagrange-identity})\qss implies\sss that\vspace{2.5pt}
\[
\quad
\sco{\dff T^{\dff *} u\fff,\qff v \dff}
\pff -\pff
\sco{\dff u\dff,\qff T^{\dff *} v \dff}
\off =\off
\bsco{\dff \Gamma_1\dff u\dff,\qff 
\Gamma_0\dff v \dff}_{\dff K\fff'\fff,\qff K}
\pff -\pff
\bsco{\dff \Gamma_0\dff u\dff,\qff 
\Gamma_1\dff v \dff}_{\dff K\fff'\fff,\qff K}
\]

\vspace{-12pt}\vspace{2.5pt}
for every $u\fff,\qff v\qff \in\qff H_{\dff 1}$\nsp.\oss
In contrast\sss with\qss (\ref{lagrange-identity})\qss
both sides of\dss this equality make sense also for\sss
$u\qff \in\qff \mathcal{D}\dff(\trf T^{\dff *}\trf)$ and\sss $v\qff \in\qff H_{\dff 1}$\nsp.\oss
By continuity\sss this equality\sss extends\sss to such $u\fff,\qff v$\nnsp.\oss
In other words\qss (\ref{lagrange-asymm})\qss holds for every\sss $u\qff \in\qff \mathcal{D}\dff(\trf T^{\dff *}\trf)\dff,\qff
v\qff \in\qff H_{\dff 1}$\nsp.\oss  \eproof

\mypar{Lemma.}{gamma-0-ker}
\emph{If\qss 
$\mathcal{D}\dff(\trf A\trf)\off =\off \kernel\fff \gamma_0$\nnsp,\oss
then
$\kernel\fff \Gamma_0\off =\off \kernel\fff \gamma_0$\dss
and\dss hence\trs
$\kernel\fff \Gamma_0\qff \subset\qff H_{\dff 1}$\nsp.\oss}

\proof
If\dss $u\qff \in\qff \kernel\fff \gamma_0$
and\sss $v\qff \in\qff \kernel\fff \Gamma_0$\nsp,\oss
then\qss (\ref{lagrange-asymm})\qss implies\sss that\sss
$\sco{\dff A\dff u\fff,\qff v \dff}
\pff -\pff
\sco{\dff u\dff,\qff T^{\dff *} v \dff}
\off =\off
0$\nnsp.\oss
It\sss follows\sss that\sss the map\sss 
$u\off \longmapsto\off \sco{\dff A\dff u\fff,\qff v \dff}$\sss
extends\sss to a continuous map\sss
$H\qff \ttoo\qff \ccc$\nnsp.\oss 
Hence\sss
$u$ belongs\sss to\sss the domain $\mathcal{D}\dff(\trf A^*\dff)$
of\dss $A^*$\nnsp.\oss
Since\sss $A$\sss is\dss self-adjoint,\oss
it\sss follows\sss that\sss
$u\qff \in\qff \mathcal{D}\dff(\trf A\trf)$
and\dss hence
$u\qff \in\qff \kernel\fff \gamma_0$\nsp.\oss
Therefore\sss 
$\kernel\fff \Gamma_0\qff \subset\qff \kernel\fff \gamma_0$\nsp.\oss
The opposite inclusion\dss is\dss obvious.\oss  \eproof

\mypar{Theorem.}{extended-isomorphism}
\emph{If\qss $\mathcal{D}\dff(\trf A\trf)\off =\off \kernel\fff \gamma_0$\sss
and $A$\sss is\dss invertible,\oss
then\sss $\Gamma_0$\sss induces a\sss topological\dss isomorphism\sss
$\kernel\dff T^{\dff *}\qff \ttoo\qff K\fff'$\dnsp,\oss
and\dss
$T^{\dff *} \oplus\dff \Gamma_0
\qff \colon\dff
\mathcal{D}\dff(\trf T^{\dff *}\dff)
\qff \ttoo\qff
H\dff \oplus\dff K\fff'$\sss
is\dss a\sss topological\dss isomorphism.\oss}

\proof
The proof\dss is\dss based on\sss ideas of\trs Lions\dss and\dss Magenes\qss \cite{lm2}.\oss
See\qss \cite{lm2},\oss the proof\dss of\trs Theorem\qss 9.2.\oss
Lemma\qss \ref{gamma-0-ker}\qss implies\sss that\sss
$A\off =\off T^{\dff *}\trf|\trf \kernel\fff \Gamma_0$\nsp.\oss
By our assumptions $A$ induces a\sss topological\sss isomorphism\sss
$\mathcal{D}\dff(\trf A\trf)\qff \ttoo\qff H_{\dff 0}$\nsp.\oss
By combining\sss this with\dss Lemma\qss \ref{basic-decomp},\oss
we see\sss that\sss it\dss is\dss sufficient\sss to prove\sss that\sss
$\Gamma_0\trf|\trf \kernel\fff T^{\dff *}$\sss 
is\dss a\sss topological\sss isomorphism\sss
$\kernel\fff T^{\dff *}\qff \ttoo\qff K\fff'$\dnsp.\oss
The kernel\sss of\dss 
$\Gamma_0\trf|\trf \kernel\fff T^{\dff *}$\sss
is\dss equal\dss to\sss the kernel\sss of\dss $A$ and\dss hence\dss is\dss
equal\sss to $0$\sss by our assumptions.\oss
Since\sss $\Gamma_0$\sss is\dss continuous,\oss it\dss is\dss
sufficient\sss to prove\sss that\sss the induced\sss map\dss is\dss surjective.\oss
Given $x\qff \in\qff K\fff'$\dnsp,\oss
let\sss us consider\sss the anti-linear functional\sss
$l\dff \colon\dff
\mathcal{D}\dff(\trf A\trf)\qff \ttoo\qff \ccc$
defined\dss by\vspace{3pt}
\[
\quad
l\dff \colon\dff
v
\off \longmapsto\off 
\sco{\dff x\dff,\qff \gamma_1\dff v \dff}_{\dff K\fff'\fff,\qff K}
\pff.
\]

\vspace{-12pt}\vspace{3pt}
Since $\gamma_1$\sss is\dss continuous,\pss $l$\sss is\dss also continuous.\oss
Since $A$ induces a\sss topological\sss isomorphism\sss
$\mathcal{D}\dff(\trf A\trf)\qff \ttoo\qff H_{\dff 0}$\nsp,\oss
the functional\sss $l$\sss is\dss equal\dss to\sss
the functional\sss
$v\off \longmapsto\off \sco{\dff u\fff,\qff A\dff v\dff}$\sss
for a unique $u\qff \in\qff H_{\dff 0}$\nsp.\oss
If\dss 
$v
\qff \in\qff 
\kernel\fff \gamma
\off =\off
\kernel\fff \gamma_0
\qff \cap\qff\fff
\kernel\fff \gamma_1$\nsp,\oss
then\sss $l\trf(\trf v\trf)\off =\off 0$
and\dss hence\sss 
$\sco{\dff u\fff,\qff A\dff v\dff}
\off =\off
0$\nnsp.\oss
In view of\dss our assumptions\sss this means\sss that\sss
$\sco{\dff u\fff,\qff T\dff v\dff}\off =\off 0$\sss
for every\sss $v\qff \in\qff \mathcal{D}\dff(\trf T\trf)$\nnsp.\oss
Therefore\sss $u\qff \in\qff \mathcal{D}\dff(\trf T^{\dff *}\trf)$
and\sss $T^{\dff *} u\off =\off 0$\nnsp,\oss
i.e.\dss $u\qff \in\qff \kernel\fff T^{\dff *}$\dnsp.\oss
Now\qss (\ref{lagrange-asymm})\qss implies\sss that\vspace{4.5pt}
\[
\quad
-\pff
\sco{\dff u\dff,\qff T^{\dff *} v \dff}
\off =\off
\bsco{\dff \Gamma_1\dff u\dff,\qff 
\Gamma_0\dff v \dff}_{\dff K\fff'\fff,\qff K}
\pff -\pff
\bsco{\dff \Gamma_0\dff u\dff,\qff 
\Gamma_1\dff v \dff}_{\dff K\fff'\fff,\qff K}
\]

\vspace{-12pt}\vspace{4.5pt}
for every\sss $v\qff \in\qff H_{\dff 1}$\nnsp.\oss
If\dss
$v\qff \in\qff \mathcal{D}\dff(\trf A\trf)\off =\off \kernel\fff \gamma_0$\nsp,\oss
then\sss
$\sco{\dff u\dff,\qff T^{\dff *} v \dff}
\off =\off
l\trf(\trf v\trf)$\nnsp,\oss
and\dss hence\vspace{4.5pt}
\[
\quad
\sco{\dff x\dff,\qff \gamma_1\dff v \dff}_{\dff K\fff'\fff,\qff K}
\off =\off
\bsco{\dff \Gamma_0\dff u\dff,\qff 
\Gamma_1\dff v \dff}_{\dff K\fff'\fff,\qff K}
\off =\off
\sco{\dff \Gamma_0\dff u\dff,\qff 
\gamma_1\dff v \dff}_{\dff K\fff'\fff,\qff K}
\pff.
\]

\vspace{-12pt}\vspace{4.5pt}
Since $\gamma\off =\off \gamma_0\dff \oplus\dff \gamma_1$\sss
is\dss a map onto\sss $K\dff \oplus\dff K$\nnsp,\oss
the boundary map $\gamma_1$ maps $\kernel\fff \gamma_0$\sss
onto\sss $K$\nnsp.\oss
Therefore\sss the\sss last\sss displayed equality\sss implies\sss that\sss
$x\off =\off \Gamma_0\dff u$\nnsp.\oss
The surjectivity\sss follows.\oss  \eproof

\myuppar{The reference operator and\dss the boundary operators.}
Let\sss us say\sss that\sss an\sss unbounded self-adjoint\sss operator $P$\sss in\sss $H_{\dff 0}$\sss
is\qss \emph{elliptic regular}\pss if\trs 
$\mathcal{D}\dff(\trf P\trf)\qff \subset\qff H_{\dff 1}$\nsp.\oss 
\emph{For\sss the rest\sss of\dss this section we will\sss assume\sss
that\sss the reference operator $A$\sss is\dss elliptic regular{}.\oss 
Moreover{},\oss we will\sss assume\sss that\sss 
$\mathcal{D}\dff(\trf A\trf)\off =\off \kernel\fff \gamma_0$\nsp.}\oss
The\sss last\sss assumption\sss has\sss technical\sss character
and\sss usually can be achieved\dss by changing\sss the boundary operators\sss
$\gamma_0\dff,\pff \gamma_1$\nsp.\oss
\emph{We will\sss also assume\sss that\sss $A$\sss is\dss invertible
and\sss that\trs 
$\mathcal{D}\dff(\trf T\trf)
\off =\off 
\kernel\fff \gamma
\off =\off
\kernel\fff \gamma_0\dff \oplus\dff \gamma_1$\nsp.\oss}

\myuppar{Where such\sss reference operators come from\fff?}
Let\sss 
$H^{\dff \partial}
\off =\off
K^{\dff \partial}\dff \oplus\dff K^{\dff \partial}$ and\sss
$H^{\dff \partial}_{\dff 1/2}
\off =\off
K\dff \oplus\dff K$\nnsp,\oss
and\dss let\sss 
$\Sigma\dff \colon\dff H^{\dff \partial}\qff \ttoo\qff H^{\dff \partial}$\sss
be\sss the operator such\dss that\vspace{-1pt}
\[
\quad
i\trf \Sigma
\off =\off\dff
\begin{pmatrix}
\off 0 &
1 \qff\off
\vspace{4.5pt} \\
\off\dff -\qff 1 &
0 \qff\off 
\end{pmatrix}
\off
\]

\vspace{-12pt}\vspace{-1pt}
with respect\sss to\sss the decomposition\sss
$H^{\dff \partial}
\off =\off
K^{\dff \partial}\dff \oplus\dff K^{\dff \partial}$\dnsp.\oss
Clearly,\pss $\Sigma$\sss leaves $H^{\dff \partial}_{\dff 1/2}$ invariant.\oss
In\sss these\sss terms\sss the\dss Lagrange\dss identity\qss (\ref{lagrange-identity})\qss
takes\sss the form\vspace{3pt}
\[
\quad
\sco{\dff T^{\dff *} u\fff,\qff v \dff}
\pff -\pff
\sco{\dff u\dff,\qff T^{\dff *} v \dff}
\off =\off
\bsco{\dff i\trf \Sigma\trf \gamma\dff u\dff,\qff 
\gamma\dff v \dff}_{\dff \partial}
\pff,
\]

\vspace{-12pt}\vspace{1.5pt}
where $u\fff,\qff v\qff \in\qff H_{\dff 1}$\nsp.\oss
Let\sss 
$\Pi\dff \colon\dff H^{\dff \partial}\qff \ttoo\qff H^{\dff \partial}$\sss
be\sss the projection onto\sss the second summand of\dss the decomposition\sss
$H^{\dff \partial}
\off =\off
K^{\dff \partial}\dff \oplus\dff K^{\dff \partial}$\dnsp.\oss
Then\sss the pair\sss $T^{\dff *}\trf|\trf H_{\dff 1}\dff,\off \Pi$\sss
is\dss a\qss \emph{boundary\dss problem}\qss in\sss the sense of\pss \cite{i2},\oss Section\qss 5.\oss
Since,\oss clearly,\pss
$\Sigma\trf(\trf \image\fff \Pi\trf)\off =\off \kernel\fff \Pi$\nnsp,\oss
this boundary\sss problem\dss is\dss self-adjoint\sss in\sss the sense of\pss \cite{i2}.\oss
The unbounded operator\sss induced\dss by\sss this boundary\sss problem\dss is\dss
nothing else but\sss 
$A\off =\off T^{\dff *}\trf|\trf \kernel\fff \gamma_0$\dss
(in\qss \cite{i2}\qss it\dss is\dss denoted\dss by\sss $A_{\dff \Gamma}$\nsp).\oss

Suppose now\sss that\sss the boundary\sss problem\sss
$T^{\dff *}\trf|\trf H_{\dff 1}\dff,\off \Pi$\sss
is\qss \emph{elliptic\sss regular}\pss in\sss the sense of\pss \cite{i2}\qss
(we will\sss not\sss need\sss the precise definition).\oss
Suppose also\sss that\sss the inclusion\sss
$H_{\dff 1}\qff \ttoo\qff H_{\dff 0}$\sss is\dss a compact\sss operator{}.\oss
If,\oss furthermore,\oss the operator\vspace{3pt}
\[
\quad
(\trf T^{\dff *}\trf|\trf H_{\dff 1}\trf)\dff \oplus\dff \gamma_0
\dff \colon\dff
H_{\dff 1}
\qff \ttoo\qff
H_{\dff 0}\dff \oplus\dff K
\]

\vspace{-12pt}\vspace{3pt}
is\dss Fredholm,\oss then\sss
$A\off =\off T^{\dff *}\trf|\trf \kernel\fff \gamma_0$\sss
is\dss an\sss unbounded self-adjoint\trs Fredholm\dss operator in\sss $H_{\dff 0}$\nsp.\oss
Moreover{},\oss it\sss has discrete spectrum and\dss is\dss an operator with compact\sss
resolvent.\oss
See\qss \cite{i2},\oss Theorem\qss 5.4.\oss
All\sss these assumptions hold when\sss the\dss Hilbert\sss spaces involved
are\dss Sobolev\dss spaces,\pss $T$\sss is\dss a differential\sss operator{},\oss
and\sss the boundary condition\sss $\gamma_0\off =\off 0$ satisfies\sss the\dss
Sha\-pi\-ro--Lopatinskii\dss condition\sss for\sss $T$\dnsp,\oss
as\sss it\sss will\dss be\sss the case in\dss 
Sections\qss \ref{differential-boundary-problems}\dss --\dss \ref{comparing}.\oss

Under\sss the assumptions of\dss the previous paragraph,\oss
the operator\sss $A$ can serve as\sss the reference operator\sss
if\dss its kernel\dss is\dss equal\dss to $0$\nnsp.\oss
Indeed,\oss since under\sss these assumptions $A$\sss is\dss self-adjoint,\oss
the injectivity of\dss $A$\sss implies\sss that\sss $A$\sss
is\dss an\sss isomorphism\sss
$\mathcal{D}\dff(\trf A\trf)\qff \ttoo\qff H_{\dff 0}$\nsp.\oss

\myuppar{The reduced\dss boundary operator.}
Let\sss us 
denote\sss the inverse of\dss the  isomorphism\sss 
$\Gamma_0\trf|\trf \kernel\fff T^{\dff *}$\sss by\sss
$\bm{\gamma}\trf(\dff 0\dff)\dff \colon\dff
K\fff'\qff \ttoo\qff \kernel\fff T^{\dff *}$\sss
and\sss set\sss
$M\trf(\dff 0\dff)
\off =\off 
\Gamma_1\dff \circ\trf \bm{\gamma}\trf(\dff 0\dff)
\dff \colon\dff
K\fff'\qff \ttoo\qff K\fff'$\dnsp.\oss
The operators\sss $\bm{\gamma}\trf(\dff 0\dff)$ and\sss $M\trf(\dff 0\dff)$
are analogues of\dss the values at\sss $z\off =\off 0$ of\dss the gamma field $\gamma\trf(\trf z\trf)$
and\sss the\dss Weyl\dss function\sss $M\trf(\trf z\trf)$\sss from\sss the\sss theory of\dss
boundary\sss triplets.\oss
Since\sss $\Gamma_1$\sss is\dss continuous in\sss
the graph\sss topology of\dss $\mathcal{D}\dff(\trf T^{\dff *}\dff)$\nnsp,\oss
the operator\sss $M\trf(\dff 0\dff)$\sss is\dss continuous.\oss
The\qss \emph{reduced\dss boundary operator}\pss is\dss
the operator\vspace{3pt}
\[
\quad
\bm{\Gamma}_1
\off =\off
\Gamma_1\qff -\qff M\trf(\dff 0\dff)\dff \circ\dff \Gamma_0
\qff \colon\qff
\mathcal{D}\dff(\trf T^{\dff *}\dff)
\qff \ttoo\qff
K\fff'
\pff.
\]

\vspace{-12pt}\vspace{3pt}
The continuity of\trs $\Gamma_0\dff,\pff \Gamma_1$\sss and\sss
$M\trf(\dff 0\dff)$\sss implies\sss that\sss
$\bm{\Gamma}_1$\sss is\dss continuous.\oss

\mypar{Lemma.}{reduced-gamma}
$\bm{\Gamma}_1
\off =\off
\Gamma_1\dff \circ\dff p$\nnsp,\oss
\emph{where $p$\sss is\dss the projection\sss
$\mathcal{D}\dff(\trf T^{\dff *}\dff)
\qff \ttoo\qff
\mathcal{D}\dff(\trf A\dff)$\nnsp,\oss
and\qss
$\image\fff \bm{\Gamma}_1\qff \subset\qff K$\nnsp.\oss}

\proof
By\trs Lemma\qss \ref{basic-decomp},\oss 
if\dss $u\qff \in\qff \mathcal{D}\trf(\trf T^{\dff *}\trf)$\nnsp,\oss
then\sss $u\off =\off p\trf(\trf u\trf)\qff +\qff z$\sss
for some $z\qff \in\qff \kernel\fff T^{\dff *}$\dnsp.\oss
Since\sss 
$\mathcal{D}\dff(\trf A\dff)
\off =\off 
\kernel\fff \Gamma_0 \qff|\qff H_{\dff 1}$\sss
and\sss $p$\sss is\dss the projection\sss to $\mathcal{D}\dff(\trf A\dff)$\nnsp,\oss
we see\sss that\sss $\Gamma_0\dff(\trf p\trf(\trf u\trf)\trf)\off =\off 0$\nnsp.\oss
Hence\vspace{4.5pt}
\[
\quad
\bm{\Gamma}_1\dff(\trf u\trf)
\off =\off
\Gamma_1\dff(\trf p\trf(\trf u\trf)\trf)
\qff +\qff
\Gamma_1\dff(\trf z\trf)
\qff -\qff 
M\trf(\dff 0\dff)\dff \circ\dff \Gamma_0\dff(\trf z\trf)
\off =\off
\Gamma_1\dff(\trf p\trf(\trf u\trf)\trf)
\]

\vspace{-12pt}\vspace{4.5pt}
because\sss
$M\trf(\dff 0\dff)\dff \circ\dff \Gamma_0\dff(\trf z\trf)
\off =\off
\Gamma_1\dff(\trf z\trf)$\sss
by\sss the definition.\oss
This proves\sss that\sss
$\bm{\Gamma}_1
\off =\off
\Gamma_1\dff \circ\dff p$\nnsp.\oss
Since\sss $p$\sss maps\sss $\mathcal{D}\dff(\trf T^{\dff *}\dff)$\sss
into\sss $\mathcal{D}\dff(\trf A\dff)\qff \subset\qff H_{\dff 1}$
and\sss $\Gamma_1$\sss maps $H_{\dff 1}$\sss to $K$\nnsp,\oss
this implies\sss that\sss
$\image\fff \bm{\Gamma}_1\qff \subset\qff K$\nnsp.\oss  \eproof

\mypar{Lemma.}{toward-rli}
\emph{If\dss $u\qff \in\qff \mathcal{D}\dff(\trf T^{\dff *}\dff)$\sss
and\sss
$v\qff \in\qff \kernel\fff T^{\dff *}$\dnsp,\oss
then}\vspace{3pt}
\[
\quad
\sco{\dff T^{\dff *} u\fff,\qff v \dff}
\off =\off
\bsco{\dff \bm{\Gamma}_1\dff u\dff,\qff 
\Gamma_0\dff v \dff}_{\dff K\fff,\qff K\fff'}
\pff.
\]

\vspace{-12pt}\vspace{3pt}
\proof
Since\sss $T^{\dff *} v\off =\off 0$\nnsp,\oss\vspace{1.5pt}
\[
\quad
\sco{\dff T^{\dff *} u\fff,\qff v \dff}
\off =\off
\sco{\dff T^{\dff *} p\trf(\trf u\trf)\fff,\qff v \dff}
\off =\off
\sco{\dff T^{\dff *} p\trf(\trf u\trf)\fff,\qff v \dff}
\pff -\pff
\sco{\dff p\trf(\trf u\trf)\fff,\qff T^{\dff *} v \dff}
\pff.
\]

\vspace{-12pt}\vspace{3pt}
Since\sss 
$p\trf(\trf u\trf)
\qff \in\qff 
\mathcal{D}\dff(\trf A\dff)
\qff \subset\qff
H_{\dff 1}$\nsp,\oss
the\dss Lagrange\dss identity\qss (\ref{lagrange-asymm})\qss
applies with $p\trf(\trf u\trf)$ in\sss the role of\dss $u$\nnsp,\oss
and\dss hence\sss the\sss last\sss expression\dss
is\dss equal\dss to\vspace{3pt}
\[
\quad
\bsco{\dff \Gamma_1\dff p\trf(\trf u\trf)\dff,\qff 
\Gamma_0\dff v \dff}_{\dff K\fff,\qff K\fff'}
\pff -\pff
\bsco{\dff \Gamma_0\dff p\trf(\trf u\trf)\dff,\qff 
\Gamma_1\dff v \dff}_{\dff K\fff,\qff K\fff'}
\pff.
\]

\vspace{-12pt}\vspace{3pt}
Since\sss $\Gamma_0\dff p\trf(\trf u\trf)\off =\off 0$\nnsp,\oss
it\dss is\dss also equal\dss to\sss\vspace{3pt}
\[
\quad
\bsco{\dff \Gamma_1\dff p\trf(\trf u\trf)\dff,\qff 
\Gamma_0\dff v \dff}_{\dff K\fff,\qff K\fff'}
\off =\off
\bsco{\dff \bm{\Gamma}_1\dff u\dff,\qff 
\Gamma_0\dff v \dff}_{\dff K\fff,\qff K\fff'}
\pff.
\]

\vspace{-12pt}\vspace{3pt}  
The\sss lemma\sss follows.\oss  \eproof

\mypar{Lemma.}{reduced-l}
\emph{Suppose\sss that\sss 
$u\fff,\qff v\qff \in\qff \mathcal{D}\dff(\trf T^{\dff *}\dff)$\nnsp.\oss
Then}\vspace{3pt}
\[
\quad
\sco{\dff T^{\dff *} u\fff,\qff v \dff}
\pff -\pff
\sco{\dff u\dff,\qff T^{\dff *} v \dff}
\off =\off
\bsco{\dff \bm{\Gamma}_1\dff u\dff,\qff 
\Gamma_0\dff v \dff}_{\dff K\fff,\qff K\fff'}
\pff -\pff
\bsco{\dff \Gamma_0\dff u\dff,\qff 
\bm{\Gamma}_1\dff v \dff}_{\dff K\fff'\fff,\qff K}
\pff
\]

\vspace{-12pt}\vspace{3pt}
\proof
Since\sss
$k\trf(\trf u\trf)\fff,\qff k\trf(\trf v\trf)
\qff \in\qff
\kernel\fff T^{\dff *}$\dnsp,\oss
Lemma\qss \ref{toward-rli}\qss
implies\sss that\vspace{3pt}
\[
\quad
\sco{\dff T^{\dff *} u\fff,\qff k\trf(\trf v\trf) \dff}
\off =\off
\bsco{\dff \bm{\Gamma}_1\dff u\dff,\qff 
\Gamma_0\qff k\trf(\trf v\trf) \dff}_{\dff K\fff,\qff K\fff'}
\hspace{0.8em}
\mbox{and}\hspace{0.8em}
\bsco{\dff k\trf(\trf u\trf)\dff,\qff 
T^{\dff *}  v \dff}
\off =\off
\bsco{\dff \Gamma_0\qff k\trf(\trf u\trf)\dff,\qff 
\bm{\Gamma}_1\dff v \dff}_{\dff K\fff'\fff,\qff K}
\qff.
\]

\vspace{-12pt}\vspace{3pt}
Also,\qss $\Gamma_0\dff p\trf(\trf u\trf)\off =\off 0$\sss
implies\sss that\sss\vspace{3pt}
\[
\quad
\Gamma_0\dff u
\off =\off 
\Gamma_0\qff k\trf(\trf u\trf)
\pff.
\]

\vspace{-12pt}\vspace{3pt}
Similarly,\pss
$\Gamma_0\dff v\off =\off \Gamma_0\qff k\trf(\trf v\trf)$\nnsp.\oss
Hence\vspace{3pt}
\[
\quad
\bsco{\dff T^{\dff *} u\dff,\qff 
k\trf(\trf v\trf) \dff}
\pff -\pff
\bsco{\dff k\trf(\trf u\trf)\dff,\qff 
T^{\dff *}  v \dff}
\off =\off
\bsco{\dff \bm{\Gamma}_1\dff u\dff,\qff 
\Gamma_0\dff v \dff}_{\dff K\fff,\qff K\fff'}
\pff -\pff
\bsco{\dff \Gamma_0\dff u\dff,\qff 
\bm{\Gamma}_1\dff v \dff}_{\dff K\fff'\fff,\qff K}
\pff.
\]

\vspace{-12pt}\vspace{3pt}
It\sss remains\sss to combine\sss this equality\sss with\trs Lemma\qss \ref{modified-l}.\oss  \eproof

\mypar{Lemma.}{surjectivity}
\emph{The map\dss
$\Gamma_0\dff \oplus\dff \bm{\Gamma}_1
\dff \colon\dff
\mathcal{D}\dff(\trf T^{\dff *}\dff)
\qff \ttoo\qff
K\fff'\dff \oplus\dff K$\sss
is\dss surjective.\oss}

\proof
Recall\dss that\sss the restriction\dss 
$\Gamma_0\qff|\qff \kernel\fff T^{\dff *}$\sss is\dss an\sss isomorphism\sss
$\kernel\fff T^{\dff *}\qff \ttoo\qff K\fff'$\dnsp.\oss
Therefore in order\sss to prove surjectivity,\oss
it\dss is\dss sufficient\sss to prove\sss that\sss $\bm{\Gamma}_1$
maps\sss $\kernel\fff \Gamma_0$ onto\sss $K$\nnsp.\oss

By our assumptions,\oss the map\sss
$\Gamma_0\dff \oplus\dff \Gamma_1$ restricted\dss to\sss $H_{\dff 1}$\sss
is\dss surjective onto\sss $K\dff \oplus\dff K$\nnsp.\oss
Therefore\sss $\Gamma_1$\sss maps\sss
$H_{\dff 1}\qff \cap\qff \kernel\fff \Gamma_0$ onto $K$\nnsp.\oss
If\dss $u\qff \in\qff H_{\dff 1}\qff \cap\qff \kernel\fff \Gamma_0$\nnsp,\oss
then\sss $u\qff \in\qff \mathcal{D}\dff(\trf A\dff)$ 
and\dss hence $p\trf(\trf u\trf)\off =\off u$\nnsp.\oss
It\sss follows\sss that\sss
$\bm{\Gamma}_1\dff(\trf u\trf)
\off =\off
\Gamma_1\dff \circ\dff p\dff(\trf u\trf)
\off =\off
\Gamma_1\dff(\trf u\trf)$\nnsp.\oss
This\sss implies\sss that\sss $\bm{\Gamma}_1$\sss maps\sss
$H_{\dff 1}\qff \cap\qff \kernel\fff \Gamma_0$ onto $K$\nnsp.\oss
Together with\sss the previous paragraph\sss this proves surjectivity.\oss  \eproof

\myuppar{The reduced\dss boundary\sss triplet.}
Let\vspace{3pt}
\[
\quad
\overline{\Gamma}_0
\off =\off 
\Lambda'\dff \circ\trf \Gamma_0
\qff \colon\qff
\mathcal{D}\dff(\trf T^{\dff *}\dff)
\qff \ttoo\qff
K^{\dff \partial}
\quad
\mbox{and}\quad 
\]

\vspace{-33pt}
\[
\quad
\overline{\bm{\Gamma}}_1
\off =\off 
\Lambda^{\fff -\dff 1}\dff \circ\trf \bm{\Gamma}_1
\qff \colon\qff
\mathcal{D}\dff(\trf T^{\dff *}\dff)
\qff \ttoo\qff
K^{\dff \partial}
\pff,
\]

\vspace{-12pt}\vspace{3pt}
where\sss $\Lambda\fff,\off \Lambda'$ are\sss the operators associated\sss
with\sss the\dss Gelfand\trs triple\sss
$K\off \subset\off K^{\dff \partial}\off \subset\off K\fff'$\dnsp.\oss
Using\trs Lemma\qss \ref{shift}\qss we can\sss rewrite\sss the\dss
Lagrange\dss identity of\qss Lemma\qss \ref{reduced-l}\qss as\vspace{3pt}
\[
\quad
\bsco{\dff T^{\dff *} u\fff,\qff v \dff}
\pff -\pff
\bsco{\dff u\dff,\qff T^{\dff *} v \dff}
\off =\off
\bsco{\dff \overline{\bm{\Gamma}}_1\dff u\dff,\qff 
\overline{\Gamma}_0\dff v \dff}_{\dff \partial}
\pff -\pff
\bsco{\dff \overline{\Gamma}_0\dff u\dff,\qff 
\overline{\bm{\Gamma}}_1\dff v \dff}_{\dff \partial}
\pff,
\]

\vspace{-12pt}\vspace{3pt}
i.e.\qss in\sss the standard\sss form of\dss the\sss theory of\dss
boundary\sss triplets.\oss
Lemma\qss \ref{surjectivity}\qss implies\sss that\sss\vspace{3pt}
\[
\quad
\overline{\Gamma}_0\qff \oplus\qff \overline{\bm{\Gamma}}_1
\qff \colon\qff
\mathcal{D}\dff(\trf T^{\dff *}\dff)
\qff \ttoo\qff
K^{\dff \partial}\dff \oplus\dff K^{\dff \partial}
\]

\vspace{-12pt}\vspace{3pt}
is\dss surjective.\oss
Therefore\dss\vspace{1.5pt}
\[
\quad
\left(\qff
K^{\dff \partial},\off
\overline{\Gamma}_0\dff,\off
\overline{\bm{\Gamma}}_1
\qff\right)
\]

\vspace{-12pt}\vspace{1.5pt}
is\dss a boundary\sss triplet\dss for\sss $T^{\dff *}$\dnsp,\oss
called\dss the\qss \emph{reduced\dss boundary\dss triplet}.\oss

\mypar{Lemma.}{kernels}
$\kernel\fff \bm{\Gamma}_1
\off =\off
\mathcal{D}\dff(\trf T \trf)
\qff \dotplus\qff
\kernel\fff T^{\dff *}$
\emph{and}\pss
$\kernel\fff \Gamma_0\dff \oplus\dff \bm{\Gamma}_1
\off =\off 
\mathcal{D}\dff(\trf T \trf)$\nnsp.\oss

\proof
Lemmas\qss \ref{basic-decomp}\qss and\qss \ref{reduced-gamma}\qss imply\sss that\sss
$\kernel\fff \bm{\Gamma}_1$\sss is\dss equal\dss to\sss
the direct\sss sum of\dss
$\kernel\fff T^{\dff *}$
and\sss the kernel\sss of\dss the restriction\sss
$\Gamma_1\trf|\qff\fff \mathcal{D}\dff(\trf A\dff)$\nnsp.\oss
Since\sss 
$\mathcal{D}\dff(\trf A\trf)
\off =\off 
\kernel\fff \Gamma_0$\nsp,\oss
the\sss latter\sss kernel\dss is\dss equal\dss to\sss
$\kernel\fff \Gamma_0\dff \oplus\dff \Gamma_1
\off =\off
\mathcal{D}\dff(\trf T\trf)$\nnsp.\oss
This proves\sss the first\sss claim of\dss the\sss lemma.\oss
Since\sss the restriction\dss 
$\Gamma_0\qff|\qff \kernel\fff T^{\dff *}$\sss is\dss injective,\oss
the second claim\sss follows.\oss  \eproof

\myuppar{Comparing\dss boundary\sss triplets.}
Let\sss us apply\sss to\sss the extension $A$ of\dss $T$ and\dss
$\mu\off =\off i$\sss the construction of\dss boundary\sss triplets
from\dss Section\qss \ref{abstract},\oss
and\dss let\vspace{1.5pt}
\begin{equation}
\label{inner-triplet}
\quad
\left(\qff \mathcal{K}_{\dff -}\dff,\off \Gamma^{\trf 0}\dff,\off \Gamma^{\trf 1} \pff\right)
\end{equation}

\vspace{-12pt}\vspace{1.5pt}
be\sss the resulting\sss boundary\sss triplet.\oss 
Then\sss $\mathcal{D}\dff(\trf A\trf)\off =\off \kernel\fff \Gamma^{\trf 0}$\dss
and\sss 
$\Gamma^{\trf 0}\qff \colon\dff
\mathcal{D}\dff(\trf T^{\dff *}\dff)\qff \ttoo\qff \mathcal{K}_{\dff -}$\sss
is\dss a continuous surjective map.\oss
Therefore\sss $\Gamma^{\trf 0}$\sss induces an\sss isomorphism\sss
between\sss the quotient\sss space\sss
$\mathcal{D}\dff(\trf T^{\dff *}\dff)\dff/\dff\mathcal{D}\dff(\trf A\trf)$
and\sss $\mathcal{K}_{\dff -}$\nsp.\oss
By\sss the same reasons\sss $\overline{\Gamma}_0$\sss
induces an\sss isomorphism\sss between\sss the same quotient\sss space\sss
$\mathcal{D}\dff(\trf T^{\dff *}\dff)\dff/\dff\mathcal{D}\dff(\trf A\trf)$\sss
and\sss $K^{\dff \partial}$\dnsp.\oss
It\dss follows\sss that\sss there\dss is\sss a unique\sss topological\sss isomorphism\sss
$D\dff \colon\dff
K^{\dff \partial}\qff \ttoo\qff \mathcal{K}_{\dff -}$\sss
such\sss that\sss the\sss triangle
\[
\quad
\begin{tikzcd}[column sep=tri, row sep=tri]
&
\mathcal{D}\dff(\trf T^{\dff *}\dff)
\arrow[dl, "\dis \overline{\Gamma}_0"']
\arrow[dr, "\dis \Gamma^{\trf 0}"]
\\
K^{\dff \partial}
\arrow[rr, "\dis D"]
&
&
\mathcal{K}_{\dff -}
\end{tikzcd}
\]

\vspace{-12pt}
is\dss commutative.\oss
Similarly,\oss the kernels of\trs 
$\Gamma^{\trf 0}\dff \oplus\dff \Gamma^{\trf 1}$\dss
and\trs
$\overline{\Gamma}_0\dff \oplus\dff \overline{\bm{\Gamma}}_1$
are equal\dss to 
$\mathcal{D}\dff(\trf T \trf)$\sss
and\dss there\dss is\sss a unique\sss topological\sss isomorphism\sss
$\mathcal{W}\dff \colon\dff
K^{\dff \partial}\dff \oplus\dff K^{\dff \partial}
\dff \qff \ttoo\qff 
\mathcal{K}_{\dff -}\dff \oplus\dff \mathcal{K}_{\dff -}$\sss
such\sss that\sss the\sss triangle
\[
\quad
\begin{tikzcd}[column sep=tri, row sep=tri]
&
\mathcal{D}\dff(\trf T^{\dff *}\dff)
\arrow[dl, "\dis \overline{\Gamma}_0\dff \oplus\dff \overline{\bm{\Gamma}}_1"']
\arrow[dr, "\dis \Gamma^{\trf 0}\dff \oplus\dff \Gamma^{\trf 1}"]
\\
K^{\dff \partial}\dff \oplus\dff K^{\dff \partial}
\arrow[rr, "\dis \mathcal{W}"]
&
&
\mathcal{K}_{\dff -}\dff \oplus\dff \mathcal{K}_{\dff -}
\end{tikzcd}
\]

\vspace{-12pt}
is\dss commutative.\oss
The uniqueness of\dss $D$\sss implies\sss that\sss the square
\[
\quad
\begin{tikzcd}[column sep=tri, row sep=tri]
K^{\dff \partial}\dff \oplus\dff K^{\dff \partial}
\arrow[rr, "\dis \mathcal{W}"]
\arrow[d]
&
&
\mathcal{K}_{\dff -}\dff \oplus\dff \mathcal{K}_{\dff -}
\arrow[d]
\\
K^{\dff \partial}
\arrow[rr, "\dis D"]
&
&
\qff
\mathcal{K}_{\dff -}\dff,
\end{tikzcd}
\]

\vspace{-12pt}
where\sss the vertical\sss arrows are projections on\sss the\sss first\sss summand,\oss
is\dss commutative.\oss
By\sss the uniqueness of\trs boundary\sss triplets\sss 
$\mathcal{W}$\sss is\dss an\sss isometry\sss between\sss the\dss Hermitian\dss
scalar product\vspace{3pt}
\[
\quad
[\trf (\trf u\fff,\qff v\trf)\fff,\qff (\trf a\fff,\qff b\trf)\trf]_{\dff \partial}
\off =\off
i\trf \sco{\dff u\fff,\qff b\dff}_{\dff \partial}
\qff -\qff
i\trf \sco{\dff v\fff,\qff a\dff}_{\dff \partial}
\qff.
\]

\vspace{-12pt}\vspace{3pt}
on\sss $K^{\dff \partial}\dff \oplus\dff K^{\dff \partial}$\sss
and similarly defined\sss product on\sss
$\mathcal{K}_{\dff -}\dff \oplus\dff \mathcal{K}_{\dff -}$\nnsp.\oss
See\qss \cite{bhs},\oss Section\qss 1.8\qss and\trs Theorem\qss 2.5.1.\oss
Together\sss with\sss the fact\sss that\sss
$\kernel\fff \Gamma^{\trf 0}
\off =\off 
\kernel\fff \overline{\Gamma}_0$\sss
this implies\sss that\sss there exists a bounded self-adjoint\sss operator\sss
$P\dff \colon\dff
\mathcal{K}_{\dff -}\qff \ttoo\qff \mathcal{K}_{\dff -}$\sss
such\sss that\vspace{4.5pt}
\begin{equation}
\label{compare}
\quad
\overline{\Gamma}_0
\off =\off
D^{\dff -\dff 1}\trf \Gamma^{\trf 0}
\quad
\mbox{and}\quad
\overline{\bm{\Gamma}}_1
\off =\off
D^{\dff *}\trf \Gamma^{\trf 1}
\qff +\qff
P\trf D^{\dff -\dff 1}\trf \Gamma^{\trf 0}
\pff.
\end{equation}

\vspace{-12pt}\vspace{4.5pt}
See\qss \cite{bhs},\oss Corollary\qss 2.5.6.\oss

\myuppar{Two examples.}
The first\sss one\dss is\dss the reference operator $A$\sss itself.\oss
By our assumptions\sss it\dss is\dss the self-adjoint\sss extension of\dss $T$\sss
defined\sss by either\sss the boundary condition\sss $\gamma_0\off =\off 0$\nnsp,\oss
or\sss by\sss the boundary condition\sss $\Gamma_0\off =\off 0$\nnsp.\oss 
The\sss latter\dss is\dss obviously equivalent\sss to\sss $\overline{\Gamma}_0\off =\off 0$\nnsp.\oss
Therefore $A$ can\sss be defined\sss in\sss terms of\dss the reduced\dss boundary\dss triplet\sss
as\sss the extension defined\dss by\sss the self-adjoint\sss relation\sss
$0\dff \oplus\dff K^{\dff \partial}
\qff \subset\qff 
K^{\dff \partial}\dff \oplus\dff K^{\dff \partial}$\dnsp.\oss
In\sss the notations of\pss \cite{s},\oss $A$\sss is\dss the extension\sss
$T_{\dff 0}$\nsp.\oss

The second example\dss is\sss  
$A'\off =\off T^{\dff *}\trf|\trf \kernel\fff \gamma_1$\nsp.\oss 
\emph{Suppose\sss that\dss $A'$
is\dss self-adjoint.}\oss 
By\trs Lemma\qss \ref{gamma-0-ker}\qss with\sss the roles of\dss
$\gamma_0$ and\sss $\gamma_1$\sss interchanged,\oss
we see\sss that\sss
$\kernel\fff \Gamma_1
\off =\off
\kernel\fff \gamma_1
\qff \subset\qff 
H_{\dff 1}$\nsp.\oss
Clearly,\vspace{3pt}
\[
\quad
\kernel\fff \Gamma_1
\off =\off
\kernel\dff \bigl(\qff
\bm{\Gamma}_1\qff +\qff M\trf(\dff 0\dff)\dff \circ\dff \Gamma_0
\qff\bigr)
\pff.
\]

\vspace{-12pt}\vspace{3pt}
If\dss $u\qff \in\qff \kernel\fff \Gamma_1$\nsp,\oss 
then\sss $u\qff \in\qff H_{\dff 1}$\sss
and\dss hence\sss 
$\Gamma_0\dff u\off =\off \gamma_0\dff u\qff \in\qff K$\nnsp.\oss
Therefore we can\sss replace\sss $M\trf(\dff 0\dff)$\sss 
by\sss $M\off =\off M\trf(\dff 0\dff)\trf|\trf K$
considered as a densely defined operator\sss in\sss $K\fff'$\dnsp.\oss
More precisely,\oss\vspace{3pt}
\[
\quad
\kernel\fff \Gamma_1
\off =\off
\kernel\dff \bigl(\qff
\bm{\Gamma}_1\qff +\qff M\dff \circ\dff \Gamma_0
\qff\bigr)
\pff.
\]

\vspace{-12pt}\vspace{3pt}
\emph{Let\sss us assume\sss that\sss the operator\sss 
$M\trf(\dff 0\dff)\dff \colon\dff
K\fff'\qff \ttoo\qff K\fff'
$\sss
leaves $K$\sss invariant\sss and\dss that\sss the induced operator\sss
$K\qff \ttoo\qff K$\sss is\trs Fredholm.\oss}
Then\sss $M\off =\off M\trf(\dff 0\dff)\trf|\trf K$\sss is\dss a closed
densely defined operator\sss from\sss $K\fff'$\sss to\sss $K$\dss
({\fff}but\sss $M$\dss is\dss not\sss closed as an operator\sss
from\sss $K\fff'$\sss to\sss $K\fff'$\nsp).\oss
Clearly,\vspace{4.5pt}
\[
\quad
\bm{\Gamma}_1\qff +\qff M\dff \circ\dff \Gamma_0
\off =\off
\Lambda\dff \circ\dff \overline{\bm{\Gamma}}_1
\qff +\qff 
M\dff \circ\dff (\trf \Lambda'\qff)^{\dff -\dff 1}\dff \circ\dff \overline{\Gamma}_0
\pff.
\]

\vspace{-12pt}\vspace{4.5pt}
Since $\Lambda$\sss is\dss an\sss isomorphism,\oss 
it\sss follows\sss that\vspace{3pt}
\[
\quad
\kernel\dff \bigl(\qff
\bm{\Gamma}_1\qff +\qff M\dff \circ\dff \Gamma_0
\qff\bigr)
\off =\off
\kernel\dff \bigl(\qff
\overline{\bm{\Gamma}}_1\qff +\qff \overline{M}\dff \circ\dff \overline{\Gamma}_0
\qff\bigr)
\pff,
\]

\vspace{-12pt}\vspace{3pt}
where\dss
$\overline{M}
\off =\off
\Lambda^{\fff -\dff 1}\dff \circ\dff
M\dff \circ\dff 
(\trf \Lambda'\qff)^{\dff -\dff 1}$\sss
is\dss a closed densely defined operator\sss in\sss $K^{\dff \partial}$\dnsp.\oss
We see\sss that\sss $A'$\sss can\sss be defined\sss in\sss terms of\dss the
reduced\dss boundary\sss triplet\sss as\sss the extension of\dss $T$\sss
corresponding\sss to\sss the operator\sss $-\qff \overline{M}$\nnsp,\oss
or{},\oss rather{},\oss its\sss graph considered as a relation\sss in\sss
$K^{\dff \partial}$\dnsp.\oss
Hence,\oss by\sss the\sss theory of\dss boundary\sss triplets,\oss 
the self-adjointness of\dss $A'$\sss implies\sss that\sss
$\overline{M}$\sss is\dss a self-adjoint\sss operator\sss
in\sss $K^{\dff \partial}$\dnsp.\oss
By\sss the discussion at\sss the end of\trs Section\qss \ref{gelfand-triples}\qss
the operator\sss $\overline{M}$\sss is\dss self-adjoint\sss if\dss
and\dss only\trs if\dss the operator\sss $M$\sss is\dss self-adjoint\sss
as an operator\sss from\sss $K\fff'$\sss to $K$\nnsp.\oss

\myuppar{Remark.}
If\dss the inclusion\sss $K\qff \ttoo\qff K\fff'$\sss is\dss a compact\sss operator,\oss
then\sss $\overline{M}$\sss is\dss an operator with compact\sss resolvent.\oss
Indeed,\oss under our assumptions $M$\sss is\trs Fredholm\dss and\dss hence\sss
$M\qff -\qff \lambda$\sss is\dss an\sss isomorphism\sss
$K\qff \ttoo\qff K$\sss for some $\lambda\qff \in\qff \rrr$\nnsp.\oss 
The inverse of\qss $\overline{M}\qff -\qff \lambda$\sss is\sss
$\Lambda'\dff \circ\dff (\trf M\qff -\qff \lambda\trf)^{\dff -\dff 1}\dff \circ\dff \Lambda$\nnsp.\oss
Hence\sss $(\trf M\qff -\qff \lambda\trf)^{\dff -\dff 1}$\sss 
is\dss the composition of\dss a\sss topological\sss isomorphism\sss
$K\qff \ttoo\qff K$\sss and\sss the compact\sss inclusion\sss $K\qff \ttoo\qff K\fff'$\dnsp.\oss
Since $\Lambda$ and\sss $\Lambda'$ are\sss topological\sss isomorphisms,\oss
this\sss implies\sss that\sss the inverse of\qss $\overline{M}\qff -\qff \lambda$\sss is\dss compact.\oss

\myuppar{Boundary problems defined\sss in\sss terms of\qss $\gamma_0\dff,\pff \gamma_1$\nsp.}
Let\sss $\gamma\off =\off \gamma_0\dff \oplus\dff \gamma_1$\nsp.\oss
Let\sss 
$\mathcal{B}\off \subset\off K^{\dff \partial}\dff \oplus\dff K^{\dff \partial}$\sss 
be a closed\sss relation and\dss
let $T_{\fff \mathcal{B}}$\sss be\sss the restriction of\dss $T^{\dff *}${\nsp}  to\sss
$\gamma^{\dff -\dff 1}\dff(\trf \mathcal{B}\trf)\qff \subset\qff H_{\dff 1}$\nsp.\oss
Recall\dss that\sss 
$H_{\dff 1}\qff \subset\qff \mathcal{D}\dff(\trf T^{\dff *}\trf)$\nnsp.\oss 
\emph{Suppose\sss that\sss $T_{\fff \mathcal{B}}$ is\dss a self-adjoint\sss operator in\sss $H_{\trf 0}$\nsp.}\oss
Then\sss $T_{\fff \mathcal{B}}$\sss can\sss be defined\sss in\sss terms of\dss
the reduced\dss boundary\sss triplet.\oss
The corresponding\sss boundary condition\dss is\vspace{3pt}
\[
\quad
\overline{\Gamma}_0\qff \oplus\qff \overline{\bm{\Gamma}}_1\qff
\bigl(\qff \mathcal{D}\dff(\trf T_{\fff \mathcal{B}}\trf) \qff\bigr)
\off =\off
\overline{\Gamma}_0\qff \oplus\qff \overline{\bm{\Gamma}}_1\qff
\bigl(\qff \gamma^{\dff -\dff 1}\dff(\trf \mathcal{B}\trf) \qff\bigr)
\pff.
\]

\vspace{-12pt}\vspace{3pt}
It\dss is\dss equal\sss to\sss the image under\sss the map\sss
$\Lambda'\dff \oplus\dff \Lambda^{\fff -\dff 1}$\sss 
of\dss\vspace{3pt}
\[
\quad
\Gamma_0\dff \oplus\dff \bm{\Gamma}_1\qff
\bigl(\qff \gamma^{\dff -\dff 1}\dff(\trf \mathcal{B}\trf) \qff\bigr)
\off \subset\off
K\fff'\dff \oplus\dff K
\pff.
\]

\vspace{-12pt}\vspace{3pt}
Let\sss $\mathcal{B}\trf|\trf K
\off =\off
\mathcal{B}\qff \cap\qff K\dff \oplus\dff K$\sss
be\sss the restriction of\dss $\mathcal{B}$\sss to $K$\nnsp.\oss
Since\sss 
$\gamma^{\dff -\dff 1}\dff(\trf \mathcal{B}\trf)\qff \subset\qff H_{\dff 1}$\nsp,\oss
\vspace{3pt}
\[
\quad
\Gamma_0\dff \oplus\dff \bm{\Gamma}_1\qff
\bigl(\qff \gamma^{\dff -\dff 1}\dff(\trf \mathcal{B}\trf) \qff\bigr)
\off =\off
\gamma_0\dff \oplus\dff \bm{\gamma}_1\qff
\bigl(\qff \gamma^{\dff -\dff 1}\dff(\trf \mathcal{B}\trf|\trf K\trf) \qff\bigr)
\pff,
\]

\vspace{-12pt}\vspace{3pt}
where\sss 
$\bm{\gamma}_1
\off =\off
\gamma_1\qff -\qff M\dff \circ\dff \gamma_0$\sss
and\sss $M\off =\off M\trf(\dff 0\dff)\trf|\trf K$\sss as above.\oss
Since $\gamma$\sss is\dss a map onto $K\dff \oplus\dff K$\nnsp,\oss\vspace{3pt}
\[
\quad
\gamma_0\dff \oplus\dff \bm{\gamma}_1\qff
\bigl(\qff \gamma^{\dff -\dff 1}\dff(\trf \mathcal{B}\trf|\trf K\trf) \qff\bigr)
\off =\off
\left\{\pff
(\trf u\fff,\qff v\qff -\qff M\dff u\trf)
\pff\bigl|\pff
(\trf u\fff,\qff v\trf)\qff \in\qff \mathcal{B}\trf|\trf K
\pff\right\}
\off =\off
\mathcal{B}\trf|\trf K\qff -\qff M
\qff.
\]

\vspace{-12pt}\vspace{3pt}
Therefore in\sss terms of\dss the reduced\dss boundary\sss triplet\sss
$T_{\fff \mathcal{B}}$\sss is\dss defined\dss by\sss the boundary condition\vspace{3pt}
\[
\quad
\Lambda'\dff \oplus\dff \Lambda^{\fff -\dff 1}\qff
\bigl(\qff
\mathcal{B}\trf|\trf K\qff -\qff M
\qff\bigr)
\pff.
\]

\vspace{-12pt}\vspace{3pt}
If\dss $\mathcal{B}\trf|\trf K$\sss is\dss the graph of\dss an operator\sss $B_{\dff K}$\nsp,\oss
then\sss this boundary condition\dss is\dss the graph of\dss\vspace{3pt}\vspace{-0.25pt}
\[
\quad
\Lambda^{\fff -\dff 1}\qff \circ\qff 
\bigl(\qff
B_{\dff K}\qff -\qff M
\qff\bigr)\qff \circ\qff 
(\trf \Lambda'\trf)^{\fff -\dff 1}
\pff.
\]

\vspace{-12pt}\vspace{3pt}\vspace{-0.25pt}
In our applications\sss $\mathcal{B}\trf|\trf K$\sss will\dss be usually a proper\sss
relation,\oss not\sss a\sss graph.\oss
Note\sss that\sss the assumption\sss of\dss self-adjointness of\dss $T_{\fff \mathcal{B}}$\sss
implies,\oss by\sss the\sss theory of\dss boundary\sss triplets and\dss the remarks
at\sss the end of\qss Section\qss \ref{gelfand-triples},\oss
that\sss the relation\sss $\mathcal{B}\trf|\trf K\qff -\qff M$\sss
is\dss self-adjoint.\oss
In applications,\oss the self-adjointness of\dss $T_{\fff \mathcal{B}}$\sss
is\dss established\sss in\sss the same way as\sss
the self-adjointness of\dss $A$\nnsp.\oss

\newpage
\mysection{Families\qss of\pss abstract\qss boundary\qss problems}{families}

\myuppar{Families of\dss extensions.}
Let\sss $W$\sss be a\sss reasonable\qss 
(say,\oss compactly\sss generated and\sss paracompact)\qss 
topological\sss space.\oss
Let $H$ be a separable\dss Hilbert\dss space,\pss
$T_w\dff,\pff w\qff \in\qff W$\sss be a family of\dss 
densely defined closed symmetric operators in $H$\nnsp,\oss
and $A_{\dff w}\dff,\pff w\qff \in\qff W$\sss be a family of\dss self-adjoint\sss operators
such\sss that\sss $T_w\qff \subset\qff A_{\dff w}\qff \subset\qff T_w^{\dff *}$\sss
for every\sss $w\qff \in\qff W$\nnsp.\oss
We will\sss assume\sss 
that\sss the family\sss $A_{\dff w}\dff,\pff w\qff \in\qff W$\sss
is\dss continuous in\sss the\sss topology of\dss uniform\sss resolvent\sss convergence.\oss
Let\vspace{3pt}
\[
\quad
\mathcal{K}_{\dff w\dff +}
\off =\off
\kernel\fff (\trf T_w^{\dff *}\qff -\qff\halfff i\qff)
\off =\off 
\image\fff (\trf T_w\qff +\qff i\qff)^{\dff \perp}
\quad
\mbox{and}\quad
\]

\vspace{-33pt}
\[
\quad
\mathcal{K}_{\dff w\dff -}
\off =\off
\kernel\fff (\trf T_w^{\dff *}\qff +\qff i\qff)
\off =\off 
\image\fff (\trf T_w\qff -\qff i\qff)^{\dff \perp}
\pff.
\]

\vspace{-12pt}\vspace{3pt}
Let\sss 
$V_w\dff \colon\dff
\mathcal{K}_{\dff w\dff +}\qff \ttoo\qff \mathcal{K}_{\dff w\dff -}$\sss
be\sss the isometry corresponding\sss to $A\off =\off A_{\dff w}$
and\sss $\mu\off =\off i$\sss as
in\dss Section\qss \ref{abstract}.\oss
By\trs Lemma\qss \ref{unitary-op}\qss the isometry\sss $V_w$\sss
is\dss equal\sss to\sss the restriction of\dss  
$U\trf(\trf A_{\dff w}\trf)$\sss
to\sss $\mathcal{K}_{\dff w\dff +}$\nsp.\oss
Similarly,\oss the constructions of\trs Section\qss \ref{abstract}\qss
lead\dss to boundary operators\vspace{3pt}
\[
\quad
\Gamma_{w\trf 0}\dff,\off \Gamma_{w\trf 1}
\qff \colon\qff
\mathcal{D}\trf(\trf T_w^{\dff *}\trf)
\qff \ttoo\qff
\mathcal{K}_{\dff w\dff -}
\pff,\quad
w\qff \in\qff W
\pff
\]

\vspace{-12pt}\vspace{3pt}
such\sss that\sss the analogue of\dss the\dss 
Lagrange\dss identity\qss (\ref{l-von-N})\qss with subscripts $w$ holds 
for every $w$\nnsp.\oss

Let\sss $\mathcal{B}_w\dff,\pff w\qff \in\qff W$\sss be a\sss family\sss
of\dss self-adjoint\sss relations\sss
$\mathcal{B}_w
\qff \subset\qff 
\mathcal{K}_{\dff w\dff -}
\dff \oplus\dff 
\mathcal{K}_{\dff w\dff -}$\nsp.\oss
We will\dss impose a continuity assumption on\sss this\sss family\sss
a\sss little bit\dss later.\oss 
The boundary\sss triplets\vspace{3pt}
\[
\quad
\bigl(\qff \mathcal{K}_{\dff w\dff -}\dff,\off \Gamma_{w\trf 0}\dff,\off \Gamma_{w\trf 1} \pff\bigr)
\]

\vspace{-12pt}\vspace{3pt}
together with\sss relations $\mathcal{B}_w$ 
define a new\sss family\sss $A'_{\dff w}\dff,\pff w\qff \in\qff W$\sss
of\dss extensions of\dss the operators $T_w\dff,\pff w\qff \in\qff W$\sss
such\sss that\sss $T_w\qff \subset\qff A'_{\dff w}\qff \subset\qff T_w^{\dff *}$\sss
for every\sss $w\qff \in\qff W$\nnsp.\oss
As we saw\sss in\dss Section\qss \ref{abstract},\oss\vspace{3pt}
\begin{equation}
\label{cayley-t-parameter}
\quad
U\trf(\trf A'_{\dff w}\trf)
\off =\off
U\trf(\trf \mathcal{B}_w\trf)_{\dff H}\trf
U\trf(\trf A_{\dff w}\trf)
\pff
\end{equation}

\vspace{-12pt}\vspace{3pt}
for every\sss $w\qff \in\qff W$\nnsp.\oss
The continuity of\dss the family\sss
$A_{\dff w}\dff,\pff w\qff \in\qff W$\sss
in\sss the\sss topology of\dss uniform\sss resolvent\sss convergence\dss
is\dss equivalent\sss to\sss the continuity of\dss
the family\sss $U\trf(\trf A_{\dff w}\trf)\dff,\pff w\qff \in\qff W$\sss
in\sss the norm\sss topology.\oss
Therefore\sss the family\sss
$U\trf(\trf A'_{\dff w}\trf)\dff,\pff w\qff \in\qff W$\sss
is\dss continuous in\sss the\sss topology of\dss uniform\sss resolvent\sss convergence\dss
if\dss and\dss only\trs if\dss the family\sss
$U\trf(\trf \mathcal{B}_w\trf)_{\dff H}\dff,\pff w\qff \in\qff W$\sss
is\dss continuous in\sss the norm\sss topology.\oss
The following assumptions ensure\sss such continuity.\oss

\myuppar{The continuity assumptions.}
Let\sss us assume\sss that\sss the subspace\sss
$\mathcal{K}_{\dff w\dff -}$
of\dss $H$ continuously depends on $w$\nnsp.\oss
Since\sss
$\mathcal{K}_{\dff w\dff +}
\off =\off
U\trf(\trf A_{\dff w}\trf)^{\dff -\dff 1}\dff
(\trf \mathcal{K}_{\dff w\dff -}\trf)$\nnsp,\oss
under\sss this assumption\sss
$\mathcal{K}_{\dff w\dff +}$\sss
also continuously depends on $w$\nnsp.\oss
Since\sss $\mathcal{K}_{\dff w\dff -}$
continuously depends on $w$\nnsp,\oss
it\dss makes sense\sss to speak about\sss the continuity of\dss the family\sss
$\mathcal{B}_w\dff,\pff w\qff \in\qff W$\nnsp,\oss
and we will\sss assume\sss that\sss it\dss is\dss continuous.\oss
These continuity assumptions imply\sss that\sss the family\sss
$U\trf(\trf \mathcal{B}_w\trf)_{\dff H}\dff,\pff w\qff \in\qff W$\sss
is\dss continuous\sss in\sss the norm\sss topology,\oss
and\dss hence\sss the family\sss
$U\trf(\trf A'_{\dff w}\trf)\dff,\pff w\qff \in\qff W$\sss
is\dss also continuous.\oss
Hence\sss the family\sss
$A'_{\dff w}\dff,\pff w\qff \in\qff W$\sss
is\dss continuous in\sss
in\sss the\sss topology of\dss uniform\sss resolvent\sss convergence.\oss

\myuppar{Fredholm\dss families of\dss operators and\sss relations.}
Suppose now\sss that\sss
$A_{\dff w}\dff,\pff w\qff \in\qff W$\sss
is\dss a\dss family of\trs Fredholm\dss operators.\oss
Since\sss this family\dss is\dss continuous\sss
in\sss the\sss topology of\dss uniform\sss resolvent\sss convergence,\oss
this\sss implies\sss that\sss $A_{\dff w}\dff,\pff w\qff \in\qff W$\sss
is\dss a\dss Fredholm\trs family\sss in\sss the sense of\pss \cite{i1}.\oss
Suppose also\sss that\sss $\mathcal{B}_w\dff,\pff w\qff \in\qff W$\sss
is\dss a\sss family of\trs Fredholm\dss relations in\sss the sense of\pss \cite{i2}.\oss
Then\sss the analytical\dss index\dss is\dss defined\sss for each of\dss
the families\sss $A_{\dff w}\dff,\pff w\qff \in\qff W$\sss
and\sss $\mathcal{B}_w\dff,\pff w\qff \in\qff W$\nnsp.\oss
We will\sss denote\sss these analytical\dss indices by\sss
$\ai\trf(\trf A_{\dff \bullet}\trf)$ and\sss
$\ai\trf(\trf \mathcal{B}_{\fff \bullet}\trf)$\sss respectively.\oss

As in\qss \cite{i2},\oss we denote by\sss $U\fred$\sss the space of\pss \emph{Fredholm-unitary}\pss
operators in $H$\nnsp,\oss i.e.\qss of\dss unitary operators in $H$ such\sss that\sss
$-\qff 1$\sss does not\sss belongs\sss to\sss the essential\sss spectrum.\oss
Equivalently,\oss an operator belongs\sss to $U\fred$\sss if\dss and\dss only\sss if\trs
it\dss is\dss equal\dss to\sss the\dss Cayley\dss transform of\dss a self-adjoint\dss
Fredholm\dss relation in $H$\nnsp.\oss
Unfortunately,\pss $U\fred$\sss is\dss not\sss closed under\sss the composition.\oss
By\sss this reason\qss (\ref{cayley-t-parameter})\qss alone\dss
is\dss not\sss sufficient\sss to conclude\sss that\sss the operators\sss
$A'_{\dff w}$ are\dss Fredholm.

As in\qss \cite{i2},\oss we denote by\sss $U\comp$\sss the group of\dss
unitary operators $H\qff \ttoo\qff H$ which differ\sss from\sss
$\id_{\trf H}$\sss by a compact\sss operator.\oss
It\dss is\dss easy\sss to see\sss that\sss $U\comp$
acts on $U\fred$\sss by composition\sss from either side,\oss
i.e.\qss that\sss $V\qff \in\qff U\comp$\dnsp,\qss $V\fff'\qff \in\qff U\fred$\sss
implies\sss $V\dff \circ\dff V\fff'$ and\sss $V\fff'\dff \circ\dff V$\sss
belong\sss to $U\fred$\dnsp.\oss
Together\sss with\qss (\ref{cayley-t-parameter})\qss
this implies\sss that\dss if\dss
$A_{\dff w}\dff,\pff w\qff \in\qff W$\sss
is\dss a\sss family of\trs Fredholm\dss operators with compact\sss resolvent,\oss
then\sss $A'_{\dff w}\dff,\pff w\qff \in\qff W$\sss is\dss a family of\trs Fredholm\dss operators.\oss
The same conclusion\sss holds\sss if\dss we impose a similar condition on\sss
$\mathcal{B}_w\dff,\pff w\qff \in\qff W$ instead of\dss $A_{\dff w}\dff,\pff w\qff \in\qff W$\nnsp.

In more details,\oss let\sss us say\sss that\sss a self-adjoint\dss relation\sss
$\mathcal{B}\qff \subset\qff H\dff \oplus\dff H$\qss
\emph{has compact\dss resolvent}\pss if\dss its operator part\dss
has compact\sss resolvent.\oss
It\dss is\dss easy\sss to see\sss that\sss $\mathcal{B}$\sss
has compact\sss resolvent\dss if\dss and\dss only\trs if\trs
$U\trf(\trf \mathcal{B}\trf)\qff \in\qff U\comp$\dnsp.\oss
Clearly,\oss if\dss
$\mathcal{B}\qff \subset\qff \mathcal{K}_{\dff -}\dff \oplus\dff \mathcal{K}_{\dff -}$\qss
is\dss a relation with compact\sss resolvent,\oss then 
$U\trf(\trf \mathcal{B}\trf)_{\dff H}$\sss belongs\sss to $U\comp$ of\dss $H$\nnsp.\oss
Therefore\sss if\dss $\mathcal{B}_w\dff,\pff w\qff \in\qff W$\sss
is\dss a\sss family of\trs Fredholm\dss relations with compact\sss resolvent,\oss
then\sss $A'_{\dff w}\dff,\pff w\qff \in\qff W$\sss is\dss a family of\trs Fredholm\dss operators.\oss

\mypar{Theorem.}{index-op}
\emph{Suppose\sss that\sss either\sss every operator\sss
$A_{\dff w}\dff,\pff w\qff \in\qff W$\sss has compact\sss resolvent,\oss
or every\sss relation $\mathcal{B}_w\dff,\pff w\qff \in\qff W$\sss has compact\sss resolvent.\oss
Then}\vspace{3pt}
\[
\quad
\ai\trf(\trf A'_{\dff \bullet}\trf)
\off =\off
\ai\trf(\trf A_{\dff \bullet}\trf)
\pff +\pff
\ai\trf(\trf \mathcal{B}_{\fff \bullet}\trf)
\pff.
\]

\vspace{-12pt}\vspace{3pt}
\proof
The discussion\sss preceding\sss the\sss theorem shows\sss that\sss
$\ai\trf(\trf A'_{\dff \bullet}\trf)$\sss is\dss well-defined.\oss
The analytical\dss index of\dss the family\sss
$\mathcal{B}_w\dff,\pff w\qff \in\qff W$\sss
is\dss equal\dss to\sss the analytical\dss index of\dss the family\sss
$U\trf(\trf \mathcal{B}_w\trf)\dff,\pff w\qff \in\qff W$\sss
of\trs Fredholm-unitary\sss operators.\oss
See\qss \cite{i2},\oss the end of\trs Section\qss 11.\oss
Clearly,\oss the families\sss 
$U\trf(\trf \mathcal{B}_w\trf)\dff,\pff w\qff \in\qff W$\sss
and\sss
$U\trf(\trf \mathcal{B}_w\trf)_{\dff H}\dff,\pff w\qff \in\qff W$\sss
of\trs Fredholm-unitary\sss operators have\sss the same analytical\dss index.\oss
Now\sss the\sss theorem\sss follows\sss from\qss (\ref{cayley-t-parameter})\qss and\qss \cite{i2},\oss
Lemma\qss 11.2.\oss  \eproof

\myuppar{Families of\dss self-adjoint\dss boundary\sss problems.\oss I.}
Let\sss us\sss pass\sss to\sss the framework of\qss Section\qss \ref{abstract-index}\qss
and\sss keep all\sss assumptions of\dss that\sss section.\oss
Let\sss 
$\mathcal{B}_{\fff w}\dff,\pff w\qff \in\qff W$\sss
be a continuous family of\dss self-adjoint\trs Fredholm\dss relations in\sss
$K^{\dff \partial}$\dnsp,\pss
$\mathcal{B}_{\fff w}
\qff \subset\qff
K^{\dff \partial}\dff \oplus\dff K^{\dff \partial}$\dss
for every\sss $w\qff \in\qff W$\nnsp.\oss
Then\sss $\mathcal{B}_{\fff w}\dff,\pff w\qff \in\qff W$\sss
is\dss a\dss Fredholm\dss family\sss
and\dss its analytical\dss index\dss is\dss defined.\oss
See\qss \cite{i2},\oss Section\qss 11.\oss
The\sss family\sss $\mathcal{B}_{\fff w}\dff,\pff w\qff \in\qff W$ 
defines a family of\dss self-adjoint\sss extensions\sss
\vspace{1.5pt}
\[
\quad
\mathcal{T}_w
\off =\off
T_{\fff \mathcal{B}_{\fff w}}\dff,\off w\qff \in\qff W
\]

\vspace{-12pt}\vspace{1.5pt}
of\dss $T$\dnsp,\oss where\sss 
$\mathcal{T}_w
\off =\off
T_{\fff \mathcal{B}_{\fff w}}$\sss
is\dss the restriction of\dss $T^{\dff *}$\sss to\sss the subspace\vspace{3pt}
\[
\quad
\mathcal{D}\dff\left(\qff T_{\fff \mathcal{B}_{\fff w}} \qff\right)
\off =\off
\left\{\qff 
x\qff \in\qff \mathcal{D}\dff(\trf T^{\dff *}\trf)
\pff \left|\off
\left(\qff \overline{\Gamma}_0\dff x\fff,\pff \overline{\bm{\Gamma}}_1\dff x \qff\right)
\qff \in\qff
\mathcal{B}_{\fff w} \right. 
\pff\right\} 
\off.
\]

\vspace{-12pt}\vspace{3pt}
We will\sss say\sss that\sss $\mathcal{T}_w$\sss is\dss defined\dss by\sss the equation\sss
$\overline{\bm{\Gamma}}_1\off =\off \mathcal{B}_{\fff w}\qff \overline{\Gamma}_0$\nsp,\oss
and similarly for other relations and\dss boundary\sss triplets.\oss
The extensions\sss $\mathcal{T}_w$\sss can\sss be also defined\sss in\sss terms of\dss the boundary\sss
triplet\qss (\ref{inner-triplet}).\oss
In\sss view of\pss (\ref{compare})\qss the corresponding equations are\vspace{3pt}
\[
\quad
D^{\dff *}\trf \Gamma^{\trf 1}
\qff +\qff
P\trf D^{\dff -\dff 1}\trf \Gamma^{\trf 0}
\off =\off
\mathcal{B}_{\fff w}\qff
D^{\dff -\dff 1}\trf \Gamma^{\trf 0}
\pff,
\]

\vspace{-12pt}\vspace{3pt}
or{},\oss equivalently,\oss
$\Gamma^{\trf 1}
\off =\off
\mathcal{C}_{\fff w}\qff \Gamma^{\trf 0}$\dnsp,\oss
where\vspace{3pt}
\[
\quad
\mathcal{C}_{\fff w}
\off =\off
\left(\qff D^{\dff *}\qff\right)^{\dff -\dff 1}
\qff
\left(\qff
\mathcal{B}_{\fff w}\qff
D^{\dff -\dff 1}\trf 
\trf -\pff
P\trf D^{\dff -\dff 1}
\qff \right)
\pff.
\]

\vspace{-12pt}\vspace{3pt}
The family\sss $\mathcal{C}_{\fff w}\dff,\pff w\qff \in\qff W$\sss
is\dss continuous\sss together\sss with\sss
$\mathcal{B}_{\fff w}\dff,\pff w\qff \in\qff W$\nnsp.\oss

\mypar{Theorem.}{index-for-const-a}
\emph{Suppose\sss that\sss $\mathcal{B}_{\fff w}$\sss is\dss a relation with compact\dss
resolvent\sss for every\sss $w\qff \in\qff W$\nnsp.\oss
Then\sss the family\sss 
$\mathcal{T}_w\dff,\off w\qff \in\qff W$\sss is\trs Fredholm\dss and\dss its
analytical\dss index\dss is\dss equal\dss to\sss the analytical\dss index of\dss
the family\sss $\mathcal{B}_{\fff w}\dff,\pff w\qff \in\qff W$\nnsp.\oss}

\proof
The relation $\mathcal{C}_{\fff w}$\sss is\dss equal\dss to\sss the difference\vspace{3pt}
\[
\quad
\left(\qff D^{\dff *}\qff\right)^{\dff -\dff 1}
\qff
\mathcal{B}_{\fff w}\qff
D^{\dff -\dff 1}\trf 
\trf -\pff
\left(\qff D^{\dff *}\qff\right)^{\dff -\dff 1}
\qff
P\trf D^{\dff -\dff 1}
\]

\vspace{-12pt}\vspace{3pt}
of\dss a relation with compact\sss resolvent\sss and\sss a bounded operator.\oss
It\sss follows\sss that\sss $\mathcal{C}_{\fff w}$\sss is\dss a relation with
compact\sss resolvent\sss and\dss hence\sss
$U\trf(\trf \mathcal{C}_{\fff w}\trf)\qff \in\qff U\comp$\sss for every $w$\nnsp.\oss
By\sss the equality\qss (\ref{cayley-t-i})\vspace{3pt}
\begin{equation}
\label{tw}
\quad
U\trf(\trf \mathcal{T}_w\trf)
\off =\off
U\trf(\trf \mathcal{C}_{\fff w}\trf)_{\dff H}\trf
U\trf(\trf A\trf)
\pff
\end{equation}

\vspace{-12pt}\vspace{3pt}
for every $w$\nnsp.\oss
Since $A$\sss is\dss a\dss Fredholm\dss operator and\sss
$U\trf(\trf \mathcal{C}_{\fff w}\trf)\qff \in\qff U\comp$\dnsp,\oss
this implies\sss that\sss $\mathcal{T}_w$\sss is\dss Fredholm\dss for every $w$\nnsp.\oss
See\qss \cite{i2},\oss the discussion\sss preceding\trs Lemma\qss 11.2.\oss
Moreover,\oss the equality\qss (\ref{tw})\qss implies\sss that\sss the family\qss
$U\trf(\trf \mathcal{T}_w\trf)\dff,\pff w\qff \in\qff W$\sss
is\dss continuous and\dss hence\sss
$\mathcal{T}_w\dff,\off w\qff \in\qff W$\sss
is\dss not\sss only a\sss family of\trs Fredholm\dss operators,\oss 
but\dss is\dss a\dss Fredholm\trs family.\oss
Since\sss $U\trf(\trf A\trf)$\sss does not\sss depend on $w$\nnsp,\oss
the equality\qss (\ref{tw})\qss together\sss with\dss Lemma\qss 11.2\qss from\qss \cite{i2}\qss
implies\sss that\sss the analytical\dss index of\dss the family\sss
$\mathcal{T}_w\dff,\off w\qff \in\qff W$\sss is\trs Fredholm\dss and\dss its
analytical\dss index\dss is\dss equal\dss to\sss the analytical\dss index of\dss
the family\sss $\mathcal{C}_{\fff w}\dff,\pff w\qff \in\qff W$\nnsp.\oss
It\sss remains\sss to prove\sss that\sss the analytical\dss indices of\dss families\sss
$\mathcal{C}_{\fff w}\dff,\pff w\qff \in\qff W$ and\sss
$\mathcal{B}_{\fff w}\dff,\pff w\qff \in\qff W$ are equal.\oss
The homotopy\vspace{3pt}
\[
\quad
\mathcal{C}_{\fff w\fff,\qff t}
\off =\off
\left(\qff D^{\dff *}\qff\right)^{\dff -\dff 1}
\qff
\left(\qff
\mathcal{B}_{\fff w}\qff
D^{\dff -\dff 1}\trf 
\trf -\pff
t\trf P\trf D^{\dff -\dff 1}
\qff \right)
\dff,\quad
t\qff \in\qff [\trf 0\fff,\qff 1\trf]
\]

\vspace{-12pt}\vspace{3pt}
connects\sss the family\sss
$\mathcal{C}_{\fff w}\dff,\pff w\qff \in\qff W$\sss
with\sss the family\vspace{3pt}\vspace{-0.75pt}
\[
\quad
\left(\qff D^{\dff *}\qff\right)^{\dff -\dff 1}
\qff
\mathcal{B}_{\fff w}\qff
D^{\dff -\dff 1}
\dff,
\quad
w\qff \in\qff W
\pff
\]

\vspace{-12pt}\vspace{3pt}\vspace{-0.75pt}
and\dss hence\sss the index of\dss
$\mathcal{C}_{\fff w}\dff,\pff w\qff \in\qff W$\sss
is\dss equal\dss to\sss the index of\dss the\sss latter\sss family.\oss
But\sss the\sss latter\sss family\dss is\dss conjugate\sss to\sss
$\mathcal{B}_{\fff w}\dff,\pff w\qff \in\qff W$\sss
and\dss hence has\sss the same index.\oss  \eproof

\myuppar{Allowing\dss $T{},\pff A\dff,\pff \gamma_0\dff,\pff \gamma_1$\sss
depending\sss on parameters.}
Again,\oss let\dss
$T_w\dff,\pff w\qff \in\qff W$\sss be a family of\dss 
densely defined closed symmetric operators in $H$\nnsp,\oss
and $A_{\dff w}\dff,\pff w\qff \in\qff W$\sss be a family self-adjoint\sss operators
such\sss that\sss $T_w\qff \subset\qff A_{\dff w}\qff \subset\qff T_w^{\dff *}$\sss
for every\sss $w\qff \in\qff W$\nnsp.\oss
Suppose\sss that\sss the family\sss $A_{\dff w}\dff,\pff w\qff \in\qff W$\sss
is\dss continuous in\sss the\sss topology of\dss the uniform\sss resolvent\sss convergence.\oss
Let\vspace{3pt}
\[
\quad
\gamma_{w\trf 0}\dff,\pff \gamma_{w\trf 1}\dff \colon\dff
H_{\dff 1}\qff \ttoo\qff K
\dff,\quad
w\qff \in\qff W
\]

\vspace{-12pt}\vspace{3pt}
be norm-continuous families of\dss bounded operators.\oss
Suppose\sss that\sss for every\sss $w\qff \in\qff W$\sss
all\sss assumptions of\trs Section\qss \ref{abstract-index}\qss
hold\sss for\sss $T\off =\off T_w$\nsp,\qss $A\off =\off A_{\dff w}$\nsp,\qss 
$\gamma_0\off =\off \gamma_{w\trf 0}$ and\sss
$\gamma_1\off =\off \gamma_{w\trf 1}$\nsp.\oss
In\sss particular{},\oss the operators\sss
$\gamma_{w\trf 0}\dff,\pff \gamma_{w\trf 1}$\sss
extend\dss by continuity\sss to bounded operators\sss 
$\mathcal{D}\dff(\trf T_w^{\dff *}\trf)\qff \ttoo\qff K\fff'$\dnsp,\oss
which we will\sss denote\sss by\sss $\Gamma_{w\trf 0}\dff,\off \Gamma_{w\trf 1}$\sss respectively.\oss
Then\sss $\Gamma_{w\trf 0}$\sss
induces a\sss topological\sss isomorphism\sss
$\kernel\fff T_w^{\dff *}
\qff \ttoo\qff
K\fff'$ 
for every\sss $w\qff \in\qff W$\nnsp.\oss
Let\sss\vspace{3pt}
\[
\quad
\bm{\gamma}_{w}\trf(\dff 0\dff)\dff \colon\dff
K\fff'\qff \ttoo\qff \kernel\fff T_w^{\dff *}
\]

\vspace{-12pt}\vspace{3pt}
be\sss its\sss inverse
and\dss let\vspace{3pt}
\[
\quad
M_{\dff w}\trf(\dff 0\dff)
\off =\off 
\Gamma_{w\trf 1}\dff \circ\trf \bm{\gamma}_w\trf(\dff 0\dff)
\dff \colon\dff
K\fff'\qff \ttoo\qff K\fff'
\pff.
\]

\vspace{-12pt}\vspace{3pt}
We need\dss the family\sss 
$\bm{\gamma}_{w}\trf(\dff 0\dff)\dff,\pff w\qff \in\qff W$\sss
be norm-continuous as\sss a\sss a family of\dss bounded operators\sss
$K\fff'\qff \ttoo\qff H$\sss and\dss the family\sss
$M_{\dff w}\trf(\dff 0\dff)\dff,\pff w\qff \in\qff W$\sss
to be norm-continuous as a family of\dss bounded operators in\sss $K\fff'$\dnsp.\oss
\emph{In\sss the present\sss abstract\sss setting\sss we will\sss simply
assume\sss that\sss this\dss is\dss the case.}\oss
Then\sss the family of\dss reduced\dss boundary operators\sss
$\bm{\Gamma}_{w\dff 1}
\off =\off
\Gamma_{w\dff 1}\qff -\qff M_{\dff w}\trf(\dff 0\dff)\dff \circ\dff \Gamma_{w\dff 0}$\sss
is\dss norm-continuous,\oss as also\sss the families\sss
$\overline{\Gamma}_{w\dff 0}$\sss and\qss 
$\overline{\bm{\Gamma}}_{w\dff 1}$\nsp,\dss $w\qff \in\qff W$\dnsp.\oss

\myuppar{The family\sss $\mathcal{D}\dff(\trf T_w^{\dff *}\trf)\dff,\qff w\qff \in\qff W$\nnsp.}
Recall\dss that\sss the domains\sss $\mathcal{D}\dff(\trf T_w^{\dff *}\trf)$
are equipped\sss with\sss the graph\sss topology.\oss
We may even equip\sss them\sss with\sss the structure of\trs Hilbert\dss spaces
induced\sss from\sss the graphs of\dss operators\sss $T_w^{\dff *}$\nsp.\oss
By\trs Lemma\qss \ref{basic-decomp}\qss\vspace{3pt}
\[
\quad
\mathcal{D}\dff(\trf T_w^{\dff *}\dff)
\off =\off
\mathcal{D}\dff(\trf A_{\dff w}\dff)
\qff \dotplus\qff
\kernel\fff T_w^{\dff *}
\pff.
\]

\vspace{-12pt}\vspace{3pt}
Moreover,\oss replacing $\dotplus$ by\sss the orthogonal\sss
direct\sss sum $\oplus$ does not\sss affect\sss the underlying\sss topology.\oss
By our assumptions,\oss the domains\sss $\mathcal{D}\dff(\trf A_{\dff w}\dff)$
are contained\sss in\sss $H_{\dff 1}$ and\sss their\sss graph\sss topology\dss
is\dss the same as\sss the\sss topology\sss induced\sss from\sss $H_{\dff 1}$\nsp.\oss
Moreover,\oss these domains are equal\dss to\sss the kernels of\dss
bounded operators\sss
$\Gamma_{w\trf 0}\dff \colon\dff
H_{\dff 1}\qff \ttoo\qff K$\sss
continuously depending on\sss $w$\nnsp.\oss
Therefore\sss $\mathcal{D}\dff(\trf A_{\dff w}\dff)\dff,\pff w\qff \in\qff W$\sss
is\dss a continuous\sss family of\trs closed subspaces of\trs $H_{\dff 1}$\nsp.\oss
On\sss the kernel\sss $\kernel\fff T_w^{\dff *}$\sss
the graph\sss topology of\dss $\mathcal{D}\dff(\trf T_w^{\dff *}\dff)$\sss
is\dss equal\sss to\sss the\sss topology\sss induced\sss from\sss $H$\nnsp.\oss
The continuity of\dss the family\sss
$\bm{\gamma}_{w}\trf(\dff 0\dff)\dff,\pff w\qff \in\qff W$\sss
implies\sss that\sss the\sss family\sss
$\kernel\fff T_w^{\dff *}\dff,\pff w\qff \in\qff W$\sss is\dss a continuous family\sss of\dss
subspaces of\dss $H$\nnsp.\oss
It\sss follows\sss that\sss the family\sss of\dss domains\sss
$\mathcal{D}\dff(\trf T_w^{\dff *}\dff)
\off =\off
\mathcal{D}\dff(\trf A_{\dff w}\dff)
\qff \dotplus\qff
\kernel\fff T_w^{\dff *}$
has a natural\sss structure of\dss a\dss Hilbert\dss bundle over $W$
with\sss the structure group $GL\trf(\trf H\trf)$ in\sss the norm\sss topology.\oss
The fibers may\sss be not\sss isometric,\oss but\sss are\sss topologically\sss isomorphic\sss
to\sss the spaces $\mathcal{D}\dff(\trf T_w^{\dff *}\dff)$\sss with\sss the\dss
Hilbert\dss space structures induced\sss from\sss graphs.\oss

\myuppar{Comparing\dss families\dss of\dss boundary\sss triplets.}
Let\sss us apply\sss to\sss the extension $A_{\dff w}$ of\dss $T_w$ and\dss
$\mu\off =\off i$\sss the construction of\dss boundary\sss triplets
from\dss Section\qss \ref{abstract},\oss
and\dss let\vspace{1.5pt}
\begin{equation}
\label{inner-triplet-family}
\quad
\left(\qff \mathcal{K}_{\dff w\dff -}\dff,\off \Gamma_w^{\trf 0}\dff,\off \Gamma_w^{\trf 1} \pff\right)
\end{equation}

\vspace{-12pt}\vspace{1.5pt}
be\sss the resulting\sss boundary\sss triplets.\oss 
\emph{Let\sss us assume\sss that\sss
$\mathcal{K}_{\dff w\dff -}\dff,\pff w\qff \in\qff W$\sss is\dss a continuous family\sss of\dss
subspaces of\dss $H$\nnsp.}\oss
The discussion of\dss comparing\dss boundary\sss triplets\sss in\dss Section\qss \ref{abstract-index}\qss
leads\sss to families of\dss topological\dss isomorphisms\sss\vspace{1.5pt}
\[
\quad
D_w\dff \colon\dff
K^{\dff \partial}\qff \ttoo\qff \mathcal{K}_{\dff w\dff -}\dff,\quad 
w\qff \in\qff W
\qquad
\mbox{and}\quad
\]

\vspace{-34.5pt}
\[
\quad
\mathcal{W}_w\trf \colon\dff
K^{\dff \partial}\dff \oplus\dff K^{\dff \partial}
\qff \ttoo\qff
\mathcal{K}_{\dff w\dff -}\dff \oplus\dff \mathcal{K}_{\dff w\dff -}\dff,
\quad
w\qff \in\qff W
\pff
\]

\vspace{-12pt}\vspace{3pt}
Moreover,\oss there exist\sss bounded self-adjoint\sss operators\sss
$P_w\dff \colon\dff
\mathcal{K}_{\dff w\dff -}\qff \ttoo\qff \mathcal{K}_{\dff \dff -}$\sss
such\sss that\vspace{3pt}
\[
\quad
\overline{\Gamma}_{w\dff 0}
\off =\off
D_w^{\dff -\dff 1}\qff \Gamma_w^{\trf 0}
\quad\dff
\mbox{and}\quad\dff
\]

\vspace{-34.5pt}
\[
\quad
\overline{\bm{\Gamma}}_{w\dff 1}
\off =\off
D_w^{\dff *}\qff \Gamma_w^{\trf 1}
\off +\off
P_w\trf D_w^{\dff -\dff 1}\qff \Gamma_w^{\trf 0}
\off,
\]

\vspace{-12pt}\vspace{4.5pt}
where\sss $\overline{\Gamma}_{w\dff 0}$ and\dss $\overline{\bm{\Gamma}}_{w\dff 1}$
are defined\sss in\sss the same way as\dss 
$\overline{\Gamma}_0$ and\dss
$\overline{\bm{\Gamma}}_1$\nsp.\oss\vspace{-0pt}

\mypar{Lemma.}{compare-cont}
\emph{The families\qss
$D_w\dff,\off \mathcal{W}_w\dff,\off P_w\dff,\off w\qff \in\qff W$\sss
are continuous in\sss the norm\sss topology.\oss}

\proof
By\trs Kuiper's\dss theorem we can\sss trivialize\sss the bundles\sss
$\mathcal{D}\dff(\trf T_w^{\dff *}\dff)\dff,\pff w\qff \in\qff W$\sss
and\sss
$\mathcal{K}_{\dff \dff -}\dff,\pff w\qff \in\qff W$\sss
as\sss the bundles with\sss the structure group\sss
$GL\trf(\trf H\trf)$ in\sss the norm\sss topology.\oss
Therefore\sss the continuity\sss of\dss the families\sss
$D_w\dff,\off \mathcal{W}_w\dff,\off w\qff \in\qff W$\sss
follows from\sss the commutative diagrams at\sss the end of\trs
Section\qss \ref{abstract-index}\qss (with\sss the parameters $w$ added\halfff).\oss
In order\sss to prove\sss the continuity of\dss
$P_w\dff,\off w\qff \in\qff W$\dnsp,\oss
let\sss us\sss write\sss $\mathcal{W}_w$ as a $2\dff \times\dff 2$ matrix.\oss
Then\sss the operator\sss $P_w\trf D_w^{\dff -\dff 1}$\sss is\dss one of\dss the
entries of\dss this matrix.\oss
Therefore\sss the continuity of\dss the family\sss
$P_w\dff,\pff w\qff \in\qff W$\sss
follows from\sss the already\sss proved continuity of\dss the families\sss
$D_w\dff,\off \mathcal{W}_w\dff,\off w\qff \in\qff W$\dnsp.\oss  \eproof

\myuppar{Families of\dss self-adjoint\dss boundary\sss problems.\oss II.}
As above,\oss let\sss $\mathcal{B}_{\fff w}\dff,\pff w\qff \in\qff W$\sss
be a continuous family of\dss self-adjoint\trs Fredholm\dss relations in\sss
$K^{\dff \partial}$\dnsp.\oss
For $w\qff \in\qff W$\dss let\sss $\mathcal{T}_w$\sss be\sss the extension of\dss $T_w$
defined\dss by\sss the equation\dss
$\overline{\bm{\Gamma}}_{w\dff 1}\off =\off \mathcal{B}_{\fff w}\qff \overline{\Gamma}_{w\dff 0}$\nsp.\oss

\mypar{Theorem.}{index-for-changing-a}
\emph{Suppose\sss that\sss either\sss every operator\sss
$A_{\dff w}\dff,\pff w\qff \in\qff W$\sss has compact\sss resolvent,\oss
or every\sss relation $\mathcal{B}_{\fff w}\dff,\pff w\qff \in\qff W$\sss has compact\sss resolvent.\oss
Then\sss the family\sss 
$\mathcal{T}_w\dff,\off w\qff \in\qff W$\sss is\trs Fredholm\dss and\dss its
analytical\dss index\dss is\dss equal\dss to\sss the analytical\dss index of\dss
the family\sss $\mathcal{B}_{\fff w}\dff,\pff w\qff \in\qff W$\nnsp.\oss}

\proof
The proof\trs is\dss similar\sss to\sss the proof\dss of\trs Theorem\qss \ref{index-for-const-a}.\oss
Now\sss the definition of\dss relations $\mathcal{C}_{\fff w}$\sss has a subscript\sss $w$ at\sss
each\dss letter.\oss
Lemma\qss \ref{compare-cont}\qss implies\sss that\sss the family\sss
$\mathcal{C}_{\fff w}\dff,\pff w\qff \in\qff W$\sss is\dss continuous.\oss
The equality\qss (\ref{tw})\qss is\dss replaced\dss by\vspace{3pt}
\[
\quad
U\trf(\trf \mathcal{T}_w\trf)
\off =\off
U\trf(\trf \mathcal{C}_{\fff w}\trf)_{\dff H}\trf
U\trf(\trf A_{\dff w}\trf)
\pff.
\]

\vspace{-12pt}\vspace{3pt}
The discussion\sss preceding\trs Theorem\qss \ref{index-op}\qss shows\sss that\sss
the family\sss $\mathcal{T}_w\dff,\pff w\qff \in\qff W$\sss is\trs Fredholm.\oss
Arguing as in\sss the proof\dss of\trs Theorem\qss \ref{index-op},\oss
we see\sss that\vspace{3pt}
\[
\quad
\ai\trf(\trf \mathcal{T}_{\dff \bullet}\trf)
\off =\off
\ai\trf(\trf A_{\dff \bullet}\trf)
\pff +\pff
\ai\trf(\trf \mathcal{C}_{\fff \bullet}\trf)
\pff.
\]

\vspace{-12pt}\vspace{3pt}
But\sss
$\ai\trf(\trf \mathcal{C}_{\fff \bullet}\trf)
\off =\off
\ai\trf(\trf \mathcal{B}_{\fff \bullet}\trf)$\sss
and\sss
$\ai\trf(\trf A_{\fff \bullet}\trf)
\off =\off
0$\sss
because\sss the operators $A_{\dff w}$ are assumed\sss to be invertible.\oss
The\sss theorem\sss follows.\oss  \eproof

\myuppar{Families of\dss boundary\sss problems
defined\sss in\sss terms of\dss $\gamma_{w\trf 0}\dff,\pff \gamma_{w\trf 1}$\nsp.}
Let\sss $\gamma_w\off =\off \gamma_{w\trf 0}\dff \oplus\dff \gamma_{w\trf 1}$\nsp.\oss
Let\sss 
$\mathcal{B}_{\fff w}
\off \subset\off 
K^{\dff \partial}\dff \oplus\dff K^{\dff \partial}$\dnsp,\oss
where $w$ runs over\sss $W$\dnsp,\oss
be a family of\dss self-adjoint\sss closed\sss relations.\oss
It\dss is\dss only\sss natural\dss to assume\sss that\sss this family\dss is\dss
continuous,\oss say,\oss in\sss the sense of\pss \cite{i2},\oss Section\qss 11.\oss
But\sss what\dss is\dss really\sss needed\dss is\dss the continuity 
of\dss the family of\dss restrictions 
$\mathcal{B}_{\fff w}\trf|\trf K
\off =\off
\mathcal{B}_{\fff w}\qff \cap\qff K\dff \oplus\dff K$\nnsp.\oss
In our\sss the applications both continuity\sss properties will\dss 
be hold\dss by\sss the same reason.\oss

For $w\qff \in\qff W$\sss let\sss 
$\mathcal{T}_{\fff w}$\sss be\sss the restriction of\dss $T^{\dff *}$  to\sss
$\gamma_w^{\dff -\dff 1}\trf(\trf \mathcal{B}_{\fff w}\trf)\qff \subset\qff H_{\dff 1}$\nsp.\oss
\emph{Suppose\sss that\sss $\mathcal{T}_{\fff w}$ is\dss a self-adjoint\sss operator\sss in\sss $H$\sss
for every $w\qff \in\qff W$\nnsp.}\oss
The discussion at\sss the end of\trs Section\qss \ref{abstract-index}\qss
shows\sss that\sss $\mathcal{T}_{\fff w}$\sss is\dss defined\sss
in\sss terms of\dss the boundary\sss triplet\sss
$(\trf K\dff,\off \overline{\Gamma}_{w\dff 0}\dff,\pff \overline{\bm{\Gamma}}_{w\dff 1}\trf)$\sss
by\sss the boundary conditions\vspace{3pt}
\[
\quad
\mathcal{R}_{\fff w}
\off =\off
\Lambda'\dff \oplus\dff \Lambda^{\fff -\dff 1}\qff
\bigl(\qff
\mathcal{B}_{\fff w}\trf|\trf K\qff -\qff M_{\dff w}
\qff\bigr)
\pff,
\]

\vspace{-12pt}\vspace{3pt}
where\sss
$M_{\dff w}
\off =\off 
M_{\dff w}\trf(\dff 0\dff)\trf|\trf K$\nnsp.\oss
By\trs Theorem\qss \ref{index-for-changing-a}\qss
the family\sss 
$\mathcal{T}_w\dff,\off w\qff \in\qff W$\sss is\trs Fredholm\dss and\dss its
analytical\dss index\dss is\dss equal\dss to\sss the analytical\dss index of\dss
$\mathcal{R}_{\fff w}\dff,\pff w\qff \in\qff W$\dss
if\dss either\sss every operator\sss
$A_{\dff w}\dff,\pff w\qff \in\qff W$\sss has compact\sss resolvent,\oss
or every\sss relation $\mathcal{R}_{\dff w}\dff,\pff w\qff \in\qff W$\sss has compact\sss resolvent.\oss

In our applications\sss to differential\dss boundary problems  
the operators\sss $A_{\dff w}$ will\dss be defined\dss by 
self-adjoint\sss elliptic boundary problems of\dss order $1$\nnsp.\oss
Such operators are known\sss to have compact\sss resolvent.\oss
While\sss this\dss is\dss sufficient\sss for our purposes,\oss
we note\sss that\sss the relations\sss $\mathcal{R}_{\dff w}$\sss 
will\sss also have compact\sss resolvent\dss
because\sss the inclusion\sss $K\qff \ttoo\qff K\fff'$\sss 
will\dss be a compact\sss operator{}.\oss

\newpage
\mysection{Differential\pss boundary\qss problems\qss of\qss order\qss one}{differential-boundary-problems}

\myuppar{Sobolev\dss spaces and\dss trace operators.}
Let\sss $X$\sss be a compact\sss manifold\sss with non-empty\sss boundary $Y$\dnsp,\oss
and\dss let\sss $X^{\dff \circ}\off =\off X\qff \smallsetminus\qff Y$\dnsp.\oss
Let\sss $E$\sss be a\dss Hermitian\dss bundle over\sss $X$\sss equipped\sss with\sss
an orthogonal\sss decomposition\sss 
$E\trf|\trf Y\off =\off F\dff \oplus\dff F$\dnsp,\oss
where $F$\sss is\dss a\dss Hermitian\dss bundle over\sss $Y$\dnsp.\oss
Let\sss $H_{\dff 0}$ and\sss $H_{\dff 1}$\sss be\sss the\dss Sobolev\dss spaces\sss
$H_{\dff 0}\dff(\trf X^{\dff \circ}\fff,\qff E\trf)$
and\sss
$H_{\dff 1}\dff(\trf X^{\dff \circ}\fff,\qff E\trf)$\sss
respectively.\oss
Let\sss\vspace{3pt}
\[
\quad 
K^{\dff \partial}
\off =\off
H_{\dff 0}\dff(\trf Y\fff,\qff F\trf)
\quad
\mbox{and}\quad
K
\off =\off
H_{\dff 1/2}\dff(\trf Y\fff,\qff F\trf)
\pff.
\]

\vspace{-12pt}\vspace{3pt}
Then\sss the anti-dual\sss space $K\fff'$\sss is\dss canonically\sss isomorphic\sss to\sss
$H_{\dff -\dff 1/2}\dff(\trf Y\fff,\qff F\trf)$
and\sss
$K\qff \subset\qff K^{\dff \partial}\qff \subset\qff K\fff'$\sss
is\dss a\dss Gelfand\trs triple.\oss
The decomposition\sss 
$E\trf|\trf Y\off =\off F\dff \oplus\dff F$
shows\sss that\sss\vspace{3pt}
\[
\quad 
H_{\dff 1/2}\dff(\trf Y\fff,\qff E\trf|\trf Y\trf)
\off =\off
K\dff \oplus\dff K
\qff,\quad
\]

\vspace{-34pt}
\[
\quad
H_{\dff -\dff 1/2}\dff(\trf Y\fff,\qff E\trf|\trf Y\trf)
\off =\off
K\fff'\dff \oplus\dff K\fff'
\qff,
\quad
\mbox{and}\quad
\]

\vspace{-34pt}
\[
\quad 
H_{\dff 0}\dff(\trf Y\fff,\qff E\trf|\trf Y\trf)
\off =\off
K^{\dff \partial}\dff \oplus\dff K^{\dff \partial}
\qff.\quad
\]

\vspace{-12pt}\vspace{3pt}
Let\sss us denote by\sss $\Pi_{\dff 0}\dff,\pff \Pi_{\dff 1}$\sss
the projections onto\sss the first\sss and\sss the second summands respectively\sss
in each of\dss these decompositions.\oss
Let\sss \vspace{3pt}
\[
\quad 
\gamma\dff \colon\dff
H_{\dff 1}\dff(\trf X^{\dff \circ}\fff,\qff E\trf)
\qff \ttoo\qff
H_{\dff 1/2}\dff(\trf Y\fff,\qff E\trf|\trf Y\trf)
\]

\vspace{-12pt}\vspace{3pt}
be\sss the\sss trace operator and\dss let\dss $\gamma_0\dff,\pff \gamma_1$\sss
be\sss its compositions with\sss the projections\sss
$\Pi_{\dff 0}\dff,\pff \Pi_{\dff 1}$\sss respectively.\oss
Then $\gamma\off =\off \gamma_0\dff \oplus\dff \gamma_1$\nsp.\oss
As\dss is\dss well\dss known,\oss the map $\gamma$\sss is\dss surjective,\oss
admits a continuous section,\oss and\sss its kernel\sss $\kernel\fff \gamma$\sss
is\dss dense in\sss 
$H_{\dff 0}\off =\off H_{\dff 0}\dff(\trf X^{\dff \circ}\fff,\qff E\trf)$\nnsp.\oss

\myuppar{The reference operator.}
Let\sss $P$\sss be
a formally self-adjoint\sss differential\sss operator\sss 
of\dss order $1$ acting on sections of\dss $E$\nnsp.\oss
Then\sss $P$\sss satisfies\sss the\dss Lagrange\dss identity\vspace{3pt}
\[
\quad
\sco{\dff P\dff u\fff,\qff v \dff}
\pff -\pff
\sco{\dff u\dff,\qff P\dff v \dff}
\off =\off
\bsco{\dff i\trf \Sigma\trf \gamma\dff u\dff,\qff 
\gamma\dff v \dff}_{\dff \partial}
\pff,
\]

\vspace{-12pt}\vspace{3pt}
where\sss $u\fff,\qff v\qff \in\qff H_{\dff 1}$\sss and\sss 
$\Sigma$\sss is\dss the coefficient\sss of\dss the normal\sss derivative\sss
$D_{\fff n}\off =\off -\qff i\qff \partial\fff x_{\dff n}$ to $Y$
in\sss $P$\qss (as usual,\pss $n\off =\off \dim\dff X$ and $x_{\dff n}$\sss
is\dss the normal\sss coordinate).\oss
Suppose\sss that\sss $P$\sss is\dss an elliptic operator
and\dss that\sss $i\trf \Sigma$\sss 
has\sss the form\vspace{0pt}
\begin{equation}
\label{i-sigma}
\quad
i\trf \Sigma
\off =\off\dff
\begin{pmatrix}
\off 0 &
1 \qff\off
\vspace{4.5pt} \\
\off\dff -\qff 1 &
0 \qff\off 
\end{pmatrix}
\off
\end{equation}

\vspace{-12pt}\vspace{0pt}
with respect\sss to\sss some decomposition\sss
$E\trf|\trf Y\off =\off F\dff \oplus\dff F$\dnsp.\oss
Without\sss any\sss loss of\dss generality\sss we can assume\sss
that\sss the decomposition\sss from\sss the previous subsection\dss
is\dss equal\dss to\sss this one.\oss
Then\sss the\dss Lagrange\dss identity\sss for\sss $P$\sss
takes\sss the standard\sss form\vspace{3pt}\vspace{-0.375pt}
\begin{equation}
\label{p-lagrange}
\quad
\sco{\dff P\dff u\fff,\qff v \dff}
\pff -\pff
\sco{\dff u\dff,\qff P\dff v \dff}
\off =\off
\bsco{\dff \gamma_1\dff u\dff,\qff 
\gamma_0\dff v \dff}_{\dff \partial}
\pff -\pff
\bsco{\dff \gamma_0\dff u\dff,\qff 
\gamma_1\dff v \dff}_{\dff \partial}
\pff
\end{equation}

\vspace{-12pt}\vspace{3pt}\vspace{-0.375pt}
from\sss the\sss theory of\dss boundary\sss triplets.\oss
Let\sss $T$\sss be\sss the restriction of\dss $P$\sss to\sss
$\kernel\fff \gamma$\nnsp.\oss
The domain $\mathcal{D}\dff(\trf T^{\dff *}\dff)$ of\dss 
the\dss Hilbert\dss space adjoint\sss $T^{\dff *}$\sss is\dss equal\dss to\sss
the space $\mathcal{D}^{\fff 0}_{\trf P}$\sss of\dss all\sss distributions $u$ such\sss that\sss 
$P\fff u\qff \in\qff H_{\dff 0}\dff(\trf X^{\dff \circ}\fff,\qff E\trf)$\nnsp,\oss
where $P\fff u$\sss is\dss understood\sss in\sss the distributional\sss sense,\oss
and\dss the action of\dss $T^{\dff *}$ agrees with\sss the action of\dss $P$ on distributions.\oss
See\qss \cite{g2},\oss Section\qss 4.1.\oss 
The\dss Lagrange\dss identity\sss implies\sss that\sss
$H_{\dff 1}\off =\off H_{\dff 1}\dff(\trf X^{\dff \circ}\fff,\qff E\trf)$\sss
is\dss contained\sss in\sss $\mathcal{D}\dff(\trf T^{\dff *}\dff)$\nnsp.\oss
Moreover,\pss $H_{\dff 1}$\sss is\dss dense in\sss $\mathcal{D}\dff(\trf T^{\dff *}\dff)$
equipped\sss with\sss the graph\sss topology.\oss
This\dss is\dss a classical\sss result\sss essentially\sss due\sss to\dss Lions\dss and\trs Magenes\qss \cite{lm1}\qss
(see\qss \cite{lm1},\oss footnote\qss (\nnsp$\stackrel{6}{}$\nnsp)\qss on\sss p.\qss 147\fff).\oss
The corresponding result\sss for\sss $P\off =\off 1\qff -\qff \Delta$\nnsp,\oss
where $\Delta$\sss is\dss the\dss Laplace\dss operator{},\oss
is\dss proved\sss in\qss \cite{g2},\oss Theorem\qss 9.8.\oss
Mutatis mutandis\sss this proof\dss applies in\sss the present\sss context.\oss

Let\sss $A$ be\sss the restriction of\dss $P$\sss to\sss
$\kernel\fff \gamma_0$\nsp.\oss
Suppose\sss that\sss each of\dss the boundary conditions\sss $\gamma_0\off =\off 0$\sss
and\sss $\gamma_1\off =\off 0$\sss
satisfies\sss the\dss Shapiro--Lopatinskii\dss condition\sss for $P$\dnsp.\oss
Then\sss the operators $A$\sss and\sss $A'\off =\off P\trf|\trf \kernel\fff \gamma_1$\sss 
are unbounded\trs Fredholm\dss operators in\sss 
$H_{\dff 0}\off =\off H_{\dff 0}\dff(\trf X^{\dff \circ}\fff,\qff E\trf)$\sss
and are self-adjoint\sss operators with\sss the domains\sss
$\mathcal{D}\dff(\trf A\trf)\fff,\off
\mathcal{D}\dff(\trf A'\trf)
\qff \subset\qff 
H_{\dff 1}$\nsp.\oss
See\qss \cite{i2},\oss Sections\qss 5\qss and\qss 7.\oss
In\sss particular,\pss $A$\sss is\dss contained\sss in\sss $T^{\dff *}$ 
and\dss the kernel\sss of\dss $A$\sss is\dss finitely dimensional.\oss
In order\sss to use $A$ as\sss the reference operator we have\sss to\qss \emph{assume\sss
that\dss $\kernel\fff A\off =\off 0$\nnsp.}\oss
Then\sss $A$\sss has a bounded everywhere defined\sss inverse\sss
$A^{\fff -\dff 1}\dff \colon\dff H\qff \ttoo\qff H$\nnsp.\oss

This completes\sss the verification of\dss assumptions of\trs Section\qss \ref{abstract-index},\oss
except\sss of\dss the assumptions concerned\sss with\sss $M\trf(\dff 0\dff)$\nnsp.\oss
Therefore\sss the results of\trs Section\qss \ref{abstract-index},\oss
in\sss particular\dss Theorems\qss \ref{extending-gamma}\qss and\qss  \ref{extended-isomorphism},\oss
and\sss the construction of\dss the reduced\dss boundary\sss triplet\sss
apply\sss in\sss the present\sss situation.\oss

\myuppar{Calder\'{o}n's\dss method.}
We would\dss like\sss to give a more explicit\sss description of\dss the reduced\dss
boundary\sss triplet\sss in\sss the above context,\oss and,\oss in\sss particular{},\oss
to determine\sss the operator $M\trf(\dff 0\dff)$\nnsp.\oss
To\sss this end\sss we will\sss apply\sss the\dss Calder\'{o}n's\sss method
as presented\dss by\dss G.\dss Grubb\dss in\qss \cite{g2},\oss Chapter\qss 11.\oss
In order\sss to be able\sss to freely\sss refer\sss to\dss G.\dss Grubb\qss \cite{g2}\qss
\emph{we will\sss assume\sss that\dss $P$ extends\sss
to an\sss invertible differential\sss operator over\sss the double\sss $\widehat{X}$\sss of\dss $X$\sss
acting on\sss the sections of\dss the double\sss $\widehat{E}$\sss of\dss the bundle $E$\nnsp.}\oss
We will\sss denote\sss the extended operator still\dss by $P$\dnsp,\oss
and\sss its\sss inverse by\sss $Q$\nnsp.\oss
We will\sss use\sss the notation\sss $\gamma$\sss also for\sss the\sss trace operator\sss
taking\sss the sections over\sss $\widehat{X}$\sss to\sss their restrictions\sss to $Y$\dnsp,\oss
and denote by\sss $r_{\dff +}$\sss the\sss operator\sss
taking\sss the sections over\sss $\widehat{X}$\sss to\sss their restrictions\sss to $X^{\dff \circ}$\dnsp.\oss
Note\sss that\sss the difference between $X$ and $X^{\dff \circ}$\sss
is\dss essential\fff:\oss some sections are generalized ones
with singularities along $Y$\nnsp;\oss the operator $r_{\dff +}$ erases such singularities.\oss
Let\vspace{3pt}\vspace{-0.375pt}
\[
\quad
\mathfrak{A}\off =\off i\trf \Sigma
\quad
\mbox{and}\quad
K_{\dff +}
\off =\off
-\qff r_{\dff +}\dff \circ\trf Q\dff \circ\dff \gamma^{\dff *} \circ\trf \mathfrak{A}
\pff,
\]

\vspace{-12pt}\vspace{3pt}\vspace{-0.375pt}
where\sss
$\gamma^{\dff *}\dff \colon\dff
H_{\dff -\dff 1/2}\dff(\trf Y\fff,\qff E\trf|\trf Y\trf)
\qff \ttoo\qff
H_{\dff -\dff 1}\dff(\qff \widehat{X}\fff,\qff \widehat{E}\trf)$\sss
is\dss the adjoint\qss (dual\fff)\qss operator of\sss $\gamma$\nnsp.\oss
Following\dss G.\dss Grubb\qss \cite{g2},\oss let\sss us set\sss
$Z_{\trf 0}\off =\off \kernel\fff T^{\dff *}$\dnsp.\oss
Let\sss $\Gamma\off =\off \Gamma_0\dff \oplus\dff \Gamma_1$\nsp.\oss
Since\sss $Z_{\trf 0}\qff \subset\qff \mathcal{D}\dff(\trf T^{\dff *}\dff)$\nnsp,\oss
the extended\sss trace operator\sss
$\Gamma$\sss is\dss well\sss defined on\sss $Z_{\dff 0}$ and\sss maps\sss
$Z_{\trf 0}$\sss into\sss
$H_{\dff -\dff 1/2}\dff(\trf Y\fff,\qff E\trf|\trf Y\trf)$\nnsp.\oss
Let\sss $N_{\dff 0}\off =\off \Gamma\trf(\trf Z_{\trf 0}\trf)$\nnsp.\oss
Note\sss that\sss the support\sss of\dss every section in\sss the image of\sss $\gamma^{\dff *}$\sss
is\dss contained\sss in $Y$\dss (they are $\delta$\dnsp-functions
in\sss the direction\sss transverse\sss to $Y$\nnsp).\oss
This implies\sss that\sss the image of\dss\vspace{3pt}
\[
\quad 
K_{\dff +}\dff \colon\dff
H_{\dff -\dff 1/2}\dff(\trf Y\fff,\qff E\trf|\trf Y\trf)
\qff \ttoo\qff
H_{\dff 0}\dff(\trf X^{\dff \circ}\fff,\qff E\trf)
\]

\vspace{-12pt}\vspace{3pt}
is\dss contained\sss in\sss the kernel\sss of\sss $P$\sss over $X^{\dff \circ}$\dnsp,\oss
or,\oss equivalently,\oss in\sss $Z_{\trf 0}$\nsp.\oss
Moreover,\pss $K_{\dff +}$ induces an\sss isomorphism\sss
$N_{\dff 0}\qff \ttoo\qff Z_{\trf 0}$\nsp,\oss
and\dss its\sss inverse is\dss induced\sss by $\Gamma$\nnsp.\oss
See\qss \cite{g2},\oss Proposition\qss 11.5.\oss

The composition\sss
$C_{\dff +}\off =\off \Gamma\trf \circ\trf K_{\dff +}$\sss 
is\dss known as\sss the\qss \emph{Calder\'{o}n\dss projector}.\oss
It\dss is\dss indeed a projection,\oss
i.e.\sss $C_{\dff +}\fff \circ\trf C_{\dff +}\off =\off C_{\dff +}$\nsp,\oss
and\dss is\dss a pseudo-differential\sss operator of\dss order $0$\nnsp.\oss
See\qss \cite{g2},\oss Proposition\qss 11.7.\oss
The symbol\sss $c_{\dff +}$ of\trs $C_{\dff +}$\sss can be expressed\sss in\sss terms of\dss
the\pss \emph{plus-integral}\pss of\dss 
the symbol\sss $q$\sss of\sss $Q$\nnsp.\oss
We refer\sss to\dss H\"{o}rmander\qss \cite{h},\oss Lemma\qss 18.2.16\qss for\sss the definition of\dss
the plus-integral\dss
$\int^{\dff +}\nsp f\trf(\trf t\trf)\qff dt$\nnsp.\oss
If\dss we write $q$\sss locally\sss in\sss terms of\dss the coordinates\sss
$(\trf y\fff,\qff x_{\dff n}\fff,\qff u\fff,\qff t\trf)$\nnsp,\oss
where $y$\sss is\dss the coordinates on $Y$\dnsp,\qss $x_{\dff n}$\sss is\dss
the normal\sss coordinate,\oss and\sss $u\fff,\qff t$ are dual\dss to\sss $y\fff,\qff x_{\dff n}$\nsp,\oss
then\vspace{1.5pt}
\[
\quad
c_{\dff +}\dff(\trf y\fff,\qff u\trf)
\off =\off
-\qff
(\trf 2\dff \pi\trf)^{\dff -\dff 1}\qff
\left(\qff\int^{\dff +} q\trf(\trf y\fff,\qff 0\fff,\qff u\fff,\qff t\trf)\qff dt
\qff\right)
\qff \circ\qff
\mathfrak{A}
\pff.
\]

\vspace{-12pt}\vspace{1.5pt}
This follows\sss from\dss H\"{o}rmander\qss \cite{h},\oss  Theorem\qss 18.2.17.\oss 
This\sss theorem determines\sss the symbol\sss of\dss operators of\dss the form\sss 
$\Gamma\dff \circ\dff r_{\dff +}\dff \circ\trf R\trf \circ\dff \gamma^{\dff *}$\dnsp,\oss
where $R$\sss is\dss a polyhomogeneous pseudo-differential\sss operator
satisfying\sss the\sss transmission condition.\oss
Since $Q$\sss is\dss the inverse of\dss a differential\sss operator,\oss
it\sss applies\sss to\sss $R\off =\off Q$\nnsp.\oss
The boundary operator $\Gamma$\dnsp,\oss being\sss the extension of\dss $\gamma$\sss
by\sss the continuity,\oss agrees with\dss H\"{o}rmander's\dss one.\oss

\myuppar{Computing\sss the plus-integral.}
At\sss the boundary $Y$\sss the symbol\sss $p$ of\dss $P$\sss has\sss the form\vspace{3pt}
\begin{equation}
\label{symbol}
\quad
p\trf(\trf y\fff,\qff 0\fff,\qff u\fff,\qff t\trf)
\off =\off
t\trf \sigma_y\qff +\qff \tau_{\fff u}\trf(\trf y\trf)
\pff,
\end{equation}

\vspace{-12pt}\vspace{3pt}
where $\sigma_y$\sss is\dss the coefficient\sss of\dss $D_{\fff n}$ at\sss $y\qff \in\qff Y$\dss
(so,\oss in\sss fact,\pss $\sigma_y\off =\off \Sigma$\nsp)\qss
and\sss $\tau_{\fff u}\trf(\trf y\trf)$\sss is\dss a self-adjoint\sss operator\sss
$E_{\dff y}\qff \ttoo\qff E_{\dff y}$\sss 
linearly depending on $u$ for every\sss $y\qff \in\qff Y$\dnsp.\oss
We will\sss write simply\sss $\tau_{\fff u}$ for\sss $\tau_{\fff u}\trf(\trf y\trf)$\nnsp.\oss
As in\qss \cite{i2},\oss let\dss 
$\rho_{\dff u}\off =\off \sigma_y^{\dff -\dff 1}\qff \tau_{\fff u}$\nsp.\oss
Then\vspace{3pt}
\[
\quad
q\trf(\trf y\fff,\qff 0\fff,\qff t\fff,\qff u\trf)
\off =\off
p\trf(\trf y\fff,\qff 0\fff,\qff u\fff,\qff t\trf)^{\dff -\dff 1}
\off =\off
(\trf t\trf \sigma_{\fff y}\qff +\qff \tau_{\dff u}\trf)^{\dff -\dff 1}
\quad
\mbox{and}
\]

\vspace{-33pt}
\[
\quad
c_{\dff +}\dff(\trf y\fff,\qff u\trf)
\off =\off
-\qff
(\trf 2\dff \pi\trf)^{\dff -\dff 1}\qff
\int^{\dff +} \frac{i\trf \sigma_{\fff y}}{t\trf \sigma_{\fff y}\qff +\qff \tau_{\dff u}}
\qff dt
\off =\off
-\qff
(\trf 2\dff \pi\trf)^{\dff -\dff 1}\qff i\qff
\int^{\dff +} \frac{1}{t\qff +\qff \rho_{\dff u}}
\qff dt
\pff.
\]

\vspace{-12pt}\vspace{3pt}
By\qss \cite{h},\oss Remark after\dss Lemma\qss 18.2.16,\oss
the\sss last\sss plus-integral\dss is\dss equal\dss to
$2\dff \pi\dff i$\sss times\sss the sum of\dss residues of\dss
$(\trf z\qff +\qff \rho_{\dff u}\trf)^{\dff -\dff 1}$\sss
in\sss the upper half-plane.\oss
The poles of\dss
$(\trf z\qff +\qff \rho_{\dff u}\trf)^{\dff -\dff 1}$\sss
in\sss the upper half-plane correspond\sss to\sss the eigenvalues of\dss $\rho_{\dff u}$\sss
in\sss the\sss lower half-plane.\oss
It\sss follows\sss that\sss the operator\sss
$c_{\dff +}\dff(\trf y\fff,\qff u\trf)$\sss
is\dss equal\dss to\qss
({\fff}because\sss 
$-\qff (\trf 2\dff \pi\trf)^{\dff -\dff 1}\dff i\dff \cdot\dff
2\dff \pi\dff i
\off =\off
1$\nnsp)\qss 
the projection onto\sss the subspace\sss
$\mathcal{L}_{\dff -}\dff(\trf \rho_{\dff u}\trf)$
having\sss $\mathcal{L}_{\dff +}\dff(\trf \rho_{\dff u}\trf)$
as\sss its kernel,\oss
where,\oss as in\qss \cite{i2},\oss we denote by\sss
$\mathcal{L}_{\dff -}\dff(\trf \rho_{\dff u}\trf)$ and\sss
$\mathcal{L}_{\dff +}\dff(\trf \rho_{\dff u}\trf)$
the sums of\dss the generalized eigenspaces of\dss $\rho_{\dff u}$\sss
corresponding\sss to eigenvalues $\lambda$\sss with\sss $\image\fff \lambda\qff <\qff 0$
and\sss $\image\fff \lambda\qff <\qff 0$\sss respectively.\oss

\myuppar{The matrix of\dss the\dss Calder\'{o}n\dss projector.} 
Let\sss us write $C_{\dff +}$ as\sss the matrix\vspace{1.5pt}
\[
\quad
C_{\dff +}
\off =\off\dff
\begin{pmatrix}
\off C_{\dff +\dff 0\dff 0} &
C_{\dff +\dff 0\dff 1} \qff\off
\vspace{4.5pt} \\
\off\dff C_{\dff +\dff 1\dff 0} &
C_{\dff +\dff 1\dff 1} \qff\off 
\end{pmatrix}
\off
\]

\vspace{-12pt}\vspace{1.5pt}
with respect\sss to\sss the decomposition $E\off =\off F\dff \oplus\dff F$\dnsp.\oss

\mypar{Lemma.}{elliptic-blocks}
\emph{Each of\dss the blocks\sss $C_{\dff +\dff i\fff j}$\sss is\dss an elliptic operator.}\oss

\proof
Cf.\qss Grubb\qss \cite{g2},\oss Lemma\qss 11.16.\oss
Let\sss $I_{\dff j}$\nsp, $j\off =\off 0\fff,\qff 1$\sss
be\sss the inclusions of\dss sections of\dss $F$\sss into\sss the sections\sss of\dss
$E\off =\off F\dff \oplus\dff F$ as\sss the sections of\dss the first\sss and\dss
the second summands respectively.\oss
Then
$C_{\dff +\dff i\fff j}\off =\off \Pi_{\dff i}\dff \circ\qff C_{\dff +}\dff \circ\qff I_{\dff j}$\nsp.\oss
The symbol\sss $s_{\dff +\dff i\fff j}\dff(\trf y\fff,\qff u\trf)$\sss 
of\dss $C_{\dff +\dff i\fff j}$\sss is\dss related\sss in\sss the same way\sss to\sss
the symbol $c_{\dff +}\dff(\trf y\fff,\qff u\trf)$ of\dss $C_{\dff +}$\nsp.\oss
We need\sss to check\sss that\sss $s_{\dff +\dff i\fff j}\dff(\trf y\fff,\qff u\trf)$\sss
is\dss an isomorphism\sss if\dss
$u\off \neq\off 0$\nnsp.\oss
Since\sss the boundary condition $\gamma_1\off =\off 0$\sss 
satisfies\sss the\dss Shapiro--Lopatinskii\dss condition,\oss
the subspace $\mathcal{L}_{\dff -}\dff(\trf \rho_{\dff u}\trf)$\sss is\dss transverse\sss to
$K\fff'\dff \oplus\dff 0$ for $u\off \neq\off 0$\nnsp.\oss
Since $P$\sss is\dss a differential\sss operator{},\pss
$\mathcal{L}_{\dff +}\dff(\trf \rho_{\dff u}\trf)
\off =\off 
\mathcal{L}_{\dff -}\dff(\trf -\qff \rho_{\dff u}\trf)
\off =\off 
\mathcal{L}_{\dff -}\dff(\trf \rho_{\dff -\dff u}\trf)$\sss
is\dss also\sss transverse\sss to $K\fff'\dff \oplus\dff 0$\nnsp.\oss
Similarly,\oss since $\gamma_0\off =\off 0$\sss 
satisfies\sss the\dss Shapiro--Lopatinskii\dss condition,\oss
the subspaces $\mathcal{L}_{\dff -}\dff(\trf \rho_{\dff u}\trf)$\sss
and\sss $\mathcal{L}_{\dff +}\dff(\trf \rho_{\dff u}\trf)$\sss
are\sss transverse\sss to\sss $0\dff \oplus\dff K\fff'$\dnsp.\oss
Together\sss with\sss the description of\dss
$c_{\dff +}\dff(\trf y\fff,\qff u\trf)$ given above\sss
these\sss transversality\dss properties\sss imply\sss
that\sss $s_{\dff +\dff i\fff j}\dff(\trf y\fff,\qff u\trf)$\sss
is\dss an isomorphism\sss if\dss
$u\off \neq\off 0$\nnsp.\oss
The\sss lemma\sss follows.\oss  \eproof

\mypar{Lemma.}{graph}
\emph{$N_{\dff 0}$\sss is\dss equal\dss to\sss the graph of\dss an operator\sss
$K\fff'\qff \ttoo\qff K\fff'$\sss induced\dss by a pseudo-differential\dss
operator\dss $\varphi$ of\dss order $0$\nnsp.\oss
Moreover{},\pss $\varphi$\sss is\dss elliptic.\oss}

\proof
First,\oss let\sss us prove\sss that\sss
$N_{\trf 0}\qff \cap\qff 0\dff \oplus\dff K\fff'\off =\off 0$\nnsp.\oss
If\dss $u\qff \in\qff N_{\trf 0}\qff \cap\qff 0\dff \oplus\dff K\fff'$\nnsp,\oss
then\sss $u\off =\off \Gamma\dff x$\sss for some\sss $x\qff \in\qff \kernel\fff T^{\dff *}$
and\sss $x\qff \in\qff \kernel\fff \Gamma_0$\nsp.\oss
Lemma\qss \ref{gamma-0-ker}\qss implies\sss that\sss
$\kernel\fff \Gamma_0\off =\off \kernel\fff \gamma_0$\nsp.\oss
It\sss follows\sss that\sss
$x
\qff \in\qff 
\kernel\fff \gamma_0
\off =\off 
\mathcal{D}\dff(\trf A\trf)$\nnsp.\oss
In\sss turn,\oss this implies\sss $x\qff \in\qff \kernel\fff A$\sss
and\dss hence $x\off =\off 0$\nnsp.\oss
This proves\sss that\sss
$N_{\trf 0}\qff \cap\qff 0\dff \oplus\dff K\fff'\off =\off 0$\nnsp.\oss
Next,\oss Theorem\qss \ref{extended-isomorphism}\qss implies\sss that\sss
the projection of\dss $N_{\trf 0}$\sss to\sss $K\fff'\dff \oplus\dff 0$\sss
is\dss surjective.\oss
It\sss follows\sss that\sss $N_{\trf 0}$\sss
is\dss equal\dss to\sss the graph of\dss an operator\sss
$\varphi\dff \colon\dff K\fff'\qff \ttoo\qff K\fff'$\dnsp.\oss

The fact\sss that\sss $C_{\dff +}$ is\dss a projection onto\sss
the graph of\sss $\varphi$\sss implies\sss that\sss
$C_{\dff +\dff 1\fff 0}\off =\off \varphi\dff \circ\trf C_{\dff +\dff 0\fff 0}$\nsp.\oss
Lemma\qss \ref{elliptic-blocks}\qss implies\sss that\sss there exists a parametrix\sss
$S$ for $C_{\dff +\dff 0\fff 0}$\nsp,\oss i.e.\oss a pseudo-differential\sss operator
of\dss order $0$\sss such\sss that\sss
$C_{\dff +\dff 0\fff 0}\dff \circ\dff S
\off =\off
1\qff -\qff R$\nnsp,\oss
where $R$\sss is\dss a smoothing operator of\dss finite rank.\oss
It\sss follows\sss that\sss
$C_{\dff +\dff 1\fff 0}\dff \circ\dff S
\off =\off
\varphi\qff -\qff C_{\dff +\dff 0\fff 0}\dff \circ\dff R$\sss
and\dss hence\sss
$\varphi
\off =\off
C_{\dff +\dff 1\fff 0}\dff \circ\dff S
\qff +\qff
C_{\dff +\dff 0\fff 0}\dff \circ\dff R$\nnsp.\oss
The operator\sss $C_{\dff +\dff 0\fff 0}\dff \circ\dff R$\sss
is\dss a smoothing operator of\dss finite rank\sss together with $R$\nnsp.\oss
Since\sss $C_{\dff +\dff 1\fff 0}\dff \circ\dff S$\sss is\dss a pseudo-differential\sss
operator of\dss order $0$\nnsp,\oss this implies\sss that\sss $\varphi$\sss is\dss
also such an operator{}.\oss
Since\sss $C_{\dff +\dff 1\fff 0}$ and $S$ are elliptic,\oss 
this implies\sss that\sss $\varphi$\sss is\dss also elliptic.\oss  \eproof%\vspace{-6pt}

\myuppar{The operators\sss
$\bm{\gamma}\trf(\dff 0\dff)\fff,\pff
M\trf(\dff 0\dff)$\sss and\sss $\varphi$\nnsp.}
Let\sss $\varphi$ be\sss the operator\dss from\trs
Lemma\qss \ref{graph},\oss
and\dss let\vspace{1.5pt}
\[
\quad 
\overline{\varphi}
\off =\off
\id\dff \oplus\qff \varphi\dff \colon\dff 
K\fff'
\qff \ttoo\qff 
K\fff'\dff \oplus\dff K\fff'
\]

\vspace{-12pt}\vspace{1.5pt}
be\sss the map\sss 
$x\off \longmapsto\off (\trf x\fff,\qff \varphi\trf(\trf x\trf)\trf)$\nnsp.\oss

Since\sss $K_{\dff +}\trf|\trf N_{\trf 0}$\sss 
is\dss the inverse of\dss 
$\Gamma$\nnsp,\oss\vspace{3pt}
\[
\quad 
\Gamma_0\dff \circ\dff K_{\dff +}\dff \circ\qff \overline{\varphi}
\off =\off
\Pi_{\dff 0}\dff \circ\dff \Gamma\dff \circ\dff K_{\dff +}\dff \circ\qff \overline{\varphi}
\off =\off
\Pi_{\dff 0}\dff \circ\qff \overline{\varphi}
\off =\off
\id
\quad
\mbox{and}
\]

\vspace{-33pt}
\[
\quad 
\Gamma_1\dff \circ\dff K_{\dff +}\dff \circ\qff \overline{\varphi}
\off =\off
\Pi_{\dff 1}\dff \circ\dff \Gamma\dff \circ\dff K_{\dff +}\dff \circ\qff \overline{\varphi}
\off =\off
\Pi_{\dff 1}\dff \circ\qff \overline{\varphi}
\off =\off
\varphi
\pff.
\]

\vspace{-12pt}\vspace{1.5pt}
Since\sss $\bm{\gamma}\trf(\dff 0\dff)$\sss 
is\dss the inverse of\dss the restriction\sss
$\Gamma_0 \trf|\trf \kernel\fff T^{\dff *}$\dnsp,\oss
it\sss follows\sss that\oss\vspace{3pt}
\[
\quad 
\bm{\gamma}\trf(\dff 0\dff)
\off =\off
\left(\qff \Gamma_0\trf|\trf \kernel\fff T^{\dff *}\qff\right)^{\dff -\dff 1}
\off =\off
K_{\dff +}\dff \circ\qff \overline{\varphi}
\]

\vspace{-12pt}\vspace{3pt}
and\dss hence\sss
$M\trf(\dff 0\dff)
\off =\off
\Gamma_1\dff \circ\dff \bm{\gamma}\trf(\dff 0\dff)
\off =\off
\Gamma_1\dff \circ\dff K_{\dff +}\dff \circ\qff \overline{\varphi}\dff
\off =\off
\varphi$\nnsp.\oss

\myuppar{The operators\sss $M\trf(\dff 0\dff)$ and\sss $M$\nnsp.}
Since\sss 
$M\trf(\dff 0\dff)\off =\off \varphi$\nnsp,\oss
Lemma\qss \ref{graph}\qss implies\sss that\sss
$M\trf(\dff 0\dff)$\sss is\dss induced\dss by a pseudo-differential\sss operator of\dss order $0$\nnsp,\oss
which we will,\oss by an abuse of\dss notations,\oss still\sss denote by\sss $M\trf(\dff 0\dff)$\nnsp.\oss
This implies,\oss in\sss particular{},\oss
that\sss $M\trf(\dff 0\dff)$\sss leaves $K$\sss invariant.\oss
The symbol\sss of\dss $M\trf(\dff 0\dff)$\sss
is\dss related\sss to\sss the symbol\sss of\dss $C_{\dff +}$\sss
in\sss the same way as $M\trf(\dff 0\dff)$\sss is\dss related\sss to $C_{\dff +}$\sss
and\dss hence\dss is\dss the bundle map\sss having as\sss its\sss graph\dss the graph of\dss
$c_{\dff +}$\nsp.\oss 
It\sss follows\sss that\sss over a point\sss $u\qff \in\qff S\fff Y$\dnsp,\oss
where\sss $S\fff Y$\sss is\dss the unite sphere bundle of\dss $Y$\dnsp,\oss
the symbol\sss of\dss $M\trf(\dff 0\dff)$\sss is\dss equal\dss to\sss
the map\sss having\sss
$\mathcal{L}_{\dff -}\dff(\trf \rho_{\dff u}\trf)$
as its graph.\oss
It\sss follows\sss that\sss this symbol\dss is\dss an\sss isomorphism,\oss
i.e.\qss the operator\sss $M\trf(\dff 0\dff)$\sss is\dss elliptic.\oss
In\sss particular{},\oss $M\trf(\dff 0\dff)$\sss is\trs Fredholm\dss
as an operator\sss $K\fff'\qff \ttoo\qff K\fff'$\dnsp.\oss

The restriction\sss $M\off =\off M\trf(\dff 0\dff)\trf|\trf K$\sss 
is\dss essentially\sss the same pseudo-differential\sss operator{},\oss
or{},\oss more precisely,\oss is\dss induced\sss by\sss it.\oss
It\sss follows\sss that\sss the operator $K\qff \ttoo\qff K$\sss
induced\dss by\sss $M\trf(\dff 0\dff)$\sss is\trs Fredholm.\oss
This verifies\sss the assumptions of\trs Section\qss \ref{abstract-index}\qss
concerned\sss with\sss $M\trf(\dff 0\dff)$\nnsp.\oss
As we saw in\dss Section\qss \ref{abstract-index},\oss
this implies\sss that\sss $M$\sss is\dss a self-adjoint\sss
as an\sss unbounded\sss operator\sss from\sss $K\fff'$\sss to $K$\nnsp.\oss

\mypar{Lemma.}{sur}
\emph{The boundary\sss map\dss $\gamma_0$\sss
maps\sss $(\trf \kernel\fff T^{\dff *}\trf)\dff \cap\dff H_{\dff 1}$ onto $K$\nnsp.\oss}

\proof
The operator\sss $K_{\dff +}$\sss is\dss known\sss to continuously\sss map\sss
$H_{\dff 1/2}\dff(\trf Y\fff,\qff E\trf|\trf Y\trf)$\sss
into\sss
$H_{\dff 1}\dff(\trf X^{\dff \circ}\fff,\qff E\trf)$\nnsp.\oss
See\dss G.\dss Grubb\qss \cite{g2},\oss (11.17).\oss
Since\dss $\overline{\varphi}$\dss together\sss with\sss $\varphi$\sss
is\dss a pseudo-differential\sss 
operator of\dss order $0$\nnsp,\oss
the operator\sss
$\bm{\gamma}\trf(\dff 0\dff)
\off =\off
K_{\dff +}\fff \circ\qff \overline{\varphi}$\dss
maps\sss 
$K
\off =\off
H_{\dff 1/2}\dff(\trf Y\fff,\qff F\trf)$\sss
into\sss
$H_{\dff 1}
\off =\off
H_{\dff 1}\dff(\trf X^{\dff \circ}\fff,\qff E\trf)$\sss
and\dss hence\sss into\sss
$(\trf \kernel\fff T^{\dff *}\trf)\dff \cap\dff H_{\dff 1}$\nsp.\oss
Since\sss $\bm{\gamma}\trf(\dff 0\dff)$\sss is\dss the inverse of\dss
$\Gamma_0\trf|\trf \kernel\fff T^{\dff *}$\sss
and\sss 
$\gamma_0\off =\off \Gamma_0\trf|\trf H_{\dff 1}$\nsp,\oss
it\sss follows\sss that\sss $\gamma_0$\sss maps
$(\trf \kernel\fff T^{\dff *}\trf)\dff \cap\dff H_{\dff 1}$ onto $K$\nnsp.\oss  \eproof

\myuppar{The operators\sss $M\trf(\dff 0\dff)$ and\sss $M$ 
and\sss the spaces of\trs Cauchy\dss data.}
Since\sss 
$M\trf(\dff 0\dff)\off =\off \varphi$\nnsp,\oss
Lemma\qss \ref{graph}\qss implies\sss that\sss the operator\sss $M\trf(\dff 0\dff)$\sss
can\sss be characterized as\sss the operator\sss
$K\fff'\qff \ttoo\qff K\fff'$\sss having as its graph\sss 
the image\sss $N_{\dff 0}$ of\dss 
$Z_{\trf 0}\off =\off \kernel\fff T^{\dff *}$\sss 
under\sss the map\sss
$\Gamma\off =\off \Gamma_0\dff \oplus\dff \Gamma_1$\nsp,\oss
i.e.\dss the space of\trs Cauchy\dss data of\dss the equation\sss $T^{\dff *} u\off =\off 0$\nnsp.\oss

The operator $M$ admits a similar characterization.\oss
Namely,\oss 
let\sss 
$Z_{\dff 1}\off =\off H_{\dff 1}\dff \cap\dff Z_{\trf 0}$
be\sss the space of\dss solutions of\trs $T^{\dff *} u\off =\off 0$
belonging\sss to $H_{\dff 1}$\nsp.\oss
Let\sss
$N_{\dff 1}\off =\off \gamma\trf(\trf Z_{\dff 1}\trf)
\qff \subset\qff 
K\dff \oplus\dff K$\sss
be\sss the space of\dss the corresponding\dss Cauchy\dss data.\oss
Lemma\qss \ref{graph}\qss implies\sss that\sss 
$N_{\dff 1}\dff \cap\dff (\trf 0\dff \oplus\dff K\trf)\off =\off 0$\nnsp,\oss
and\trs Lemma\qss \ref{sur}\qss implies\sss that\sss
the projection of\dss $N_{\dff 1}$\sss to\sss $K\dff \oplus\dff 0$\sss
is\dss surjective.\oss
Therefore\sss $N_{\dff 1}$\sss is\dss equal\dss to\sss the graph of\dss an operator\sss
$K\qff \ttoo\qff K$\nnsp.\oss
Clearly,\oss this\dss is\dss
the operator\sss induced\dss by $M\trf(\dff 0\dff)$\nnsp,\oss
i.e.\dss $M$\nnsp.\oss
It\sss follows\sss that\sss $M$\sss considered as an\sss operator\sss from\sss $K\fff'$\sss to $K$\sss
can\sss be characterized as\sss the operator\sss  having as its graph\sss 
the image\sss $N_{\dff 1}$ of\dss 
$Z_{\dff 1}\off =\off H_{\dff 1}\dff \cap\dff \kernel\fff T^{\dff *}$\sss 
under\sss the map\sss
$\gamma\off =\off \gamma_0\dff \oplus\dff \gamma_1$\nsp.\oss

\myuppar{Families of\dss differential\dss boundary\sss problems.}
Suppose\sss that\sss the manifold\sss $X\off =\off X_{\dff w}$\nsp,\oss 
the bundle\sss $E\off =\off E_{\dff w}$\nsp,\oss
the operator $P\off =\off P_{\dff w}$\nsp,\oss
the extension of\dss $P$\sss to $\widehat{X}$\nnsp,\oss etc.\qss
continuously depend on a parameter $w\qff \in\qff W$ as in\qss \cite{i2}.\oss
The explicit\sss definition of\dss the operators\sss
$K_{\dff +}\off =\off K_{\dff +\dff w}$ shows\sss that\sss they
continuously depend on $w$ in\sss the norm\sss topology.\oss
This\sss implies\sss that\sss the kernels 
$Z_{\trf 0}\off =\off Z_{\trf 0\dff w}\off =\off \kernel\fff T^{\dff *}$
and\sss the operators $\varphi\off =\off \varphi_{\dff w}$ 
continuously depend on $w$ in\sss the norm\sss topology.\oss
In\sss turn,\oss this implies\sss that\sss
$\bm{\gamma}_{w}\trf(\dff 0\dff)$
and\sss $M_{\dff w}\trf(\dff 0\dff)$
also continuously depend on $w$ in\sss the norm\sss topology.\oss
By applying\trs Calder\'{o}n's\dss method\sss to\sss the operators 
$P_{\dff w}\qff +\qff i$\sss
instead of\dss $P_{\dff w}$ we see\sss that\sss the kernels\sss 
$\mathcal{K}_{\dff w\dff -}
\off =\off
\kernel\fff (\trf T_w^{\dff *}\qff +\qff i\qff)$
continuously depend on $w$\nnsp.\oss
This verifies\sss the continuity assumption of\trs Section\qss \ref{families}\qss
in\sss the present\sss situation.\oss

%\newpage
\mysection{Dirac-like\qss boundary\qss problems}{dirac-like}

\myuppar{Graded\dss boundary\sss problems.}
Let\sss $X\fff,\qff Y,\qff E\fff,\qff F$\sss be as 
in\dss Section\qss \ref{differential-boundary-problems}.\oss
Suppose\sss that\sss the decomposition\sss
$E\trf|\trf Y\off =\off F\dff \oplus\dff F$\sss
extends\sss to a decomposition\sss
$E\off =\off G\dff \oplus\dff G$\sss
over\sss the whole manifold\sss $X$\nnsp.\oss
So,\oss in\sss particular{},\pss $F\off =\off G\trf|\trf Y$\dnsp.\oss
Let\sss $P$\sss be
a formally self-adjoint\sss elliptic differential\sss operator\sss 
of\dss order $1$ acting on sections of\dss $E$
and\sss having\sss the form\vspace{3pt}\vspace{-1.5pt}
\begin{equation}
\label{p}
\quad
P
\off =\off\dff
\begin{pmatrix}
\off 0 &
R' \qff\off
\vspace{4.5pt} \\
\off\dff R &
0 \qff\off 
\end{pmatrix}
\off
\end{equation}

\vspace{-12pt}\vspace{3pt}\vspace{-1.55pt}
with respect\sss to\sss the decomposition\sss $E\off =\off G\dff \oplus\dff G$\nnsp,\oss
i.e.\qss being\qss \emph{odd}\pss with respect\sss to\sss this decomposition.\oss
Suppose\sss that\sss the coefficient\sss
$\Sigma$\sss of\dss the normal\sss derivative\sss in\sss $P$\sss
has\sss the form\vspace{3pt}\vspace{-1.5pt}
\[
\quad
\Sigma
\off =\off\dff
\begin{pmatrix}
\off 0 &
1 \qff\off
\vspace{4.5pt} \\
\off\dff 1 &
0 \qff\off 
\end{pmatrix}
\off
\]

\vspace{-12pt}\vspace{3pt}\vspace{-1.5pt}
with respect\sss to\sss the decomposition\sss
$E\trf|\trf Y\off =\off F\dff \oplus\dff F$\dnsp.\oss
This convention agrees with\qss \cite{i2},\oss Section\qss 15,\oss
and does not\sss agree with\sss the convention of\trs 
Section\qss \ref{differential-boundary-problems}.\oss
Below we will\dss pass\sss to another direct\sss sum decomposition of\dss
$E$ which\sss will\dss match\dss Section\qss \ref{differential-boundary-problems}.\oss
Let\sss $\gamma$\sss be\sss the\sss trace operator as in\dss
Section\qss \ref{differential-boundary-problems},\oss
and\dss let\sss $\pi_{\dff 0}$ and\sss $\pi_{\dff 1}$\sss
be its compositions with\sss the projections\sss
$\Pi_{\dff 0}\dff,\pff \Pi_{\dff 1}$\sss
onto\sss the first\sss and\sss the second summands of\dss
the decomposition\vspace{1.5pt}
\[
\quad 
H_{\dff 1/2}\dff(\trf Y\fff,\qff E\trf|\trf Y\trf)
\off =\off
H_{\dff 1/2}\dff(\trf Y\fff,\qff F\trf)
\qff \oplus\qff 
H_{\dff 1/2}\dff(\trf Y\fff,\qff F\trf)
\qff,\quad
\]

\vspace{-12pt}\vspace{1.5pt}
Let\sss $f\dff \colon\dff F\qff \ttoo\qff F$\sss be a bundle map
and\dss let\sss
$\pi_{\fff f}
\off =\off
\pi_{\dff 1}
\qff -\qff 
f\dff \circ\dff \pi_{\dff 0}$\nsp.\oss
We would\dss like\sss to combine $P$ with\sss the boundary condition\sss
$\pi_{\fff f}\off =\off 0$\nnsp,\oss
i.e.\qss to consider\sss the unbounded operator\sss $P\trf|\trf \kernel\fff \pi_{\fff f}$\nsp.\oss
If\dss $f$\sss is\dss skew-adjoint,\oss then\sss this boundary condition\dss
is\dss self-adjoint.\oss
See\qss \cite{i2},\oss the discussion before\dss Lemma\qss 15.3.\oss
Let\sss us write\sss the symbol\sss of\dss the operator $P$ in\sss the form\qss (\ref{symbol}).\oss
Since $P$\sss is\dss odd,\oss the operators\sss
$\tau_{\fff u}\off =\off \tau_{\fff u}\trf(\trf y\trf)$
have\sss the form\vspace{3pt}
\[
\quad
\tau_u
\off =\off\dff
\begin{pmatrix}
\off 0  &
\bm{\tau}_u^{\dff *} \dff\off
\vspace{4.5pt} \\
\off \bm{\tau}_u &
0 \dff\off 
\end{pmatrix}
\off
\]

\vspace{-12pt}\vspace{3pt}
for some operators\sss
$\bm{\tau}_u\dff \colon\dff F_{\dff y}\qff \ttoo\qff F_{\dff y}$\nnsp.\oss
As\sss in\qss \cite{i2},\oss Section\qss 15,\oss
we will\sss say\sss that $f$ is\qss \emph{equivariant}\pss if\dss endomorphisms\sss
$f_{\dff y}\dff \colon\dff F_{\dff y}\qff \ttoo\qff F_{\dff y}$
induced\dss by\sss $f$\sss commute with\sss $\bm{\tau}_u$\sss
for every\sss $y\fff,\qff u$\nnsp.\oss
If\dss $f$\sss is\dss equivariant,\oss then\sss the boundary condition\sss
$\pi_{\fff f}\off =\off 0$\sss
is\dss a\dss Shapiro--Lopatinskii\dss boundary condition,\oss
as also\sss the boundary\sss condition\sss
$\pi_{\dff 1}
\qff +\qff 
f\dff \circ\dff \pi_{\dff 0}
\off =\off
0$\nnsp.\oss
See\qss \cite{i2},\oss Lemma\qss 15.3.\oss
For\sss the rest\sss of\dss this section\sss we will\sss
assume\sss that\sss $f$\sss is\dss skew-adjoint\sss and equivariant.\oss

In\sss general,\pss $P\trf|\trf \kernel\fff \pi_{\fff f}$\sss is\dss not\sss invertible
and\dss hence cannot\sss be used as a reference operator{}.\oss
But\sss under natural\sss assumptions
a simple modification of\dss $P\trf|\trf \kernel\fff \pi_{\fff f}$\sss is\dss invertible.\oss
Let\vspace{3pt}
\[
\quad
\varepsilon
\off =\off\dff
\begin{pmatrix}
\off 1 &
0 \qff\off
\vspace{4.5pt} \\
\off\dff 0 &
-\qff 1 \qff\off 
\end{pmatrix}
\off
\]

\vspace{-12pt}\vspace{0pt}
be\sss the endomorphism of\dss $E\off =\off G\dff \oplus\dff G$ defined\sss
by\sss the above matrix.\oss
We will\sss consider $\varepsilon$ as a differential\sss operator of\dss order $0$\nnsp.\oss
It\sss tuns out\sss that\sss if\dss the endomorphism $i f$\sss is\dss positive definite,\oss
then\sss $(\trf P\qff +\qff \varepsilon\trf)\trf|\trf \kernel\fff \pi_{\fff f}$\sss
is\dss an\sss isomorphism\sss
$\kernel\fff \pi_{\fff f}\qff \ttoo\qff H_{\dff 0}$\nsp.\oss
Similarly,\oss if\dss the endomorphism $i f$\sss is\dss negative definite,\oss
then\sss $(\trf P\qff -\qff \varepsilon\trf)\trf|\trf \kernel\fff \pi_{\fff f}$\sss
is\dss an\sss isomorphism.\oss
See\qss \cite{i2},\oss Theorem\qss 15.1.\oss
A similar result\sss holds also for closed\sss manifolds,\oss
as\sss the following\dss lemma shows.\oss

\mypar{Lemma.}{closed-invertible}
\emph{Suppose\sss that\sss $\widehat{P}$\sss is\dss
a formally self-adjoint\sss elliptic differential\sss operator\sss 
of\dss order $1$ acting in\sss a bundle\sss
$E\off =\off G\dff \oplus\dff G$\sss over a closed\dss manifold\dss $\widehat{X}$\nnsp.\oss
If\qss $\widehat{P}$ is\dss odd,\oss
then\sss $\widehat{P}\qff -\qff \varepsilon$\sss induces an\sss isomorphism\dss
$H_{\dff 1}\dff(\trf X\fff,\qff E\trf)
\qff \ttoo\qff 
H_{\dff 0}\dff(\trf X\fff,\qff E\trf)$\nnsp.\oss}

\proof
Under\sss these assumptions\sss the operator\sss
$H_{\dff 1}\dff(\trf X\fff,\qff E\trf)
\qff \ttoo\qff 
H_{\dff 0}\dff(\trf X\fff,\qff E\trf)$\sss
induced\dss by\sss $\widehat{P}\qff -\qff \varepsilon$\sss 
is\trs Fredholm\dss and\dss is\dss self-adjoint\sss
as an\sss unbounded operator\sss in\sss
$H_{\dff 0}\dff(\trf X\fff,\qff E\trf)$\nnsp.\oss
Therefore\sss it\dss is\dss sufficient\sss to prove\sss
that\sss its\sss kernel\dss is $0$\nnsp.\oss
If\dss we consider $\widehat{P}$ as an\sss unbounded operator\sss in\sss
$H_{\dff 0}\dff(\trf X\fff,\qff E\trf)$
and\sss write it\sss in\sss the form\qss (\ref{p}),\oss
then\sss $R'$\sss will\dss be equal\dss to\sss the adjoint\sss
operator $R^{\fff *}$ of\dss $R$\nnsp.\oss 
Therefore,\oss if\dss 
$(\trf u\fff,\qff v\trf)
\qff \in\qff 
\kernel\fff (\qff \widehat{P}\qff -\qff \varepsilon\trf)$\nnsp,\oss
then\sss
$-\qff u\qff +\qff R^{\fff *}\fff v\off =\off 0$
and\sss
$R\dff u\qff +\qff v\off =\off 0$\nnsp.\oss
This implies\sss that\sss
$u\qff +\qff R^{\fff *}\fff R\dff u\off =\off 0$
and\dss hence\sss
$0
\off =\off
\sco{\dff u\qff +\qff R^{\fff *}\fff R\dff u\fff,\qff u\dff}
\off =\off
\sco{\dff u\fff,\qff u\dff}
\qff +\qff
\sco{\dff R\dff u\fff,\qff R\dff u\dff}$\nnsp.\oss
It\sss follows\sss that\sss $\sco{\dff u\fff,\qff u\dff}\off =\off 0$\sss
and\dss hence\sss $u\off =\off 0$\nnsp.\oss
This proves\sss that\sss
$\kernel\fff (\qff \widehat{P}\qff -\qff \varepsilon\trf)
\off =\off 
0$\nnsp.\oss  \eproof

\myuppar{Changing\sss the decomposition $E\trf|\trf Y\off =\off F\dff \oplus\dff F$\dnsp.}
The\dss Lagrange\dss identity for\sss $P$\sss is\vspace{3pt}
\[
\quad
\sco{\dff P\dff u\fff,\qff v \dff}_{\dff 0}
\qff -\qff
\sco{\dff u\dff,\qff P\dff v \dff}_{\dff 0}
\off =\off
\bsco{\dff i\trf \pi_{\dff 1}\dff u\dff,\qff \pi_{\dff 0}\dff v\dff}_{\dff \partial}
\qff -\qff
\bsco{\dff \pi_{\dff 0}\dff u\dff,\qff i\trf \pi_{\dff 1}\dff v\dff}_{\dff \partial}
\pff.
\]

\vspace{-12pt}\vspace{3pt}
Let\sss us define\sss the new boundary operators\sss $\gamma_0$ 
and\sss $\gamma_1$ as\vspace{4.5pt}
\[
\quad
\gamma_0
\off =\off 
(\trf \pi_{\dff 1}\qff -\qff i\dff \pi_{\dff 0}\trf)\bigl/\sqrt{2} 
\quad
\mbox{and}\quad
\gamma_1
\off =\off 
(\trf \pi_{\dff 1}\qff +\qff i\dff \pi_{\dff 0}\trf)\bigl/\sqrt{2}
\off.
\]

\vspace{-12pt}\vspace{4.5pt}
The calculation\vspace{3pt}
\[
\quad
\bsco{\dff \pi_{\dff 1}\dff u\qff +\qff i\dff \pi_{\dff 0}\dff u\dff,\qff
\pi_{\dff 1}\dff v\qff -\qff i\dff \pi_{\dff 0}\dff v\dff}_{\dff \partial}
\pff -\pff
\bsco{\dff \pi_{\dff 1}\dff u\qff -\qff i\dff \pi_{\dff 0}\dff u\dff,\qff
\pi_{\dff 1}\dff v\qff +\qff i\dff \pi_{\dff 0}\dff v\dff}_{\dff \partial}
\]

\vspace{-33pt}
\[
\quad
=\off
\bsco{\dff \pi_{\dff 1}\dff u\dff,\qff -\qff i\dff \pi_{\dff 0}\dff v\dff}_{\dff \partial}
\qff +\qff
\bsco{\dff i\dff \pi_{\dff 0}\dff u\dff,\qff \pi_{\dff 1}\dff v\dff}_{\dff \partial}
\qff -\qff
\bsco{\dff \pi_{\dff 1}\dff u\dff,\qff i\dff \pi_{\dff 0}\dff v\dff}_{\dff \partial}
\qff -\qff
\bsco{\dff -\qff i\dff \pi_{\dff 0}\dff u\dff,\qff \pi_{\dff 1}\dff v\dff}_{\dff \partial}
\]

\vspace{-33pt}
\[
\quad
=\off
2\dff
\bsco{\dff i\dff \pi_{\dff 0}\dff u\dff,\qff \pi_{\dff 1}\dff v\dff}_{\dff \partial}
\qff -\qff
2\dff
\bsco{\dff \pi_{\dff 1}\dff u\dff,\qff i\dff \pi_{\dff 0}\dff v\dff}_{\dff \partial}
\off =\off
2\dff
\bsco{\dff i\dff \pi_{\dff 1}\dff u\dff,\qff \pi_{\dff 0}\dff v\dff}_{\dff \partial}
\qff -\qff
2\dff
\bsco{\dff \pi_{\dff 0}\dff u\dff,\qff i\dff \pi_{\dff 1}\dff v\dff}_{\dff \partial}
\pff
\]

\vspace{-12pt}\vspace{3pt}
shows\sss that\sss we can rewrite\sss the\dss Lagrange\dss identity\sss
in\sss the standard\sss form\qss (\ref{p-lagrange}).\oss 
Clearly,\oss  
$\gamma_0$ and\sss $\gamma_1$ are\sss the boundary operators related\sss
as in\trs Section\qss \ref{differential-boundary-problems}\qss
with a unique decomposition $E\trf|\trf Y\off =\off F\fff'\dff \oplus\dff F\fff'$\dnsp.\oss
The bundle\sss $F\fff'$\sss is\dss canonically\sss isomorphic\sss to $F$\dnsp,\oss
but\sss the decomposition of\dss $E\trf|\trf Y$\sss 
is\dss
different\sss from\sss the original\sss one.\oss
In\sss terms of\dss this new decomposition\sss 
$i\trf \Sigma$\sss takes\sss the form\qss (\ref{i-sigma})\qss
and\sss the boundary condition\sss
$\pi_{\dff 1}\qff -\qff f\dff \circ\dff \pi_{\dff 0}\off =\off 0$\sss
takes\sss the form\vspace{1.5pt}
\[
\quad
\gamma_1\qff +\qff \gamma_0
\off =\off
f\dff \circ\qff\left(\qff
\frac{\gamma_1\qff -\qff \gamma_0}{i}
\qff\right)
\pff,
\]

\vspace{-12pt}\vspace{1.5pt}
or{},\oss equivalently,\oss either\sss 
$\left(\trf f\qff -\qff i\trf\right)\dff \circ\dff \gamma_1
\off =\off
\left(\trf f\qff +\qff i\trf\right)\dff \circ\dff \gamma_0$\nsp,\oss 
or \vspace{1.5pt}
\[
\quad
\gamma_1
\off =\off
\frac{f\qff +\qff i}{f\qff -\qff i}
\qff \circ\qff
\gamma_0
\pff.
\]

\vspace{-12pt}\vspace{1.5pt}
In\sss general,\oss this boundary condition\dss is\dss defined\dss by a\sss
relation\sss between\sss $\gamma_0\dff,\pff \gamma_1$\nsp.\oss
In\sss particular{},\oss the boundary condition\sss
$\gamma_0\off =\off 0$\sss corresponds\sss to\sss $f$\sss being\sss the multiplication\sss
by $i$\sss and\dss is\dss defined\dss by a relation.\oss
The boundary condition\sss $\gamma_1\off =\off 0$\sss corresponds\sss 
to\sss $f$\sss being\sss the multiplication\sss by $-\qff i$\nnsp.\oss
Both of\dss them satisfy\sss the\dss Shapiro--Lopatinskii\dss condition.\oss

\myuppar{Replacing\sss $P$\sss by\sss $P\qff -\qff \varepsilon$\nnsp.}
Let\sss $f_{\dff 0}$\sss be\sss the operator of\dss multiplication\sss by\sss $i$\nnsp.\oss
The corresponding boundary condition\dss is\sss
$\pi_{\dff 1}\qff -\qff i\dff \pi_{\dff 0}\off =\off 0$\nnsp,\oss
or{},\oss equivalently,\pss
$\gamma_0\off =\off 0$\nnsp.\oss
Since\sss the operator\sss 
$i\fff f_{\dff 0}\off =\off -\qff \id$\sss 
is\dss negative definite,\oss the operator\sss 
$A\off =\off (\trf P\qff -\qff \varepsilon\trf)\trf|\trf \kernel\fff \gamma_0$\sss
is\dss an\sss isomorphism\sss
$\kernel\fff \gamma_0\qff \ttoo\qff H_{\dff 0}$\sss
and\dss hence can\sss be used as a reference operator.\oss
This amounts\sss to replacing\sss $P$\sss by\sss $P\qff -\qff \varepsilon$\nnsp.\oss
Such a replacement\sss does not\sss affect\sss the\dss Lagrange\dss identity\sss
and\dss the boundary conditions defined\sss in\sss terms of\dss $\gamma_0\dff,\pff \gamma_1$\nsp.\oss
In\sss particular,\oss a boundary condition\dss is\dss self-adjoint\sss or\sss has\sss
the\dss Shapiro--Lopatinskii\dss property\sss for\sss $P$\sss if\dss and\dss only\trs if\dss
it\dss has\sss the same property\sss for\sss $P\qff -\qff \varepsilon$\nnsp.\oss
At\sss the same\sss time,\oss when our manifold,\oss operators,\oss etc.\dss
depend on a parameter,\oss replacing\sss $P$\sss by\sss $P\qff -\qff \varepsilon$\sss
does not\sss affect\sss the analytical\dss index\sss because $P$ can\sss be connected\sss
with\sss $P\qff -\qff \varepsilon$\sss by\sss the homotopy\sss
$P\qff -\qff t\trf \varepsilon$\nnsp,\qss $t\qff \in\qff [\trf 0\fff,\qff 1\trf]$\nnsp.\oss
By\sss these reasons we can\sss work with $P\qff -\qff \varepsilon$ instead of\dss $P$\dnsp.\oss

\myuppar{The reduced\dss boundary\sss triplet.}
As suggested\dss by\sss the previous subsection,\oss
we will\dss take\sss the operator\sss
$A\off =\off (\trf P\qff -\qff \varepsilon\trf)\trf|\trf \kernel\fff \gamma_0$\sss
as\sss the reference operator{}.\oss
Naturally,\oss we will\sss also\sss take\sss the operators\sss
$\gamma_0\dff,\pff \gamma_1\dff \colon\dff H_{\dff 1}\qff \ttoo\qff K$\sss
as\sss the boundary operators.\oss
Then all\sss assumptions of\trs Section\qss \ref{abstract-index}\qss hold.\oss
In\sss particular{},\oss the operators\sss $M\trf(\dff 0\dff)$\sss
and\sss $M\off =\off M\trf(\dff 0\dff)\trf|\trf K$\sss are defined,\oss
and\sss $M$\sss is\dss a self-adjoint\sss densely defined operator\sss
from\sss $K\fff'$\sss to\sss $K$\nnsp.\oss

\myuppar{The boundary conditions\sss in\sss terms of\dss the reduced\dss boundary\sss triplet.}
The boundary conditions corresponding\sss to $f$ are defined\sss in\sss
terms of\sss $\gamma_0\fff,\qff \gamma_1$\sss by an obvious\sss relation\sss
$\mathcal{F}\qff \subset\qff K^{\dff \partial}\dff \oplus\dff K^{\dff \partial}$\nsp\dnsp.\oss
In\sss terms of\dss the reduced\dss boundary\sss triplet\sss these boundary conditions\sss
take\sss the form\vspace{3pt}
\[
\quad
\Lambda'\dff \oplus\dff \Lambda^{\fff -\dff 1}\qff
\bigl(\qff
\mathcal{F}\trf|\trf K\qff -\qff M
\qff\bigr)
\pff,
\]

\vspace{-12pt}\vspace{3pt}
where\sss $\Lambda\dff \colon\dff K^{\dff \partial}\qff \ttoo\qff K$\sss
and\sss $\Lambda'\dff \colon\dff K\fff'\qff \ttoo\qff K^{\dff \partial}$ 
are\sss the operators from\sss the\sss theory of\trs Gelfand\trs triples.\oss
As\dss is\dss well\dss known,\oss both\sss $\Lambda$ and\sss $\Lambda'$\sss
are pseudo-differential\sss operators of\dss order $-\qff 1/2$\sss with\sss
the symbol\sss equal\sss to\sss the identity\sss  
over\sss the unit\sss sphere bundle of\dss $Y$\dnsp.\oss
When $f\qff -\qff i$\sss is\dss invertible,\oss
these boundary conditions can be written as\sss the equation\vspace{3pt} 
\[
\quad
\overline{\bm{\Gamma}}_1
\off =\off
\Lambda^{\fff -\dff 1}\qff \circ\pff  
\left(\qff
\frac{f\qff +\qff i}{f\qff -\qff i}
\qff -\qff
M
\qff\right)
\qff \circ\qff
(\trf \Lambda'\trf)^{\fff -\dff 1}\qff \circ\qff  
\overline{\Gamma}_0
\pff.
\]

\vspace{-12pt}\vspace{3pt}
In\sss particular{},\oss the boundary conditions are given\sss by a
pseudo-differential\sss operator of\dss order $1$\nnsp.\oss
Informally,\oss one can use\sss this form of\dss boundary conditions
even when $f\qff -\qff i$\sss is\dss not\sss invertible.\oss
In\sss fact,\oss the case when $f$ has only $i$ and $-\qff i$ as eigenvalues\dss
is\dss the most\sss important\sss one.\oss

\myuppar{The decomposition of\dss $F$ defined\dss by $f$\dnsp.}
Since the bundle map $f$\sss is\dss skew-adjoint,\oss it\sss defines\sss
the decomposition\sss 
$F
\off =\off 
\mathcal{L}_{\dff +}\dff(\trf f\trf)
\dff \oplus\dff
\mathcal{L}_{\dff -}\dff(\trf f\trf)$\sss
into subbundles generated\dss by\sss eigenvectors of\dss $f$\sss
corresponding\sss to eigenvalues $\lambda$\sss with\sss
$\image\fff \lambda\qff >\qff 0$
and\sss
$\image\fff \lambda\qff <\qff 0$ respectively.\oss
Clearly,\pss $i f$\sss is\dss negative definite on\sss $\mathcal{L}_{\dff +}\dff(\trf f\trf)$
and\dss is\dss positive definite on\sss $\mathcal{L}_{\dff -}\dff(\trf f\trf)$\nnsp.\oss
Since $f$\sss is\dss equivariant,\oss the subbundles\sss
$\mathcal{L}_{\dff +}\dff(\trf f\trf)\fff,\off \mathcal{L}_{\dff -}\dff(\trf f\trf)$\sss
are invariant\sss under\sss the operators\sss $\bm{\tau}_u$\nsp.\oss
Let\sss $\bm{\tau}_u^{\dff +}\dff,\off \bm{\tau}_u^{\dff -}$\sss
be\sss the operators induced\dss by\sss $\bm{\tau}_u$\sss
in\sss the fibers of\dss
$\mathcal{L}_{\dff +}\dff(\trf f\trf)\fff,\off \mathcal{L}_{\dff -}\dff(\trf f\trf)$\sss
respectively.\oss
The families of\dss operators\sss $\bm{\tau}_u^{\dff +}\dff,\off \bm{\tau}_u^{\dff -}$\sss
can\sss be considered as symbols of\dss some pseudo-differential\sss operators
$\bm{\tau}^{\dff +}\dff,\off \bm{\tau}^{\dff -}$\sss
of\dss order $1$ acting\sss in\sss bundles
$\mathcal{L}_{\dff +}\dff(\trf f\trf)\fff,\off \mathcal{L}_{\dff -}\dff(\trf f\trf)$\sss
respectively.\oss

\myuppar{Dirac-like\sss boundary\sss problems.}
Suppose\sss that\sss $P$ and $f$ define a\qss \emph{Dirac-like}\pss
boundary problem\sss in\sss the sense of\pss $\cite{i2}$\nnsp,\oss
i.e.\qss that\sss the operators\sss $\bm{\tau}_u$ are skew-adjoint.\oss
A routine calculation shows\sss that\sss
in\sss the decomposition\sss
$E\trf|\trf Y\off =\off F\fff'\dff \oplus\dff F\fff'$ 
the operators
$\rho_{\fff u}
\off =\off
\sigma_y^{\dff -\dff 1}\dff \tau_{\fff u}$
have\sss the form\vspace{0.75pt} 
\[
\quad
\rho_{\fff u}
\off =\off\dff
\begin{pmatrix}
\off 0  &
-\qff \bm{\tau}_u \dff\off
\vspace{4.5pt} \\
\off -\qff \bm{\tau}_u &
0 \dff\off 
\end{pmatrix}
\off.
\]

\vspace{-12pt}\vspace{3pt}
Let\sss 
$\num{\bm{\tau}_u}
\off =\off
(\trf \bm{\tau}_u^{\dff *}\dff\bm{\tau}_u\trf)^{\dff 1/2}$
and\dss let\sss $\bm{\upsilon}_{\fff u}$\sss
be such\dss that\sss
$i\trf \bm{\tau}_u
\off =\off 
-\qff \bm{\upsilon}_{\fff u}\qff \num{\bm{\tau}_u}$\nsp.\oss

\mypar{Lemma.}{dirac-m}
\emph{The symbol\sss of\dss $M$\sss is\dss
equal\dss to\sss the bundle map defined\dss by operators\dss 
$\bm{\upsilon}_{\fff u}$\nsp.\oss}

\proof
Lemma\qss \ref{closed-invertible}\qss implies\sss that\sss the results of\trs
Section\qss \ref{differential-boundary-problems}\qss apply\sss to\sss the present\sss situation.\oss
By\sss the discussion at\sss the end of\trs Section\qss \ref{differential-boundary-problems},\oss
it\dss is\dss sufficient\sss to prove\sss that\sss 
$\mathcal{L}_{\dff -}\dff(\trf \rho_{\fff u}\trf)$\sss
is\dss equal\dss to\sss the graph of\dss $\bm{\upsilon}_{\fff u}$\nsp.\oss
If\dss $a$\sss is\dss an eigenvector\sss of\dss $\bm{\tau}_u$\sss
with an eigenvalue\sss $i\trf \lambda$\nnsp,\oss
then\vspace{1.5pt} 
\[
\quad
\rho_{\fff u}\trf
\bigl(\trf 
a\fff,\pff \bm{\upsilon}_{\fff u}\dff(\trf a\trf)\trf
\bigr)
\off =\off
\bigl(\trf 
-\qff \bm{\tau}_u\dff \bm{\upsilon}_{\fff u}\dff(\trf a\trf)\fff,\pff
-\qff \bm{\tau}_u\dff(\trf a\trf)
\trf\bigr)
\]

\vspace{-34.5pt}
\[
\quad
\phantom{\rho_{\fff u}\trf
\bigl(\trf 
a\fff,\pff \bm{\upsilon}_{\fff u}\dff(\trf a\trf)\trf
\bigr)
\off }
=\off
\bigl(\trf -\qff i\trf \num{\lambda}\dff a\fff,\pff
-\qff i\trf \lambda\dff a
\trf\bigr)
\off =\off
\bigl(\trf -\qff i\trf \num{\lambda}\dff a\fff,\pff
-\qff i\trf \num{\lambda}\dff \bm{\upsilon}_{\fff u}\dff(\trf a\trf)
\trf\bigr)
\]

\vspace{-34.5pt}
\[
\quad
\phantom{\rho_{\fff u}\trf
\bigl(\trf 
a\fff,\pff \bm{\upsilon}_{\fff u}\dff(\trf a\trf)\trf
\bigr)
\off =\off
\bigl(\trf -\qff i\trf \num{\lambda}\dff a\fff,\pff
-\qff i\trf \lambda\dff a
\trf\bigr)
\off }
=\off
-\qff i\trf \num{\lambda}\qff 
\bigl(\trf a\fff,\pff
\bm{\upsilon}_{\fff u}\dff(\trf a\trf)
\trf\bigr)
\]

\vspace{-12pt}\vspace{1.5pt}
and\dss hence\sss
$(\trf a\fff,\pff \bm{\upsilon}_{\fff u}\dff(\trf a\trf)\trf)
\qff \in\qff
\mathcal{L}_{\dff -}\dff(\trf \rho_{\fff u}\trf)$\nnsp.\oss
Therefore
$\mathcal{L}_{\dff -}\dff(\trf \rho_{\fff u}\trf)$\sss
is\dss equal\dss to\sss the graph of\dss $\bm{\upsilon}_{\fff u}$\nsp.\oss  \eproof

\myuppar{Families of\trs Dirac-like\sss boundary\sss problems.}
Suppose\sss that\sss the manifold\sss $X\off =\off X_{\dff w}$\nsp,\oss 
the bundle\sss $E\off =\off E_{\dff w}$\nsp,\oss
the operator $P\off =\off P_{\dff w}$\nsp,\oss the bundle map\sss
$f\off =\off f_{\dff w}$\sss etc.\qss
continuously depend on a parameter $w\qff \in\qff W$ as in\qss \cite{i2},\oss
and\sss for each value of\sss $w$\sss have all\dss the properties assumed above.\oss
The subscript\sss $w$ in $f_{\dff w}$ should\sss not\sss be confused\sss with\sss the
subscript\sss $y$ used\sss above:\qss
$w\qff \in\qff W$\nnsp,\oss but\sss $y\qff \in\qff Y$\dnsp.\oss
For each $w$\sss let\sss $A_{\dff w}$\sss be\sss the self-adjoint\sss operator
defined\dss by\sss $P_{\dff w}$ 
and\sss the boundary condition corresponding\sss to $f_{\dff w}$\nsp.\oss 
Clearly,\oss the bundles\sss
$\mathcal{L}_{\dff +}\dff(\trf f_{\dff w}\trf)\fff,\off 
\mathcal{L}_{\dff -}\dff(\trf f_{\dff w}\trf)$\sss
continuously depend on $w$\sss
and\sss one can define continuous families 
$\bm{\tau}^{\dff +}_{\dff w}\trf,\off \bm{\tau}^{\dff -}_{\dff w}$\sss
of\dss skew-adjoint\dss pseudo-differential\sss operators\sss
of\dss order $1$\nnsp.\oss

\mypar{Theorem.}{index-theorem}
\emph{The analytical\dss index of\dss the family\sss
$A_{\dff w}\dff,\pff w\qff \in\qff W$\sss
is\dss equal\dss to\sss the analytical\dss index of\dss
the family\sss
$i\trf \bm{\tau}^{\dff -}_{\dff w}\dff,\pff w\qff \in\qff W$\nnsp,\oss
as also of\dss the family\sss
$-\qff i\trf \bm{\tau}^{\dff +}_{\dff w}\dff,\pff w\qff \in\qff W$\nnsp.\oss}

\proof
As we pointed out\sss above,\oss we can replace\sss the operators\sss 
$P_{\dff w}$\sss by\sss the operators\sss
$P_{\dff w}\qff -\qff \varepsilon$\nnsp.\oss
Also,\oss without\sss affecting\sss the analytical\dss indices,\oss
we can\sss deform\sss the bundle maps\sss $f_{\dff w}$\sss
to skew-adjoint\sss bundle maps having only $i$ and $-\qff i$ as eigenvalues.\oss
Then\vspace{-1.5pt} 
\[
\quad
\frac{f_{\dff w}\qff +\qff i}{f_{\dff w}\qff -\qff i}
\pff
\]

\vspace{-12pt}\vspace{-1.5pt}
is\dss equal\dss to $0$\sss on\sss
$\mathcal{L}_{\dff -}\dff(\trf f_{\dff w}\trf)$
and\sss to $\infty$\sss on\sss
$\mathcal{L}_{\dff +}\dff(\trf f_{\dff w}\trf)$\nnsp.\oss
Of\dss course,\oss the boundary condition\sss
$\gamma_1\off =\off \infty\qff \gamma_0$\sss
should\dss be interpreted as\sss
$\gamma_0\off =\off 0$\nnsp.\oss
Next,\oss let\sss us\sss pass\sss to\sss
the reduced\dss boundary\sss triplets
and\dss rewrite\sss the boundary conditions 
in\sss the form\vspace{0.5pt}%\vspace{-2pt}
\[
\quad
\overline{\bm{\Gamma}}_1
\off =\off
\Lambda^{\fff -\dff 1}\qff \circ\pff  
\left(\qff
\frac{f_{\dff w}\qff +\qff i}{f_{\dff w}\qff -\qff i}
\qff -\qff
M_{\dff w}
\qff\right)
\qff \circ\qff
(\trf \Lambda'\trf)^{\fff -\dff 1}\qff \circ\qff  
\overline{\Gamma}_0
\pff,
\]

\vspace{-12pt}\vspace{0.5pt}%\vspace{-0.25pt}
where we omitted\dss the dependence on $w$ of\qss $\overline{\Gamma}_0\dff,\off \overline{\bm{\Gamma}}_1$ 
and\sss $\Lambda\fff,\qff \Lambda'$\nnsp.\oss
The continuity assumptions of\trs Section\qss \ref{families}\qss were
verified at\sss the end of\trs Section\qss \ref{differential-boundary-problems}.\oss
Also,\oss the operators $A_{\dff w}$ are operators with compact\sss resolvent.\oss
Hence\sss the results of\trs Section\qss \ref{families}\qss apply.\oss
It\sss follows\sss that\sss the analytical\dss index of\dss the family\sss
$A_{\dff w}\dff,\pff w\qff \in\qff W$\sss
is\dss equal\dss to\sss the analytical\dss index of\dss the family
of\dss relations\vspace{1.5pt} 
\[
\quad
\Lambda^{\fff -\dff 1}\qff \circ\pff  
\left(\qff
\frac{f_{\dff w}\qff +\qff i}{f_{\dff w}\qff -\qff i}
\qff -\qff
M_{\dff w}
\qff\right)
\qff \circ\qff
(\trf \Lambda'\trf)^{\fff -\dff 1}
\qff,\quad
w\qff \in\qff W
\pff.
\]

\vspace{-12pt}\vspace{1.5pt}
This family\dss is\dss equal\dss to\sss the direct\sss sum of\dss two families,\oss
one in\sss the bundles\sss 
$\mathcal{L}_{\dff +}\dff(\trf f_{\dff w}\trf)$
and\dss the other one in\sss the bundles\sss 
$\mathcal{L}_{\dff -}\dff(\trf f_{\dff w}\trf)$\nnsp.\oss
The family\sss in\sss the bundles $\mathcal{L}_{\dff +}\dff(\trf f_{\dff w}\trf)$\sss
is\vspace{1.5pt} 
\[
\quad
\Lambda^{\fff -\dff 1}\qff \circ\pff  
\bigl(\qff
\infty
\qff -\qff
M^{\dff +}_{\dff w}
\qff\bigr)
\qff \circ\qff
(\trf \Lambda'\trf)^{\fff -\dff 1}
\qff,\quad
w\qff \in\qff W
\pff,
\]

\vspace{-12pt}\vspace{1.5pt}
where $M^{\dff +}_{\dff w}$\sss is\dss induced\sss by\sss $M_{\dff w}$\nsp.\oss
As a relation,\qss
$\infty\qff -\qff M^{\dff +}_{\dff w}$\sss is\dss equal\dss to $
\infty$\nnsp,\oss
and\sss the boundary conditions defined\dss by\sss the above relations should\dss be interpreted as\dss 
$\overline{\Gamma}_0\off =\off 0$\nnsp.\oss
The index of\dss this\sss family of\dss relations\dss is\dss equal\dss to $0$\nnsp.\oss
The family\sss in\sss the bundles $\mathcal{L}_{\dff -}\dff(\trf f_{\dff w}\trf)$\sss
is\vspace{1.5pt} 
\[
\quad
-\pff
\Lambda^{\fff -\dff 1}\qff \circ\pff  
M^{\dff -}_{\dff w}
\qff \circ\qff
(\trf \Lambda'\trf)^{\fff -\dff 1}
\qff,\quad
w\qff \in\qff W
\pff,
\]

\vspace{-12pt}\vspace{3pt}
where\sss $M^{\dff -}_{\dff w}$\sss is\dss induced\dss by\sss $M_{\dff w}$\nsp.\oss
Lemma\qss \ref{dirac-m}\qss implies\sss that\sss for every $w\qff \in\qff W$\sss 
the symbol\sss of\dss the above operator\dss is\dss equal\dss to\sss 
$-\qff \bm{\upsilon}^{\dff -}_{\dff w}$\nsp,\oss
where\sss $\bm{\upsilon}^{\dff -}_{\dff w}$\sss is\dss
determined\dss by\sss the equality\vspace{3pt} 
\[
\quad
i\qff \bm{\tau}^{\dff -}_{\dff w}
\off =\off
-\qff 
\bm{\upsilon}^{\dff -}_{\dff w}\off \bnum{\bm{\tau}^{\dff -}_{\dff w}}
\pff.
\]

\vspace{-12pt}\vspace{3pt}
Clearly,\oss the family\sss 
$-\qff \bm{\upsilon}^{\dff -}_{\dff w}\dff,\pff w\qff \in\qff W$\sss
is\dss canonically\sss homotopic\sss to\sss the family\sss
$i\qff \bm{\tau}^{\dff -}_{\dff w}\dff,\pff w\qff \in\qff W$\dnsp.\oss
The first\sss statement\sss of\dss the\sss theorem\sss follows.\oss

In order\sss to prove\sss the second statement,\oss 
it\dss is\dss sufficient\sss to prove\sss that\sss the sum
of\dss the analytical\dss indices of\dss the families\sss
$i\trf \bm{\tau}^{\dff -}_{\dff w}\dff,\pff w\qff \in\qff W$
and\sss
$i\trf \bm{\tau}^{\dff +}_{\dff w}\dff,\pff w\qff \in\qff W$\sss
is\dss equal\dss to $0$\nnsp.\oss
This sum\dss is\dss equal\dss to\sss the analytical\dss index of\dss
the direct\sss sum of\dss these families,\oss
i.e.\qss to\sss the analytical\dss index of\dss the family\sss
$i\trf \bm{\tau}_{\dff w}\dff,\pff w\qff \in\qff W$\dnsp.\oss
If\dss we\sss take as\sss $f_{\fff w}$\sss for every $w$ the multiplication\sss by $-\qff i$\nnsp,\oss
then\sss $\mathcal{L}_{\dff -}\dff(\trf f_{\fff w}\trf)\off =\off F_{\fff w}$\sss
and\sss $\bm{\tau}^{\dff -}_{\dff w}\off =\off \bm{\tau}_{\dff w}$\sss 
for every $w$\nnsp.\oss
Since\sss the bundle maps\sss $i f_{\fff w}\off =\off \id$\sss
are positive definite,\oss the index of\dss the corresponding\sss family\sss
$A_{\dff w}\dff,\pff w\qff \in\qff W$\sss
is\dss equal\dss to $0$\nnsp.\oss
At\sss the same,\oss by\sss the already\sss proved\dss first\sss statement\sss
of\dss the\sss theorem,\oss the index of\dss this family\dss is\dss equal\dss to\sss
the index of\dss the family\sss
$i\trf \bm{\tau}^{\dff -}_{\dff w}
\off =\off
i\trf \bm{\tau}_{\dff w}\dff,\pff w\qff \in\qff W$\nnsp.\oss
Therefore\sss the index of\dss the family\sss
$i\trf \bm{\tau}_{\dff w}\dff,\pff w\qff \in\qff W$\sss
is\dss indeed equal\dss to $0$\nnsp.\oss
The second statement\sss of\dss the\sss theorem\sss follows.\oss  \eproof

%\newpage
\mysection{Comparing\qss two\qss boundary\qss conditions}{comparing}

\myuppar{The operator\sss
$S\off =\off T\dff \oplus\dff -\qff T$\dnsp.}
Suppose\sss that\sss the assumptions of\trs Section\qss \ref{abstract-index}\qss hold,\oss
except\sss of\dss the assumptions concerned\sss with\sss the reference 
operator $A$\sss and\dss the\dss Lagrange\dss identity\qss (\ref{lagrange-identity}).\oss
We will\sss assume\sss that\sss a weaker\sss form of\dss the\dss Lagrange\dss identity\sss holds,\oss
namely,\oss that\vspace{3pt}
\[
\quad
\sco{\dff T^{\dff *} u\fff,\qff v \dff}
\pff -\pff
\sco{\dff u\dff,\qff T^{\dff *} v \dff}
\off =\off
\bsco{\dff i\trf \Sigma\qff \gamma\fff u\dff,\qff \gamma\fff v\dff}_{\dff \partial}
\off
\]

\vspace{-12pt}\vspace{3pt}
for every $u\fff,\qff v\qff \in\qff H_{\dff 1}$\sss
and some self-adjoint\sss invertible bounded operator\vspace{1.5pt}
\[
\quad
\Sigma\dff \colon\dff
K^{\dff \partial}\dff \oplus\dff K^{\dff \partial}
\qff \ttoo\qff
K^{\dff \partial}\dff \oplus\dff K^{\dff \partial}
\]

\vspace{-12pt}\vspace{1.5pt}
leaving $K\dff \oplus\dff K$\sss invariant.\oss
Cf.\qss \cite{i2},\oss Section\qss 5.\oss
Let\vspace{3pt}
\[
\quad
\widehat{H}_{\trf 0}
\off =\off
H_{\trf 0}\dff \oplus\dff H_{\trf 0}\qff,\quad
\]

\vspace{-34.5pt}
\[
\quad
\widehat{H}_{\dff 1}
\off =\off 
H_{\dff 1}\dff \oplus\dff H_{\dff 1}\qff,\quad
\]

\vspace{-34.5pt}
\[
\quad
\widehat{K}^{\dff \partial}
\off =\off 
K^{\dff \partial}\dff \oplus\dff K^{\dff \partial}\qff,\quad
\]

\vspace{-34.5pt}
\[
\quad
\widehat{K\dff \oplus\dff K}
\off =\off 
\bigl(\trf K\dff \oplus\dff K\trf\bigr)
\qff \oplus\qff
\bigl(\trf K\dff \oplus\dff K\trf\bigr)
\qff,\quad
\]

\vspace{-34.5pt}
\[
\quad
\widehat{\gamma}
\off =\off 
\gamma\dff \oplus\dff \gamma
\qff \colon\qff
\widehat{H}_{\dff 1}
\qff \ttoo\qff
\widehat{K\dff \oplus\dff K}
\qff.
\]

\vspace{-12pt}\vspace{3pt}
We are interested\sss in\sss the operator\sss
$S\off =\off T\dff \oplus\dff -\qff T$\dnsp.\oss
It\dss is\dss an\sss unbounded operator\sss in\sss $\widehat{H}_{\trf 0}$\sss
having as its domain\sss
$\kernel\fff \widehat{\gamma}
\off \subset\off
\widehat{H}_{\dff 1}$\nnsp.\oss
Let\sss
$\widehat{\Sigma}\off =\off \Sigma\dff \oplus\dff -\qff \Sigma$\nnsp.\oss
The\dss Lagrange\dss identity\sss for $S^*$\sss is\vspace{3pt}\vspace{-0.125pt}
\[
\quad
\sco{\dff S^{*} u\fff,\qff v \dff}
\pff -\pff
\sco{\dff u\dff,\qff S^{*} v \dff}
\off =\off
\bsco{\dff i\qff \widehat{\Sigma}\pff \widehat{\gamma}\dff u\dff,\qff \widehat{\gamma}\dff v\dff}_{\dff \partial}
\off.
\]

\vspace{-12pt}\vspace{3pt}\vspace{-0.125pt}
We will\dss need also\sss the operator\sss $S^*\qff +\qff \varepsilon$\nnsp,\oss
where $\varepsilon$\sss is\dss given\sss by\sss the matrix\vspace{0pt}\vspace{-0.125pt}
\[
\quad
\varepsilon
\off =\off\dff
\begin{pmatrix}
\off 0 \dff &
1 \off
\vspace{4.5pt} \\
\off 1 &
0\qff \off 
\end{pmatrix}
\off
\]

\vspace{-12pt}\vspace{0pt}\vspace{-0.125pt}
in\sss the decomposition\sss
$\widehat{H}_{\trf 0}
\off =\off
H_{\trf 0}\dff \oplus\dff H_{\trf 0}$\nsp.\oss

\myuppar{A\sss boundary\sss condition for operators $S^*$ and\sss $S^*\qff +\qff \varepsilon$\nnsp.}
Let\sss 
$\Pi_{\dff 0}\dff,\pff \Pi_{\dff 1}\dff \colon\dff
\widehat{H}_{\dff 1}\qff \ttoo\qff H_{\dff 1}$\sss
be\sss the projections onto\sss the first\sss and\sss the second summands respectively.\oss
Let\sss\vspace{3pt} 
\[
\quad
\pi_{\dff 0}
\off =\off 
\gamma\dff \circ\dff \Pi_{\dff 0}
\qff,\quad
\pi_{\dff 1}
\off =\off 
\gamma\dff \circ\dff \Pi_{\dff 1}
\qff,
\quad
\mbox{and}\quad
\pi
\off =\off 
\pi_{\dff 0}\dff \oplus\dff \pi_{\dff 1}
\pff.
\]

\vspace{-12pt}\vspace{3pt}
Let\sss 
$\Upsilon\dff \colon\dff 
K^{\dff \partial}\dff \oplus\dff K^{\dff \partial}
\qff \ttoo\qff 
K^{\dff \partial}\dff \oplus\dff K^{\dff \partial}$\sss
be a unitary operator\sss leaving\sss $K\dff \oplus\dff K$\sss invariant,\oss
and\dss let\sss us impose on $S^*$\sss the boundary condition\sss
$\pi_{\dff 1}\off =\off \Upsilon\fff \circ\trf \pi_{\dff 0}$\nsp,\oss
i.e.\qss consider\sss the restriction $P$ of\dss $S^*$\sss to\sss 
$\kernel\trf(\trf \pi_{\dff 1}\qff -\qff \Upsilon\fff \circ\trf \pi_{\dff 0}\trf)$\nsp.\oss
Equivalently,\pss 
$P$\sss is\dss the restriction of\dss $S^*$\sss to\sss
$\pi^{\dff -\dff 1}\dff(\dff \mathcal{C}\trf)$\nnsp,\oss
where\sss
$\mathcal{C}$\sss is\dss the space\sss of\dss pairs of\dss the form\sss
$(\trf u\fff,\qff \Upsilon\trf(\trf u\trf)\trf)$\nnsp.\oss
We will\dss further assume\sss that\sss $\Upsilon$\sss commutes with $\Sigma$\nnsp.\oss
The next\dss lemma shows\sss that\sss then\sss this boundary condition\dss
is\dss self-adjoint\sss in\sss the sense of\qss \cite{i2},\oss Section\qss 5.\oss
\emph{We will\sss assume\sss that,\oss 
moreover{},\oss the operator $P$\sss is\dss 
self-adjoint\sss and\dss Fredholm.\oss}

\mypar{Lemma.}{self-adjointness}
\emph{The above boundary condition\dss is\dss self-adjoint,\oss
i.e.\qss $\widehat{\Sigma}\qff(\trf \mathcal{C}\trf)$\sss
is\dss equal\dss to\sss the orthogonal\sss complement\sss of\dss
$\mathcal{C}$\dnsp.\oss}

\proof
Let\sss
$(\trf u\fff,\qff \Upsilon\trf(\trf u\trf)\trf)
\dff,\off
(\trf v\fff,\qff \Upsilon\trf(\trf v\trf)\trf)
\qff \in\qff \mathcal{C}$\dnsp.\oss
Then\vspace{3pt}
\[
\quad
\bsco{\dff \widehat{\Sigma}\qff(\trf u\fff,\qff \Upsilon\trf(\trf u\trf)\trf)\fff,\off
(\trf v\fff,\qff \Upsilon\trf(\trf v\trf)\trf)\dff}_{\dff \partial}
\]

\vspace{-33pt}
\[
\quad
=\off
\bsco{\dff (\trf \Sigma\trf u\fff,\qff -\qff \Sigma\dff \circ\dff \Upsilon\trf(\trf u\trf)\trf)\fff,\off
(\trf v\fff,\qff \Upsilon\trf(\trf v\trf)\trf)\dff}_{\dff \partial}
%\off =\off
%\bsco{\dff (\trf \Sigma\trf u\fff,\qff -\qff \Upsilon\dff \circ\dff \Sigma\trf(\trf u\trf)\trf)\fff,\off
%(\trf v\fff,\qff \Upsilon\trf(\trf v\trf)\trf)\dff}_{\dff \partial}
\]

\vspace{-33pt}
\[
\quad
=\off
\bsco{\dff \Sigma\trf u\fff,\qff v\dff}_{\dff \partial}
\qff -\qff
\bsco{\dff \Upsilon\dff(\trf \Sigma\trf(\trf u\trf)\trf)\fff,\qff
\Upsilon\dff(\trf v\trf)\dff}_{\dff \partial}
\off =\off
\bsco{\dff \Sigma\trf u\fff,\qff v\dff}_{\dff \partial}
\qff -\qff
\bsco{\dff \Sigma\trf u\fff,\qff v\dff}_{\dff \partial}
\off =\off
0
\]

\vspace{-12pt}\vspace{3pt}
because $\Upsilon$\sss commutes with $\Sigma$ and\dss is\dss unitary.\oss
It\sss follows\sss that\sss $\widehat{\Sigma}\qff(\trf \mathcal{C}\trf)$\sss
is\dss orthogonal\dss to $\mathcal{C}$\dnsp.\oss
Conversely,\oss suppose\sss that\sss
$(\trf a\fff,\qff b\trf)$\sss is\dss orthogonal\dss to\sss $\mathcal{C}$\dnsp.\oss
Then\sss for every\sss $u\qff \in\qff K\dff \oplus\dff K$\vspace{3pt}
\[
\quad
0
\off =\off
\bsco{\dff (\trf a\fff,\qff b\trf)\fff,\qff
(\trf u\fff,\off \Upsilon\trf(\trf u\trf)\trf)\dff}
\off =\off
\bsco{\dff a\fff,\qff u\dff}
\qff +\qff
\bsco{\dff b\fff,\qff \Upsilon\trf(\trf u\trf)\trf)\dff}
\pff.
\]

\vspace{-12pt}\vspace{3pt}
Since\sss $\Upsilon$\sss is\dss unitary and\sss
$K$\sss is\dss dense in\sss $K^{\dff \partial}$\dnsp,\oss
this implies\sss that\sss
$b\off =\off -\qff \Upsilon\dff(\trf a\trf)$\nnsp.\oss
Since\sss $\Sigma$\sss is\dss invertible,\oss
this implies\sss that\sss
$(\trf a\fff,\qff b\trf)
\off =\off
(\trf a\fff,\qff -\qff \Upsilon\dff(\trf a\trf)\trf)
\qff \in\pff
\widehat{\Sigma}\qff(\trf \mathcal{C}\trf)$\nnsp.\oss
The\sss lemma follows.\oss  \eproof

\mypar{Theorem.}{reference-operator}
\emph{If\qss the operator\sss 
$-\qff i\trf \Sigma^{\dff *}\dff \circ\dff \Upsilon$\sss
is\trs positive\sss definite,\oss
then\dss the operator\trs $P\qff +\qff \varepsilon$\sss
is\dss an\sss isomorphism\trs
$\mathcal{D}\trf(\trf P\trf)\qff \ttoo\qff \widehat{H}_{\dff 0}$\nsp.\oss}

\proof
Since $P$\sss is\dss assumed\sss to be self-adjoint\sss and $\varepsilon$\sss
is\dss bounded,\pss $P\qff +\qff \varepsilon$\sss is\dss
a self-adjoint\sss operator\sss in $\widehat{H}_{\dff 0}$\nsp.\oss
Therefore it\dss is\dss sufficient\sss to prove\sss that\sss the kernel\sss
of\dss  $P\qff +\qff \varepsilon$\sss is\dss equal\dss to $0$\nnsp.\oss
Suppose\sss that\sss
$(\trf a\fff,\qff b\trf)\qff \in\qff \kernel\fff (\trf P\qff +\qff \varepsilon\trf)$\nnsp.\oss
Then\sss $T^{\dff *}\fff a\qff +\qff b\off =\off 0$
and\sss $a\qff -\qff T^{\dff *}\fff b\off =\off 0$\nnsp.\oss
Let\sss us apply\sss the\dss Lagrange\dss identity\sss for $S^*$\sss to\sss
$u\off =\off (\trf a\fff,\qff a\trf)$ and\sss
$v\off =\off (\trf b\fff,\qff -\qff b\trf)$\nnsp.\oss
Since\sss $\varepsilon$\sss is\dss self-adjoint,\pss
$\varepsilon$\sss does not\sss affect\sss the\sss left\sss hand\sss side of\dss the\dss
Lagrange\dss identity and\dss hence\sss the\sss latter\dss
is\dss equal\dss to\vspace{4.5pt}
\[
\quad
\sco{\dff T^{\dff *} a\fff,\qff b \dff}
\pff -\pff
\sco{\dff a\dff,\qff T^{\dff *} b \dff}
\pff +\pff
\sco{\dff -\qff T^{\dff *} a\fff,\qff -\qff b \dff}
\pff -\pff
\sco{\dff a\dff,\qff -\qff T^{\dff *}\dff(\dff -\qff b\trf) \dff}
\]

\vspace{-33pt}
\[
\quad
=\off
2\trf \sco{\dff T^{\dff *} a\fff,\qff b \dff}
\pff -\pff
2\trf \sco{\dff a\dff,\qff T^{\dff *} b \dff}
\off =\off
-\qff 2\trf \sco{\dff b\fff,\qff b \dff}
\pff -\pff
2\trf \sco{\dff a\fff,\qff a \dff}
\pff,
\]

\vspace{-12pt}\vspace{4.5pt}
where at\sss the\sss last\sss step we used\sss the fact\sss that\sss
$T^{\dff *}\fff a\off =\off -\qff b$ and\sss
$T^{\dff *}\fff b\off =\off a$\nnsp.\oss
In\sss particular{},\oss the\sss left\sss hand side\dss is\dss $\leq\qff 0$\nnsp.\oss
The right\sss hand side of\dss the\dss Lagrange\dss identity\dss
is\dss equal\dss to\vspace{4.5pt}\vspace{-0.5pt}
\[
\quad
\bsco{\dff i\qff \Sigma\trf \gamma\fff a\dff,\qff \gamma\dff b\dff}_{\dff \partial}
\pff +\pff
\bsco{\dff i\trf(\dff -\qff \Sigma\trf)\trf \gamma\fff a\dff,\qff -\qff \gamma\dff b\dff}_{\dff \partial}
\off =\off
2\trf \bsco{\dff i\qff \Sigma\trf \gamma\fff a\dff,\qff \gamma\dff b\dff}_{\dff \partial}
\pff.
\]

\vspace{-12pt}\vspace{4.5pt}\vspace{-0.5pt}
At\sss the same\sss time\sss
$(\trf a\fff,\qff b\trf)\qff \in\qff \mathcal{D}\trf(\trf P\trf)$\sss
and\dss hence\sss
$\gamma\fff b\off =\off \Upsilon\dff (\dff \gamma\fff a\trf)$\nnsp.\oss
Therefore\vspace{4.5pt}\vspace{-0.5pt}
\[
\quad
\bsco{\dff i\qff \Sigma\trf \gamma\fff a\dff,\qff 
\gamma\dff b\dff}_{\dff \partial}
\off =\off
\bsco{\dff i\qff \Sigma\trf \gamma\fff a\dff,\qff 
\Upsilon\dff \gamma\fff a\dff}_{\dff \partial}
\off =\off
\bsco{\dff \gamma\fff a\dff,\qff 
-\qff i\trf \Sigma^{\dff *}\dff \circ\dff \Upsilon\dff (\dff \gamma\fff a\trf)\dff}_{\dff \partial}
\pff.
\]

\vspace{-12pt}\vspace{4.5pt}\vspace{-0.5pt}
It\sss follows\sss that\sss the right\sss hand side\dss is\dss $\geq\qff 0$\dss
if\dss $-\qff i\trf \Sigma^{\dff *}\dff \circ\dff \Upsilon$\sss
is\trs positive\sss definite.\oss
Since\sss the\sss left\sss hand side\dss is\dss $\leq\qff 0$\nnsp,\oss
in\sss this case both sides are equal\dss to $0$\nnsp.\oss
This implies\sss that\sss 
$\sco{\dff b\fff,\qff b \dff}
\qff +\qff
\sco{\dff a\fff,\qff a \dff}
\off =\off
0$
and\dss hence $a\off =\off b\off =\off 0$\nnsp.\oss
It\sss follows\sss that\sss
$\kernel\fff (\trf P\qff +\qff \varepsilon\trf)\off =\off 0$\nnsp.\oss  \eproof

\myuppar{Changing\sss the boundary operators.}
Suppose\sss that\sss $\Sigma$\sss is\dss unitary and\dss 
let\sss us\sss take $\Upsilon\off =\off i\trf \Sigma$\nnsp.\oss
Then\sss
$-\qff i\trf \Sigma^{\dff *}\dff \circ\dff \Upsilon
\off =\off
\Sigma^{\dff *}\fff \circ\trf \Sigma$\sss
is\trs positive\sss definite and\sss $\Upsilon$\sss is\dss skew-adjoint.\oss
Let\vspace{3pt}
\[
\quad
\Gamma_0
\off =\off
(\trf \pi_{\dff 0}\qff -\qff \Upsilon^{\dff -\dff 1}\fff \circ\trf \pi_{\dff 1}\trf)\bigl/\sqrt{2}
\quad
\mbox{and}\quad
\Gamma_1
\off =\off
(\trf \Upsilon\fff \circ\trf \pi_{\dff 0}\qff +\qff \pi_{\dff 1}\trf)\bigl/\sqrt{2}
\off.
\]

\vspace{-12pt}\vspace{3pt}
Using\sss the facts\sss that\sss $\Upsilon$\sss is\dss skew-adjoint\sss and\sss unitary,\oss
we see\sss that\vspace{3pt}
\[
\hspace*{-0.6em}
2\dff
\bsco{\dff \Gamma_{\dff 1}\dff u\dff,\qff
\Gamma_{\dff 0}\dff v\dff}_{\dff \partial}
\pff -\pff
2\dff
\bsco{\dff \Gamma_{\dff 0}\dff u\dff,\qff
\Gamma_{\dff 1}\dff v\dff}_{\dff \partial}
\]

\vspace{-33pt}
\[
\hspace*{-0.6em}
=\off
\bsco{\dff \Upsilon\fff \circ\trf \pi_{\dff 0}\dff u\qff +\qff \pi_{\dff 1}\dff u\dff,\qff
\pi_{\dff 0}\dff v\qff -\qff \Upsilon^{\dff -\dff 1}\fff \circ\trf \pi_{\dff 1}\dff v\dff}_{\dff \partial}
\pff -\pff
\bsco{\dff \pi_{\dff 0}\dff u\qff -\qff \Upsilon^{\dff -\dff 1}\fff \circ\trf \pi_{\dff 1}\dff u\dff,\qff
\Upsilon\fff \circ\trf \pi_{\dff 0}\dff v\qff +\qff \pi_{\dff 1}\dff v\dff}_{\dff \partial}
\]

\vspace{-33pt}
\[
=\off
\bsco{\dff \Upsilon\vphantom{^{-\dff 1}}\fff \circ\trf \pi_{\dff 0}\dff u\dff,\qff 
\pi_{\dff 0}\dff v\dff}_{\dff \partial}
\qff -\qff
\bsco{\dff \Upsilon\fff \circ\trf \pi_{\dff 0}\dff u\dff,\qff 
\Upsilon^{\dff -\dff 1}\fff \circ\trf \pi_{\dff 1}\dff v\dff}_{\dff \partial}
\qff +\qff
\bsco{\dff \vphantom{^{-\dff 1}}\pi_{\dff 1}\dff u\dff,\qff \pi_{\dff 0}\dff v\dff}_{\dff \partial}
\qff -\qff
\bsco{\dff \pi_{\dff 1}\dff u\dff,\qff \Upsilon^{\dff -\dff 1}\fff \circ\trf \pi_{\dff 1}\dff v\dff}_{\dff \partial}
\]

\vspace{-33pt}
\[
\hspace*{-0.6em}
-\qff
\bsco{\dff \pi_{\dff 0}\dff u\dff,\qff \Upsilon\vphantom{^{-\dff 1}}\fff \circ\trf \pi_{\dff 0}\dff v\dff}_{\dff \partial}
\qff -\qff
\bsco{\dff \vphantom{^{-\dff 1}}\pi_{\dff 0}\dff u\dff,\qff \pi_{\dff 1}\dff v\dff}_{\dff \partial}
\qff +\qff
\bsco{\dff \Upsilon^{\dff -\dff 1}\fff \circ\trf \pi_{\dff 1}\dff u\dff,\qff 
\Upsilon\fff \circ\trf \pi_{\dff 0}\dff v \dff}_{\dff \partial}
\qff +\qff
\bsco{\dff \Upsilon^{\dff -\dff 1}\fff \circ\trf \pi_{\dff 1}\dff u\dff,\qff 
\pi_{\dff 1}\dff v \dff}_{\dff \partial}
\]

\vspace{-33pt}
\[
\hspace*{-0.6em}
=\off
2\dff
\bsco{\dff \Upsilon\vphantom{^{-\dff 1}}\fff \circ\trf \pi_{\dff 0}\dff u\dff,\qff \pi_{\dff 0}\dff v\dff}_{\dff \partial}
\qff +\qff
2\dff
\bsco{\dff \Upsilon^{\dff -\dff 1}\fff \circ\trf \pi_{\dff 1}\dff u\dff,\qff \pi_{\dff 1}\dff v\dff}_{\dff \partial}
\]

\vspace{-33pt}
\[
\hspace*{-0.6em}
=\off
2\dff
\bsco{\dff \Upsilon\vphantom{^{-\dff 1}}\fff \circ\trf \pi_{\dff 0}\dff u\dff,\qff \pi_{\dff 0}\dff v\dff}_{\dff \partial}
\qff -\qff
2\dff
\bsco{\dff \Upsilon\vphantom{^{-\dff 1}}\fff \circ\trf \pi_{\dff 1}\dff u\dff,\qff \pi_{\dff 1}\dff v\dff}_{\dff \partial}
\off\qff =\off\qff
2\dff
\bsco{\dff i\qff \widehat{\Sigma}\pff \widehat{\gamma}\dff u\dff,\qff \widehat{\gamma}\dff v\dff}_{\dff \partial}
\off.
\]

\vspace{-12pt}\vspace{3pt}
This shows\sss that\sss the\dss Lagrange\dss identity\sss in\sss terms of\dss
$\Gamma_0\dff,\pff \Gamma_1$\sss
has\sss the standard\sss form.\oss
The operator\sss $A\off =\off P\qff +\qff \varepsilon$\sss
is\dss defined\dss by\sss the boundary condition\sss $\Upsilon\dff \circ\dff \Gamma_0\off =\off 0$\nnsp,\oss
which\dss is\dss equivalent\sss to\sss $\Gamma_0\off =\off 0$\nnsp.\oss
Theorem\qss \ref{reference-operator}\qss implies\sss that\sss 
we can\sss take\sss
$A\off =\off P\qff +\qff \varepsilon$\sss
as\sss the reference operator{}.\oss

\myuppar{The difference of\dss two boundary conditions.}
Now we will\sss return\sss to\sss the situation of\trs Section\qss \ref{abstract-index}\qss
and assume\sss that\sss the\dss Lagrange\dss identity\sss for $T^{\dff *}$\sss has\sss
the form\qss (\ref{lagrange-identity}).\oss
Equivalently,\oss\vspace{1.5pt} 
\[
\quad
i\trf \Sigma
\off =\off\dff
\begin{pmatrix}
\off 0 &
1 \qff\off
\vspace{4.5pt} \\
\off\dff -\qff 1 &
0 \qff\off 
\end{pmatrix}
\off
\]

\vspace{-12pt}\vspace{1.5pt}
with respect\sss to\sss the direct\sss sum\sss
$K^{\dff \partial}\dff \oplus\dff K^{\dff \partial}$\dnsp.\oss
Let\sss $\mathcal{B}_{\trf 0}\dff,\off \mathcal{B}_{\dff 1}
\off \subset\off 
K^{\dff \partial}\dff \oplus\dff K^{\dff \partial}$\sss 
be\sss two closed\sss relations.\oss 
Let $T_{\dff 0}\dff,\qff T_{\fff 1}$\sss be\sss the restrictions of\dss $T^{\dff *}${\nsp}  to\sss
$\pi^{\dff -\dff 1}\dff(\trf \mathcal{B}_{\trf 0}\trf)\fff,\off
\pi^{\dff -\dff 1}\dff(\trf \mathcal{B}_{\dff 1}\trf)
\qff \subset\qff 
H_{\dff 1}$\sss respectively.\oss
\emph{Let\sss us assume\sss that\dss 
$T_{\dff 0}\dff,\qff T_{\fff 1}$ 
are self-adjoint\sss operators in\sss $H_{\trf 0}$\nsp.\oss}
But\sss we will\sss not\sss assume\sss that\sss
$T_{\dff 0}$ or\sss $T_{\fff 1}$\sss is\dss invertible.\oss
The direct\sss sum\sss
$T_{\dff 0}\dff \oplus\dff -\qff T_{\fff 1}$\sss
is\dss the\qss \emph{formal\sss difference}\qss of\dss operators $T_{\dff 0}$ and\sss $T_{\fff 1}$\nsp.\oss
%When parameters are present,\oss the index of\dss
%$T_{\dff 0}\dff \oplus\dff -\qff T_{\fff 1}$\sss
%is\dss equal\dss to\sss the difference of\dss indices of\dss $T_{\dff 0}$ and $T_{\fff 1}$\nsp,\oss
%and we will\dss assume\sss that\sss everything continuously depends on parameters,\oss
%but\sss will\sss omit\sss the parameters from\sss notations.\oss

The formal\sss difference\sss
$T_{\dff 0}\dff \oplus\dff -\qff T_{\fff 1}$\sss
is\dss
the self-adjoint\sss extension of\dss
$S\off =\off T\dff \oplus\dff -\qff T$
defined\dss by\sss the boundary condition\sss
$\mathcal{B}_{\trf 0}\dff \oplus\dff \mathcal{B}_{\dff 1}$\nsp.\oss
This extension\sss
can be defined also in\sss terms of\dss boundary operators\sss
$\Gamma_0\dff,\pff \Gamma_1$\nsp.\oss
Namely,\oss let\sss
$\Gamma\off =\off \Gamma_0\dff \oplus\dff \Gamma_1$
and\dss let\sss $\Phi$\sss be\sss the automorphism of\vspace{1.5pt}\vspace{1.5pt}
\[
\quad
\bigl(\qff
K^{\dff \partial}\dff \oplus\dff K^{\dff \partial}
\qff\bigr)
\qff \oplus\qff
\bigl(\qff
K^{\dff \partial}\dff \oplus\dff K^{\dff \partial}
\qff\bigr)
\]

\vspace{-12pt}\vspace{1.5pt}\vspace{1.5pt}
defined\sss by\sss the formula\sss %\vspace{1.5pt}
%\[
%\quad
$\Phi\trf(\trf a\fff,\qff b\trf)
\off =\off
\bigl(\trf 
a\qff -\qff \Upsilon^{\dff -\dff 1}\dff(\trf b\trf)\dff,\pff
\Upsilon\dff(\trf a\trf)\qff +\qff b\trf\bigr)$\nnsp.\oss
%\pff.
%\]
%
%\vspace{-12pt}\vspace{1.5pt}
Then\sss
$\Phi\trf(\trf \mathcal{B}_{\trf 0}\dff \oplus\dff \mathcal{B}_{\dff 1}\trf)$\sss
is\dss a closed\sss self-adjoint\sss relation and\sss
$T_{\dff 0}\dff \oplus\dff -\qff T_{\fff 1}$\sss
is\dss the restriction of\dss $S^*$\sss to\vspace{3pt}
\[
\quad
\Gamma^{\dff -\dff 1}\dff\bigl(\qff
\Phi\trf(\trf \mathcal{B}_{\trf 0}\dff \oplus\dff \mathcal{B}_{\dff 1}\trf)
\qff\bigr)
\pff.
\]

\vspace{-12pt}\vspace{3pt}
In other words,\pss $T_{\dff 0}\dff \oplus\dff -\qff T_{\fff 1}$\sss
is\dss equal\dss to\sss the extension of\dss $S$ defined\dss in\sss terms of\dss
$\Gamma_0\dff,\pff \Gamma_1$ by\sss the boundary condition\sss
$\Phi\trf(\trf \mathcal{B}_{\trf 0}\dff \oplus\dff \mathcal{B}_{\dff 1}\trf)$\nnsp.\oss
When parameters are present,\oss the index of\dss
$T_{\dff 0}\dff \oplus\dff -\qff T_{\fff 1}$\sss
is\dss equal\dss to\sss the difference of\dss indices of\dss $T_{\dff 0}$ and $T_{\fff 1}$\nsp,\oss
and we will\dss assume\sss that\sss everything continuously depends on parameters,\oss
but\sss will\sss omit\sss the parameters from\sss notations.\oss

For\sss the purposes of\dss determining\sss the index,\oss
we can replace $S^*$\sss by\sss $S^*\qff +\qff \varepsilon$\nnsp,\oss
or{},\oss equivalently,\oss to consider\qss ({\fff}the family of\dss operators)\dss
$(\trf T_{\dff 0}\dff \oplus\dff -\qff T_{\fff 1}\trf)\qff +\qff \varepsilon$\nnsp.\oss
Then we can use\sss $A\off =\off P\qff +\qff \varepsilon$\sss
as\sss the reference operators and construct\sss the reduced\dss boundary\sss triplets.\oss
We can also define\sss the\qss (families of\trf)\qss operators\sss
$M\trf(\dff 0\dff)$ and\sss $M$\nnsp.\oss
In order\sss to be able\sss to apply\sss the results of\trs Section\qss \ref{families}\qss
we need\dss to assume\sss that\sss either\sss the operators\sss
$P\qff +\qff \varepsilon$ have compact\sss resolvent,\oss
or\sss that\sss the relations\sss
$\mathcal{B}_{\trf 0}\dff,\off \mathcal{B}_{\dff 1}$\sss
are self-adjoint\trs Fredholm\dss relations with compact\sss resolvent.\oss
We need also\sss the continuity assumptions from\sss Section\qss \ref{families}.\oss
Then\sss the index of\dss the family\sss
$(\trf T_{\dff 0}\dff \oplus\dff -\qff T_{\fff 1}\trf)\qff +\qff \varepsilon$\nnsp,\oss
and\dss hence\sss the index of\dss the family\sss
$T_{\dff 0}\dff \oplus\dff -\qff T_{\fff 1}$\nsp,\oss
is\dss equal\dss to\sss the index of\dss the family\vspace{4.5pt}
\[
\quad
\Lambda'\dff \oplus\dff \Lambda^{\fff -\dff 1}\qff
\bigl(\qff
\Phi\trf(\trf \mathcal{B}_{\trf 0}\dff \oplus\dff \mathcal{B}_{\dff 1}\trf)
\trf\bigl|\trf 
(\trf K\dff \oplus\dff K\trf)
\qff -\qff M
\qff\bigr)
\pff.
\]

\vspace{-12pt}\vspace{4.5pt}
In\sss this computation of\dss the index 
one can return\sss to\sss the original\sss direct\sss sum\sss
$\mathcal{B}_{\trf 0}\dff \oplus\dff \mathcal{B}_{\dff 1}$\nsp.\oss
In order\sss to do\sss this,\oss let\sss us consider\sss $M$\sss
as a closed\sss relation\sss in\sss
$K\dff \oplus\dff K$\sss
and\dss let\sss
$M_{\dff \oplus}
\off =\off
\Phi^{\dff -\dff 1}\dff(\trf M\trf)$\nnsp.\oss
Then\sss the index of\dss the family\sss
$T_{\dff 0}\dff \oplus\dff -\qff T_{\fff 1}$\sss
is\dss equal\dss to\sss the index of\dss the family\vspace{4.5pt}
\begin{equation}
\label{difference-1}
\quad
\Lambda'\dff \oplus\dff \Lambda^{\fff -\dff 1}\qff
\bigl(\qff
(\trf \mathcal{B}_{\trf 0}\dff \oplus\dff \mathcal{B}_{\dff 1}\trf)
\trf\bigl|\trf 
(\trf K\dff \oplus\dff K\trf)
\qff -\qff M_{\dff \oplus}
\qff\bigr)
\pff,
\end{equation}

\vspace{-12pt}\vspace{4.5pt}
and\dss hence\sss the difference of\dss the indices of\dss the families\sss
$T_{\dff 0}$ and\sss $T_{\fff 1}$\sss
is\dss also equal\dss to\sss this index.\oss

At\sss the end of\trs Section\qss \ref{gelfand-triples}\qss we defined
self-adjoint\sss relations
from\sss $K$\sss to\sss $K\fff'$\dnsp.\oss
One can also define\sss the index of\dss families of\dss such self-adjoint\sss relations.\oss
We\sss leave\sss this\sss to\sss the reader{}.\oss
By\trs Lemma\qss \ref{lambda-lambda}\qss the operators $\Lambda$ and\sss $\Lambda'$
are adjoint\sss to each other{}.\oss
It\sss follows\sss that\sss the\sss latter\sss index,\oss and\dss hence\sss
the difference of\dss the indices,\oss is\dss equal\dss to\sss the index of\dss
the family\vspace{4.5pt}
\begin{equation}
\label{difference-2}
\quad
(\trf \mathcal{B}_{\trf 0}\dff \oplus\dff \mathcal{B}_{\dff 1}\trf)
\trf\bigl|\trf 
(\trf K\dff \oplus\dff K\trf)
\qff -\qff M_{\dff \oplus}
\pff
\end{equation}

\vspace{-12pt}\vspace{4.5pt}
of\dss relations from\sss $K\dff \oplus\dff K$\sss
to\sss $K\fff'\dff \oplus\dff K\fff'$\dnsp.\oss
Here\sss the relation\sss $M_{\dff \oplus}$ mixes\sss the summands preserved\dss
by\sss $\mathcal{B}_{\trf 0}\dff \oplus\dff \mathcal{B}_{\dff 1}$\nsp,\oss
and\sss the fact\sss that\sss we are dealing\sss with\sss the difference\dss
is\dss reflected only\sss in\sss $M_{\dff \oplus}$\nsp.\oss
Of\dss course,\oss there\dss is\dss no similar result\sss for\sss the sum.\oss
The above arguments do not\sss work\sss for\sss the sum\sss because\sss
there\dss is\dss no analogue of\trs Theorem\qss \ref{reference-operator}\qss
if\dss $T\dff \oplus\dff -\qff T$\sss is\dss replaced\dss by\sss $T\dff \oplus\dff T$\dnsp.\oss

\myuppar{The differential\dss boundary\sss problems of\dss order one.}
Suppose now\sss that\sss we are in\sss the situation of\trs 
Section\qss \ref{differential-boundary-problems}\qss
Let\sss $\rho_{\dff u}\off =\off \sigma_y^{\dff -\dff 1}\qff \tau_{\fff u}$\sss be related\dss to\sss $T$\dnsp,\oss
and\dss let\sss $\rho'_{\fff u}$ and 
$\widehat{\rho}_{\dff u}\off =\off \rho_{\dff u}\dff \oplus\dff \rho'_{\fff u}$\sss 
be similar operators related\dss to\sss $-\qff T$ and\sss
$S\off =\off T\dff \oplus \dff -\qff T$\sss respectively.\oss
Then\sss $\rho'_{\fff u}
\off =\off 
(\dff -\qff \sigma_y^{\dff -\dff 1}\dff)\qff (\dff -\qff \tau_{\fff u}\dff)
\off =\off
\rho_{\dff u}$
and\dss hence\sss
$\widehat{\rho}_{\dff u}\off =\off \rho_{\dff u}\dff \oplus\dff \rho_{\dff u}$\nsp.\oss
It\sss follows\sss that\sss\vspace{0pt}
\[
\quad
\mathcal{L}_{\dff -}\dff(\trf \widehat{\rho}_{\dff u}\trf)
\off =\off
\mathcal{L}_{\dff -}\dff(\trf \rho_{\dff u}\trf)
\dff \oplus\dff
\mathcal{L}_{\dff -}\dff(\trf \rho_{\dff u}\trf)
\qff.
\]

\vspace{-12pt}\vspace{0pt}
Recall\dss that\sss the subspace 
$\mathcal{L}_{\dff -}\dff(\trf \rho_{\dff u}\trf)$\sss
is\dss lagrangian 
and\dss hence\sss
$i\trf \Sigma\trf(\trf \mathcal{L}_{\dff -}\dff(\trf \rho_{\dff u}\trf)\trf)
\qff \cap\qff
\mathcal{L}_{\dff -}\dff(\trf \rho_{\dff u}\trf)
\off =\off
0$\nnsp.\oss
It\sss follows\sss that\sss the kernel\sss of\dss the map\sss
$(\trf u\fff,\qff v\trf)
\off \longmapsto\off
i\trf \Sigma\trf(\trf u\trf)\qff -\qff v$\sss
intersects\sss
$\mathcal{L}_{\dff -}\dff(\trf \widehat{\rho}_{\dff u}\trf)$
only\sss by $0$\nnsp.\oss
Similarly,\oss the kernel\sss of\dss the map\sss
$(\trf u\fff,\qff v\trf)
\off \longmapsto\off
i\trf \Sigma\trf(\trf u\trf)\qff +\qff v$\sss
intersects\sss
$\mathcal{L}_{\dff -}\dff(\trf \widehat{\rho}_{\dff u}\trf)$
only\sss by $0$\nnsp.\oss
This implies\sss that\sss the boundary conditions\sss $\Gamma_0\off =\off 0$\sss
and\sss $\Gamma_1\off =\off 0$\sss
satisfy\sss the\dss Shapiro--Lopatinskii\dss condition.\oss
Since\sss the operators $P$\sss is\dss defined\dss by\sss the boundary condition\sss
$\Gamma_0\off =\off 0$\nnsp,\oss
it\sss follows\sss that\sss the operators $P$\sss are self-adjoint\sss and\dss Fredholm.\oss
Similarly,\oss the operators\sss $P'$\sss defined\dss by\sss the boundary conditions 
$\Gamma_1\off =\off 0$\sss are also self-adjoint\sss and\dss Fredholm.\oss

Suppose now\sss that\sss everything depends on a parameter $w\qff \in\qff W$
as at\sss the end of\trs Section\qss \ref{differential-boundary-problems}.\oss
The continuity assumptions of\trs Section\qss \ref{families}\qss hold\dss
by\sss the same reasons as in\dss Section\qss \ref{differential-boundary-problems}.\oss
Also,\oss the operators $A$ are operators with compact\sss resolvent.\oss
Hence,\oss as in\sss the proof\dss of\trs Theorem\qss \ref{index-theorem},\oss 
the results of\trs Section\qss \ref{families}\qss apply\sss
and\dss the difference of\dss the indices of\dss the families\sss
$T_{\dff 0}$ and\sss $T_{\fff 1}$\sss
is\dss equal\dss to\sss the index of\dss the family\qss (\ref{difference-1}),\oss
as also\sss to\sss the index of\dss the family\qss (\ref{difference-2}).\oss

In\sss the present\sss context\sss one can describe $M$ and\sss $M_{\dff \oplus}$\sss
in\sss terms of\dss the Cauchy\dss data of\dss
the equation\sss $(\trf S^*\qff +\qff \varepsilon\trf)\trf u\off =\off 0$\nnsp.\oss
Namely,\pss by\sss the results of\trs Section\qss \ref{differential-boundary-problems}\pss
the operator\sss $M$\sss is\dss the operator\sss
having\sss as its graph\sss the image of\trs 
$\widehat{H}_{\dff 1}\qff \cap\qff \kernel\dff (\trf S^*\qff +\qff \varepsilon \trf)$\sss
under\sss the map\sss $\Gamma\off =\off \Gamma_0\dff \oplus\dff \Gamma_1$\nsp.\oss
Clearly,\pss
$\Gamma\off =\off \Phi\dff \circ\dff \pi$\nnsp.\oss
It\dss follows\sss that\sss
$M_{\dff \oplus}
\off =\off
\Phi^{\dff -\dff 1}\dff(\trf M\trf)$\sss
is\dss the relation\sss having\sss as\sss its graph\qss ({\fff}i.e.\qss equal\dss to)\qss
the image of\dss
$\widehat{H}_{\dff 1}\qff \cap\qff \kernel\dff (\trf S^*\qff +\qff \varepsilon \trf)$\sss
under\sss the map\sss $\pi$\nnsp.\oss

\mysection{Rellich\qss example}{rellich}

\myuppar{Rellich\sss example.}
Let\sss $H_{\dff 0}\off =\off L_{\dff 2}\dff [\trf 0\fff,\qff 1\trf]$\nnsp.\oss
Let\sss $T$\sss be\sss the differential\sss operator\sss 
$-\qff d^{\dff 2}/\dff d x^{\dff 2}$ with\sss the domain\sss
$H_{\trf 2}^{\trf 0}\qff [\trf 0\fff,\qff 1\trf]$\nnsp,\oss
the subspace of\dss the\dss Sobolev\dss space\sss
$H_{\trf 2}\qff [\trf 0\fff,\qff 1\trf]$\sss
defined\dss by\sss the boundary conditions\sss
$u\trf(\dff 0\dff)\off =\off
u\trf(\dff 1\dff)\off =\off
u'\trf(\dff 0\dff)\off =\off
u'\trf(\dff 1\dff)\off =\off
0$\nnsp.\oss
For\sss $\kappa\qff \in\qff \rrr\trf \cup\dff \{\trf \infty \trf\}$\sss
let\sss $T\dff(\trf \kappa\trf)$\sss be\sss the restriction of\dss $T^{\dff *}$
defined\sss by\sss the boundary conditions\sss
$u\trf(\dff 0\dff)\off =\off 0$\nnsp,\qss
$\kappa\dff u'\trf(\dff 1\dff)\off =\off u\trf(\dff 1\dff)$\nnsp.\oss
For\sss $\kappa\off =\off \infty$\sss
the condition\sss
$\kappa\dff u'\trf(\dff 1\dff)\off =\off u\trf(\dff 1\dff)$\sss
is\dss interpreted as\sss 
$u'\trf(\dff 1\dff)\off =\off 0$\nnsp.\oss
The family of\dss operators\sss $T\dff(\trf \kappa\trf)$\sss is\dss
an\sss important\sss example in\sss the\sss theory of\dss self-adjoint\sss operators.\oss
See\qss Kato\qss \cite{k},\oss Example\qss V.4.14.\oss
This family\dss is\dss continuous in\sss the\sss topology of\dss
the uniform\sss resolvent\sss convergence,\oss even at\sss $\kappa\off =\off \infty$\sss
if\dss we consider\sss $\rrr\trf \cup\dff \{\trf \infty \trf\}$\sss
as\sss the one-point\sss compactification of\dss $\rrr$\nnsp.\oss
This one-point\sss compactification can\sss be identified\sss with\sss the circle $S^{\fff 1}$\dnsp,\oss
say,\oss by\sss the stereographic\sss projection.\oss
Hence\sss the index of\dss this family\sss belongs\sss to\sss
$K^{\dff 1}\dff(\trf S^{\fff 1}\trf)\off =\off \zzz$\nnsp.\oss
Alternatively,\oss this index\dss is\dss equal\dss to\sss the spectral\dss flow
of\dss the family\sss $T\dff(\trf \kappa\trf)$\nnsp,\dss
$\kappa\qff \in\qff \rrr\trf \cup\dff \{\trf \infty \trf\}$\nnsp.\oss
The behavior of\dss eigenvalues in\sss this family\dss is\dss known.\oss
See\sss the picture in\qss \cite{k},\oss loc.\dss cit.\oss
It\dss is\dss clear\sss from\sss this picture\sss that\sss the spectral\dss flow,\oss
and\dss hence\sss the index,\oss is\dss equal\dss to $1$\nnsp.\oss
It\dss is\dss instructive\sss to deduce\sss this not\sss from explicit\sss computations
outlined\sss in\qss \cite{k},\oss but\sss from\sss our general\dss theory.\oss

\myuppar{The boundary\sss triplet.}
In order\sss to apply\sss the results of\trs Sections\dss \ref{abstract-index}\dss and\dss \ref{families},\oss
let\sss us\sss take\sss the\dss So\-bolev\dss space
$H_{\trf 2}\qff [\trf 0\fff,\qff 1\trf]$  as $H_{\dff 1}$\nsp,\oss
and $\ccc^{\dff 2}$
with\sss the standard\dss Hermitian\dss structure as $K^{\dff \partial},\qff K$\nnsp.\oss
Let\vspace{3pt}
\[
\quad
\gamma_0\dff(\trf u\trf)
\off =\off
(\qff u\trf(\dff 0\dff)\fff,\qff u\trf(\dff 1\dff)\qff)
\quad
\mbox{and}\quad
\gamma_1\dff(\trf u\trf)
\off =\off
(\qff u'\trf(\dff 0\dff)\fff,\qff -\qff u'\trf(\dff 1\dff)\qff)
\pff.
\]

\vspace{-12pt}\vspace{3pt}
The integration\sss by parts shows\sss that\sss the identity\qss (\ref{lagrange-identity})\qss holds.\oss
See\qss \cite{s},\oss Example\qss 14.2.\oss
Let\sss us\sss take as\sss the reference operator $A_{\trf \kappa}$\sss
for every value of\dss $\kappa$\sss 
the restriction\dss 
$A\off =\off T^{\dff *}\trf|\trf \kernel\fff \gamma_0$\nsp.\oss
As\dss is\dss well\dss known,\oss the operator $A$\sss is\dss self-adjoint\sss
and\sss invertible.\oss
Since\sss the\dss Hilbert\dss spaces\sss $K^{\dff \partial},\qff K$
are finitely dimensional\sss and equal,\oss
there\dss is\dss no\sss need\dss to pass\sss to\sss the reduced\dss boundary\sss triplet.\oss
In\sss fact,\oss already\trs Theorem\qss \ref{index-op}\qss implies\sss that\sss
the index of\dss the family\sss
$T\dff(\trf \kappa\trf)$\nnsp,\dss
$\kappa\qff \in\qff \rrr\trf \cup\dff \{\trf \infty \trf\}$\sss
is\dss equal\dss to\sss the index of\dss the family of\dss
the corresponding\sss boundary conditions.\oss
The computation of\dss the\sss latter\dss is\dss a finitely dimensional\dss problem.\oss
The relation\sss $\mathcal{R}\dff(\trf \kappa\trf)$ defining\sss $T\dff(\trf \kappa\trf)$\sss
is\vspace{3pt}
\[
\quad
\mathcal{R}\dff(\trf \kappa\trf)
\off =\off
\bigl\{\pff (\trf 0\fff,\qff b\fff,\qff c\fff,\qff -\qff \kappa\dff b\trf)
\qff \bigl|\qff
b\fff,\qff c\qff \in\qff \ccc
\off\bigr\}
\pff.
\]

\vspace{-12pt}\vspace{3pt}
It\dss is\dss equal\dss to\sss the direct\sss sum of\dss relations\sss
$0\dff \oplus\dff \ccc$\sss and\sss
$\mathcal{R}_{\dff 1}\dff(\trf \kappa\trf)
\off =\off
\bigl\{\pff (\trf b\fff,\qff -\qff \kappa\dff b\trf)
\qff \bigl|\qff
b\qff \in\qff \ccc
\off\bigr\}$\nnsp.\oss
It\sss follows\sss that\sss the index\dss is\dss equal\dss to\sss the index
of\dss the family\sss
$\mathcal{R}_{\dff 1}\dff(\trf \kappa\trf)$\nnsp,\dss
$\kappa\qff \in\qff \rrr\trf \cup\dff \{\trf \infty \trf\}$\nnsp.\oss
It\dss is\dss easy\sss to see\sss that\sss the\sss latter\sss
generates\sss the group\sss $K^{\dff 1}\dff(\trf S^{\fff 1}\trf)$\nnsp.\oss
Therefore\sss the index\dss is\dss equal\dss to $1$\nnsp,\oss
up\sss to\sss the choice of\dss the identification\sss
$K^{\dff 1}\dff(\trf S^{\fff 1}\trf)\off =\off \zzz$\nnsp.\oss
It\dss is\dss worth\sss to point\sss out\sss that\sss
$\mathcal{R}_{\dff 1}\dff(\trf \infty\trf)\off =\off 0\dff \oplus\dff \ccc$\sss
is\dss only a relation,\oss not\sss the graph of\dss an operator{}.\oss
In\sss finite dimension\sss the index of\dss any\sss family of\dss
self-adjoint\sss operators\dss is\dss equal\dss to $0$\sss because all\sss
such families are homotopic.\oss

\myuppar{The reduced\dss boundary\sss triplet.}
While\sss it\dss is\dss not\sss needed\sss for\sss the computation
of\dss the index,\oss the reduced\dss boundary\sss triplet\dss is\dss still\sss defined.\oss
The kernel\sss $\kernel\fff T^{\dff *}$\sss consists of\dss polynomials\sss of\dss degree $1$\nnsp.\oss
Clearly,\qss $\bm{\gamma}\trf(\dff 0\dff)$ maps\sss $(\trf a\fff,\qff b\trf)$\sss
to\sss the polynomial\sss 
$x\off \longmapsto\off (\trf b\qff -\qff a\trf)\trf x\qff +\qff a$\nnsp.\oss
Therefore\sss
$M\trf(\dff 0\dff)
\off =\off 
\Gamma_1\dff \circ\trf \bm{\gamma}\trf(\dff 0\dff)$\sss
maps\sss $(\trf a\fff,\qff b\trf)$\sss to\sss
$(\trf b\qff -\qff a\fff,\qff a\qff -\qff b\trf)$\sss
and\dss hence\sss the graph $\mathcal{M}$ of\dss $M\trf(\dff 0\dff)$\sss
is\dss equal\dss to\sss
$\bigl\{\pff (\trf a\fff,\qff b\fff,\qff b\qff -\qff a\fff,\qff a\qff -\qff b\trf)
\qff \bigl|\qff
a\fff,\qff b\qff \in\qff \ccc
\off\bigr\}$\nnsp.\oss
It\sss follows\sss that\vspace{3pt}
\[
\quad
\mathcal{R}\dff(\trf \kappa\trf)\qff -\qff \mathcal{M}
\off =\off
\bigl\{\pff (\trf 0\fff,\qff b\fff,\qff c\qff -\qff b\fff,\qff b\qff -\qff \kappa\dff b\trf)
\qff \bigl|\qff
b\fff,\qff c\qff \in\qff \ccc
\off\bigr\}
\pff.
\]

\vspace{-12pt}\vspace{3pt}
This relation\dss is\dss equal\dss to\sss the direct\sss sum of\dss relations\sss
$0\dff \oplus\dff \ccc$\sss and\sss
$\bigl\{\pff (\trf b\fff,\qff (\trf 1\qff -\qff \kappa\trf)\trf b\trf)
\qff \bigl|\qff
b\qff \in\qff \ccc
\off\bigr\}$\nnsp.\oss
Clearly,\oss the index of\dss this family of\dss relations\dss is\dss
the same as\sss the index of\dss the family\sss
$\mathcal{R}_{\dff 1}\dff(\trf \kappa\trf)$\nnsp,\dss
$\kappa\qff \in\qff \rrr\trf \cup\dff \{\trf \infty \trf\}$\nnsp.\oss
Not\sss surprisingly,\oss the reduced\dss boundary\sss triplet\dss
leads to\sss the same answer{}.\oss

%\newpage

\begin{flushright}
First\sss version\qss --\qss June\dss 16\fff,\oss 2023

Present\sss version\qss --\qss July\dss 20\fff,\oss 2023
 
https\halfff:/\!/\!nikolaivivanov.com

E-mail\halfff:\oss nikolai.v.ivanov{\fff}@{\dff}icloud.com,\oss ivanov{\fff}@{\dff}msu.edu

Department\sss of\qss Mathematics,\oss Michigan\sss State\sss University
\end{flushright}

\end{document}